\documentclass[11pt]{amsart}
\usepackage{amsthm,amsfonts,amssymb,amsmath,oldgerm}
\numberwithin{equation}{section}
\usepackage{fullpage}
\usepackage{amsmath}
\usepackage{setspace}
\usepackage{stmaryrd}
\usepackage{nicematrix}
\usepackage{mathrsfs}
\usepackage{fancyhdr}  
\usepackage{graphicx}
\usepackage{enumerate}
\usepackage{bm}
\usepackage{bbm}
\usepackage{dsfont}
\usepackage{multirow}
\usepackage{psfrag}
\usepackage[font=scriptsize]{caption}
\usepackage[font=scriptsize]{subcaption}
\usepackage{listings}
\usepackage{tikz}
\usepackage{empheq}
\usepackage{cases}
\usetikzlibrary{matrix} 
\usepackage[framemethod=TikZ]{mdframed}
\mdfdefinestyle{MyFrame}{backgroundcolor=gray!20!white}
\usepackage[makeroom]{cancel}

\usepackage[top=20mm,bottom=20mm,left=16mm,right=16mm,a4paper]{geometry}

\usepackage{hyperref,bookmark}

\usepackage[normalem]{ulem}
\definecolor{violet}{rgb}{0.580,0.,0.827}

\usepackage[textwidth= \marginparwidth,textsize=tiny]{todonotes}

\hypersetup{
    colorlinks,          
    citecolor=blue,        
    filecolor=blue,      
    urlcolor=blue,           
    linkcolor=blue,
}

\newcommand\dD{\mathrm{d}}

\def\eps{\varepsilon }

\newcommand{\J}{{\mathcal J}}

\newcommand{\beq}{\begin{equation}}
\newcommand{\eeq}{\end{equation}}
\newcommand{\beqa}{\begin{eqnarray}}
\newcommand{\eeqa}{\end{eqnarray}}
\newcommand\br{\begin{remark}}
\newcommand\er{\end{remark}}
\newcommand\bp{\begin{pmatrix}}
\newcommand\ep{\end{pmatrix}}
\newcommand{\be}{\begin{equation}}
\newcommand{\ee}{\end{equation}}
\newcommand\ba{\begin{equation}\begin{aligned}}
\newcommand\ea{\end{aligned}\end{equation}}
\newcommand\ds{\displaystyle}

\newcommand{\beg}{\begin{example}}
\newcommand{\eeg}{\end{exaplem}}
\newcommand{\bpr}{\begin{proposition}}
\newcommand{\epr}{\end{proposition}}
\newcommand{\bt}{\begin{theorem}}
\newcommand{\et}{\end{theorem}}
\newcommand{\bc}{\begin{corollary}}
\newcommand{\ec}{\end{corollary}}
\newcommand{\bl}{\begin{lemma}}
\newcommand{\el}{\end{lemma}}
\newcommand{\bd}{\begin{definition}}
\newcommand{\ed}{\end{definition}}
\newcommand{\brs}{\begin{remarks}}
\newcommand{\ers}{\end{remarks}}

\newtheorem{theorem}{Theorem}[section]
\newtheorem{proposition}[theorem]{Proposition}
\newtheorem{corollary}[theorem]{Corollary}
\newtheorem{lemma}[theorem]{Lemma}
\newtheorem{remark}[theorem]{Remark}
\newtheorem{definition}[theorem]{Definition}

\newtheorem{example}[theorem]{Example}

\newcommand{\N}{{\mathbb N}}

\newcommand{\R}{{\mathbb R}}
\newcommand{\T}{{\mathbb T}}

\newcommand\bL{{\mathbf L}}

\newcommand\bN{{\mathbf N}}

\newcommand\bR{{\mathbf R}}
\newcommand\bS{{\mathbf S}}

\newcommand\bU{{\mathbf U}}

\newcommand\cB{{\mathcal B}}
\newcommand\cC{{\mathcal C}}
\newcommand\cD{{\mathcal D}}
\newcommand\cE{{\mathcal E}}

\newcommand\cJ{{\mathcal J}}

\newcommand\cL{{\mathcal L}}
\newcommand\cM{{\mathcal M}}
\newcommand\cN{{\mathcal N}}

\usetikzlibrary{fit}
\usetikzlibrary{arrows}
\numberwithin{equation}{section}
\numberwithin{figure}{section}
\hypersetup{linkbordercolor=red}

\title[Numerical scheme for the quasineutral limit of the Vlasov-Poisson system]{A structure and asymptotic preserving scheme for the quasineutral limit of the Vlasov-Poisson system} 



\begin{document}

\maketitle

\centerline{\scshape Alain Blaustein\footnote{alain.blaustein@inria.fr, Centre Inria de l'Universit\'e de Lille, France.}, \scshape Giacomo Dimarco\footnote{giacomo.dimarco@unife.it, Universit\`a di Ferrara,
    Department of Mathematics and Computer Science \& Center for Modeling, Computing and Statistics (CMCS), Ferrara, Italy.},  \scshape Francis
  Filbet\footnote{francis.filbet@math.univ-toulouse.fr, Institut de
    Math\'ematiques de Toulouse, Universit\'e de Toulouse,
    France.} and \scshape Marie-H\'el\`ene
  Vignal\footnote{mhvignal@math.univ-toulouse.fr, Institut de
    Math\'ematiques de Toulouse, Universit\'e de Toulouse,
    France.}}
\medskip

\bigskip
 
\begin{abstract}
In this work, we propose a new numerical method for the Vlasov–Poisson system that is both asymptotically consistent and stable in the quasineutral regime, {\it i.e.} when the Debye length $\lambda$ is small compared to the characteristic spatial scale of the physical domain. Our approach consists in reformulating the Vlasov–Poisson system as a hyperbolic problem by applying a spectral expansion in the basis of Hermite functions in the velocity space and in designing a structure-preserving scheme for the time and spatial variables. Through this Hermite formulation, we establish a convergence result for the electric field toward its quasineutral limit together with explicit error estimates in $\lambda$. Following this path, we then propose a fully discrete numerical method for the Vlasov-Poisson system, inspired by the approach in \cite{abff24}, and rigorously prove that it is uniformly consistent in the quasineutral limit regime. Finally, we present several numerical simulations to illustrate the behavior of the proposed scheme. These results demonstrate the capability of our method to describe quasineutral plasmas, illustrate its stability and confirm its asymptotic preserving property.

\end{abstract}

\tableofcontents

\section{Introduction}
\label{sec1}
\setcounter{equation}{0}
\setcounter{figure}{0}
\setcounter{table}{0}

Plasma flows are multiscale problems, involving interactions ranging  from microscopic to macroscopic scales \cite{Chen}. As a result, the numerical simulation of plasmas presents a significant challenge for both the computational physics and the applied mathematics communities. Depending on the physical regime of interest, plasmas are typically described by two main classes of mathematical models: fluid and kinetic  models \cite{Chen, Degond}. Fluid models describe the evolution of macroscopic quantities, such as density, temperature, and mean velocity and remain valid when the plasma is close to thermodynamical equilibrium. In contrast, kinetic models detail the time evolution of a distribution function in six-dimensional phase space, representing the probability of a particle occupying a given state at any moment. While these models are able to describe a broader range of physical regimes, they are significantly more computationally demanding.

In this paper, we focus on plasmas far from equilibrium and therefore we adopt a kinetic perspective based on the Vlasov–Poisson system. Specifically, we  investigate  the challenges related to the  quasineutrality, which refers to the assumption that, on macroscopic scales, the net charge density in the plasma is effectively zero. Under this assumption, it is not possible to solve the electric field generated by local charge separations. The macroscopic behavior of the system can be described using asymptotic models, where the electric field is determined by additional constraints on the charge and current densities. However, in some scenarios, quasineutrality and charge imbalances can coexist. This makes challenging the development of numerical methods capable of accurately capturing both effects. Addressing this issue remains an active area of investigation \cite{ACS, BCDS:09, CDV:07, CDV:16, DDD, DDNSV:10, DLV}.

Here we combine a theoretical analysis of the behavior of the Vlasov-Poisson system in the quasineutral regime with the design of a new numerical method  capable of resolving the challenges related to its description. The physical model studied involves several interacting scales, among which the Debye length $\lambda_D$ and the electron plasma period $\omega_p^{-1}$ play a fundamental role:  
\be  
\label{debye}  
\lambda_D := \sqrt{\frac{\eps_0\, k_B\, T_0}{n_0\, e^2}},  
\qquad   \omega_p := \sqrt{\frac{n_0\, e^2}{\eps_0\, m}},  
\ee  
where $n_0$ and $T_0$ are the characteristic density and temperature, $e$ is the elementary charge, $\eps_0$ the vacuum permittivity, $k_B$ the Boltzmann constant, and $m$ the electron mass.  The Debye length defines the characteristic scale at which charge imbalances occur, while the electron plasma period represents the characteristic oscillation time associated with electrostatic forces that restore electric neutrality when charge imbalances arise at the Debye length scale \cite{DD:17}. When both the Debye length and the plasma period are small compared to the macroscopic scales of interest, the quasineutral regime holds true and the plasma appears broadly electrically neutral. In the sequel we introduce the details of the model employed.

\subsection{The mathematical model and review}
\label{sec11}
Let us first introduce the mathematical model and analyze the different scales involved. The Vlasov-Poisson system provides a kinetic description of a gas constituted of charged particles interacting through an electric field:
\begin{equation}
    \label{eq:vpfp3D}
  \left\{
    \begin{array}{l}
\ds\frac{\partial f}{\partial t}\,+\,v\cdot\nabla_x f
      \,+\, \frac{q}{m}\,E\cdot\nabla_v f \,=\,
      0\,,
      \\[1.1em]
      \ds E = -\nabla_x \phi\;\;;\;\;-\eps_0\Delta_{x}\phi \,=\, q\,(\rho- \rho_i)\;\;;\;\; \rho =\int_{\R^3} f \dD v \,.
      \end{array}\right.
  \end{equation}
In \eqref{eq:vpfp3D}, $f(t,x,v)$ is the distribution of electrons over the phase space $(x,v)\in\T^3\times\R^3$ at time $t\geq0$. Interactions are taken into account thanks to a coupling between kinetic and Poisson equations (first and second line of \eqref{eq:vpfp3D} respectively). The constant $q$ expresses the elementary charge $q=-e\,<\,0$, the scalar function $\phi$ is the electric potential, while the macroscopic densities of electrons and ions are denoted respectively by $\rho(t,x)$ and $\rho_i(t,x)=\rho_i(x)$. Let us note that in the system \eqref{eq:vpfp3D},  ions constitute a fixed background, which simplifies the dynamics without omitting the fundamental aspects of the quasineutrality mechanism.

We now introduce  the rescaled parameter $\lambda$  given by
  $$
\lambda \,=\,
\frac{\lambda_D}{L}
\,=\, \left( \frac{\eps_0\,k_B T_0}{L^2 e^2\,n_0}\right)^{1/2},
$$
where $L$ represents the characteristic length at which the physical phenomena are studied, in our case it corresponds to the length of the electric field interactions in plasmas. The quasineutral regime corresponds to the limit of the Vlasov-Poisson system when $\lambda\rightarrow 0$. This asymptotic limit has been extensively studied in the past by several authors from both the numerical and theoretical points of view. For a theoretical analysis in the context of kinetic and fluid models, see for example \cite{B:00, Grenier96, VHR:14, Han-Kwan_Hauray15, HK16,  HKI, Bobylev_Potapenko19, BP, GPI} and the references therein. The study of the quasineutral limit is often based on reformulation of the Vlasov-Poisson system with respect to the dimensionless parameter $\lambda$ and then in passing to the limit $\lambda$ going to zero in  this rescaled system. Here, we restrict ourselves to the rescaled one dimensional in space and in velocity setting. This leads to the following equations, see \cite{CDV:16, DD:17} for details, 
\begin{equation}
	\label{vpfp:0}
	\left\{
	\begin{array}{l}
		\ds\partial_t f^\lambda \,+\, v\,\partial_x f^\lambda
		\,+\,E^\lambda \,\partial_v f^\lambda \,=\,0\,,
		\\[1.1em]
		\ds E^\lambda\,=\, -\partial_{x}\phi^\lambda\;;\;-\lambda^2\partial_{x}^2\phi^\lambda \,=\,\rho^\lambda-1\;;\;\rho^\lambda (t,x)\,=\,\int_{\R}\,f^\lambda(t,x,v)\,\dD v\,,
	\end{array}\right.
\end{equation}
where $f^\lambda(t,x,v)$ is the distribution at time $t\geq 0$ over the phase-space $(x,v)\in\T\times \R$ and where the density of the ions is assumed to be constant and equal to one everywhere in space. 
The following condition ensures the uniqueness of $\phi^\lambda$
\begin{equation}
	\label{vpfp:1}
	\int_\T \phi^\lambda(t,x) \,\dD x\,=\; 0\,.
\end{equation}
The model \eqref{vpfp:0}-\eqref{vpfp:1} possesses the following key properties: conservation of {total} charge, of total flux (set to $0$ for simplicity and without loss of generality) and of the global energy. These read
\begin{equation}
	\label{conservation:flux:mass}
	\left\{
	\begin{array}{l}
		\ds \int_{\T}\,\rho^\lambda(t,x)\,\dD x \,=\,  \int_{\T}\,\rho^\lambda(0,x)\,\dD x\,=\, \left|\T\right|\,,\\[1.1em]
		\ds \int_{\T}\,j^\lambda(t,x)\,\dD x\,=\,  \int_{\T}\,j^\lambda(0,x)\,\dD x\,=\,0\,,\\[1.1em]
		\ds \frac{1}{2}\int_{\T}K^\lambda(t,x)+\lambda^2\left|E^\lambda(t,x)\right|^2\dD x= \frac{1}{2}\int_{\T}K^\lambda(0,x)+\lambda^2\left|E^\lambda(0,x)\right|^2\dD x\,,
	\end{array}\right.
\end{equation}
where the flux $j^\lambda$, also called the current density,  and the kinetic energy $\frac12 K^\lambda$ are given by
\begin{equation}
	\label{flux:kin}
	\left\{
	\begin{array}{l}
		\ds j^\lambda(t,x)\,=\,
		\int_{\R}\,vf^\lambda(t,x,v)\,\dD v \,,\\[1.1em]
		\ds K^\lambda(t,x)=
		\int_{\R}|v|^2f^\lambda(t,x,v)\,\dD v\,.
	\end{array}\right.
\end{equation}

From a numerical perspective, traditional explicit schemes applied to \eqref{vpfp:0} must resolve the microscopic scales related to the small parameter $\lambda$ to ensure stability and consistency in the limit $\lambda \to 0$. However, this requirement necessitates extremely small time steps and phase space discretizations in order to remain stable. At the same time, numerical simulations must be performed down to the macroscopic scale to capture meaningful phenomena. This situation creates significant computational challenges and makes these methods impractical for realistic applications. While asymptotic models can be derived to describe macroscopic regimes, additional challenges arise in scenarios where quasineutral and non-quasineutral regions coexist, making simulations particularly difficult.

To address this complexity, domain decomposition approaches or hybrid
methods can be employed \cite{Dim2008, Dim2014, FR:15, FX:18}. However, integrating different models and numerical methods requires to accurately identify interfaces, which remains a challenging task. Therefore, it is crucial to design numerical methods capable of handling multiple regimes simultaneously, without being constrained by small-scale dynamics. The development of schemes that work without such limitations has been the focus of
extensive research in the recent years. This is precisely the area in which Asymptotic Preserving (AP) methods have been created, see \cite{jin1, jin2} and \cite{jin3} for a recent review on the subject. These methods can bypass the above discussed restrictions and automatically achieve transition using consistent discretizations of limiting models as the parameters characterizing microscopic behavior approach zero. The use of Asymptotic Preserving schemes in the context of quasineutrality has been investigated for instance in \cite{CDV:16}. The authors proposed a first-order-in-time scheme and a linear stability analysis proving indeed stability for small values of the Debye length of the proposed method. The concept revolved around the reformulation and the consequent discretization of the Poisson equation. This idea was initially introduced in \cite{CDV:07} and later used
in \cite{BCDS:09} and \cite{ACS,DDNSV:10, DLV} or more recently for the Vlasov-Maxwell system \cite{DDD, JinAndCo, Lapenta}. For a general review of numerical methods capturing the quasineutral limit for both kinetic and fluid models, we refer to \cite{DD:17}. In this work, we adopt an alternative approach that we develop in two parts. First, we study relevant analytical properties of the Vlasov–Poisson system close to the quasineutral limit. Then, we propose a new numerical method based on this theoretical result and  perform a discrete counterpart analysis. This allows us to highlight the main properties of the scheme such as stability and preservation of the asymptotic state. Unlike previous approaches in the literature \cite{CDV:07, CDV:16, DD:17}, our discretization does not rely on reformulating the Poisson equation as an equivalent harmonic oscillator equation. We instead discretize the Poisson equation concurrently with the equations governing the time evolution of the first moments of the distribution function. \\

In the following section, we first recall the formal limit of the Vlasov–Poisson system as the Debye length approaches zero (Section \ref{sec12}). Then, in Section \ref{sec2}, we reformulate the Vlasov–Poisson system as a hyperbolic system with source terms by expanding the distribution function in velocity using a basis of Hermite functions. This reformulation is particularly well-suited for studying the quasineutral limit, as we will demonstrate. Next, we establish a convergence result for the electric field toward its quasineutral limit by explicitly characterizing the oscillatory component of the solution (see Proposition \ref{Prop:dev:asymptotique}) and highlighting the impact of the initial conditions on the size of these oscillations.
Starting from this formulation, we propose a numerical scheme based on a time-splitting strategy in Section \ref{sec3}. We rigorously prove that the scheme captures the asymptotic behavior of the solution in Section \ref{disc:QN}. In Section \ref{sec4}, we present several numerical simulations to illustrate the behavior of the proposed discretization in relation to the theoretical results of Proposition \ref{Prop:dev:asymptotique}, including an analysis of the error estimates. 
Furthermore, a key aspect of our numerical investigation focuses on scenarios where our analysis does not hold. By studying these cases, we aim to identify the necessary precautions to avoid physically irrelevant results.

\subsection{Formal quasineutral limit}
\label{sec12}
Let us now observe that in the regime $\lambda\ll 1$, the second line in \eqref{vpfp:0} forces the {global} neutrality condition expressed by the first line of \eqref{conservation:flux:mass} to become a {pointwise} constraint in $x$. Thus, one expects the solution $f^\lambda$ to formally converge to the solution $f$ of the quasineutral system
\begin{equation}
	\label{vpfp:neutral}
	\left\{
	\begin{array}{l}
		\ds\partial_t f\,+\, v\,\partial_x f
		\,+\,E\,\partial_v f\,=\,0\,,
		\\[1.1em]
		\ds \rho(t,x)\,=\,1\;;\;\rho (t,x)\,=\,\int_{\R}\,f(t,x,v)\,\dD v\,,
	\end{array}\right.
\end{equation}
where $f(t,x,v)$ is the limiting distribution over the phase-space $(x,v)\in\T\times \R$. The latter system may be confusing at first glance, as the equation on the limiting field $E$ seems to be missing. In fact, $E$ remains fully determined even in \eqref{vpfp:neutral}, as it may be interpreted as the Lagrange multiplier associated to the pointwise neutrality constraint from the second line of equations \eqref{vpfp:neutral}. Indeed, the quasineutral limit described by system \eqref{vpfp:neutral} together with the conditions \eqref{conservation:flux:mass} imposes the limiting electric field and the limiting electric flux to satisfy the following relations
\begin{equation}
	\label{neutral:field:flux}
	\left\{
	\begin{array}{l}
		\ds j(t,x)\,=\,0\,,
		\\[1.1em]
		\ds E(t,x)\,=\,\partial_x K(t,x)\,,
	\end{array}\right.
\end{equation}
where the kinetic energy $\frac 12 K$ is given by \eqref{flux:kin}. Hence, global conservation of the total charge and the total flux in \eqref{conservation:flux:mass} both become pointwise constraints in the quasineutral regime, according to the second line of \eqref{vpfp:neutral} and the first line of \eqref{neutral:field:flux} respectively.
\begin{remark}
Let us observe that \eqref{neutral:field:flux} is obtained  multiplying the Vlasov equation  \eqref{vpfp:neutral} by $(1,v)^t$ and integrating over $v\in \R$. We then use the pointwise neutrality constraint and total conservation of flux for $j(t,x)$, and both pointwise neutrality and conservation of flux to obtain the expression of the electric field $E(t,x)$.
\end{remark}
Before continuing, we point out that the quasineutral system may not be globally well-posed for all initial data. This was observed in \cite{Bardos14} for the related Vlasov-Dirac-Benney system and in \cite{Han-Kwan_Hauray15} near the so-called Penrose unstable profiles for the Vlasov-Poisson system. Indeed, in such cases, small perturbations from the quasineutral state can lead to the development of instabilities in the plasma. In what follows, we restrict our analysis to configurations near Penrose stable profiles. However, in the numerical section, we will simulate both stable and unstable cases to illustrate these phenomena.\\

One of the key challenges in the numerical analysis and simulation of the quasineutral regime is that the strong convergence of  $f^\lambda$ to $f$ may fail due to the emergence of fast oscillations with frequency $1/\lambda$, that is, period $2\pi\lambda$, when $\lambda \ll 1$, see for instance \cite{Grenier96, Han-Kwan_Hauray15, HK16, Bobylev_Potapenko19}.  These oscillations result from the fact that the pointwise neutrality and conservation of flux in \eqref{vpfp:neutral}-\eqref{neutral:field:flux} may not be verified by the non-neutral initial distribution $f^\lambda(t=0)$. To uncover these oscillations, we determine their amplitude thanks to formal energy considerations. In particular, we focus on initial configurations with uniformly bounded total energy as $\lambda \ll 1$, that is,
\[
\sup_{\lambda > 0}\left(
\int_{\T}K^\lambda(0,x)+\lambda^2\left|E^\lambda(0,x)\right|^2\dD x\right)
\,<\,+\infty\,.
\]
Since the total energy is conserved according to the third line of \eqref{conservation:flux:mass}, this yields for all $t\geq 0$
\[
\sup_{\lambda > 0}\left(
\int_{\T}K^\lambda(t,x)+\lambda^2\left|E^\lambda(t,x)\right|^2\dD x\right)
\,<\,+\infty\,.
\]
In particular, this estimate indicates that $E^\lambda$ is at most of order $O(\lambda^{-1})$:
\[
\sup_{\lambda>0}\left(
\lambda\left\|E^\lambda(t)\right\|_{L^2\left(\T\right)}\right)
\,<\,+\infty\,;\quad\textrm{that is,}\quad
\left|E^\lambda(t)\right|\lesssim\,\frac{1}{\lambda}\,,\quad \textrm{as}\quad \lambda \ll 1\,,
\]
which means that the amplitude of the electric field oscillations is at most of order $O\left(\lambda^{-1}\right)$. In the following, we restrict our investigation to cases characterized by smaller amplitudes of the electric field, namely
\begin{equation}\label{hyp:amplitude}
\left|E^\lambda(t)\right|\lesssim\,\frac{1}{\lambda^\alpha}\,,\quad\textrm{for some}\quad 0\,\leq\,\alpha\,<\,1\,,\quad \textrm{as}\quad \lambda \ll 1.
\end{equation}
It is worth mentioning that the analysis of the critical case $\alpha\,=\,1$ significantly differs from the cases where $\alpha < 1$ \cite{Grenier96,Han-Kwan_Hauray15}. We will discuss this critical case in the section devoted to numerical simulations. Now, to identify the period of oscillations, we analyze the coupled system current density $j^\lambda$ and electric field $E^\lambda$, which reads
\begin{equation}\label{amp:eq}
		\left\{
	\begin{array}{l}
		\ds\lambda^2 \partial_t \,E^\lambda \,+\,j^\lambda \,=\,0\,,
		\\[1.1em]
		\ds \partial_t \,j^\lambda \,+\, \partial_x K^\lambda
		\,- \rho^\lambda\,E^\lambda  \,=\,0\,.
	\end{array}\right.
\end{equation}
This is obtained multiplying the Vlasov equation in \eqref{vpfp:0} by $(1,v)^t$, integrating over $v\in \R$ and then substituting the Poisson equation in the continuity equation. We now rewrite the coupled system \eqref{amp:eq} as
\begin{equation}
	\label{vect:Amp:flux}
	\partial_t
	\begin{pmatrix}
		\ds E^\lambda - \partial_x K^\lambda\\[0.7em]
		\ds \lambda^{-1}\,j^\lambda
	\end{pmatrix}\,+\,
	\frac{1}{\lambda}\,J\,\cdot
	\begin{pmatrix}
		\ds E^\lambda - \partial_x K^\lambda\\[0.7em]
		\ds \lambda^{-1}\,j^\lambda
	\end{pmatrix}
	\,+\,
	\bR^\lambda\,=\,0\,,
\end{equation}
where the matrix $J$ and the vector $\bR^\lambda$ are given as follows
\begin{equation*}
	J\,=\,
	\begin{pmatrix}
		\ds 0 & \;1\; \\[0.5em]
		\ds -1 & \;0
	\end{pmatrix}\,,\quad \textrm{and}\quad
	\bR^\lambda\,=\,
	\begin{pmatrix}
		\ds\,\partial_{tx} K^\lambda\\[1.em]
		\ds  \lambda^{-1}\,\left(1-\rho^\lambda\right)E^\lambda 
	\end{pmatrix}\,.
\end{equation*}
Equation \eqref{vect:Amp:flux} shows that the electric field and the flux behave like harmonic oscillators with frequency $1/\lambda$ up to the source term $\bR^\lambda$. Indeed, multiplying \eqref{vect:Amp:flux} by the matrix $\exp{(tJ/\lambda)}$ and integrating in time, we obtain the following Duhamel formula
\[
\begin{pmatrix}
	\ds E^\lambda - \partial_x K^\lambda\\[0.7em]
	\ds \lambda^{-1}\,j^\lambda
\end{pmatrix}(t)\,=\,
\exp{\left(-\,\frac{t}{\lambda}\,J\right)}\cdot\,\begin{pmatrix}
	\ds E^\lambda - \partial_x K^\lambda\\[0.7em]
	\ds \lambda^{-1}\,j^\lambda
\end{pmatrix}(0)\,\,-\,\,
\int_{0}^{t}
\exp{\left(\frac{s-t}{\lambda}\,J\right)}\cdot\,\bR^\lambda(s)\,\dD s\,,
\]
for all time $t\geq 0$, 
where $\exp{(tJ)}$ is the rotation matrix with angle $t$:
\begin{equation}\label{mat:rot}
	\exp{\left(tJ\right)}\,=\,
	\begin{pmatrix}
		\ds \;\;\;\cos\left(t\right) & \ds\sin\left(t\right) \\[1.em]
		\ds -\sin\left(t\right) & \ds\cos\left(t\right)
	\end{pmatrix}\,.
\end{equation}
The difficulty is then to prove that, under the conditions \eqref{hyp:amplitude},
the remainder in the Duhamel formula may be neglected in the regime
$\lambda \ll 1$ (see also \cite{Grenier96, Han-Kwan_Hauray15,Bobylev_Potapenko19, HK16}), which yields
\begin{equation}\label{expansion:field:flux}
	\left\{
	\begin{array}{lll}
		\ds E^\lambda(t)\,&\ds\underset{\lambda \ll 1}{\simeq}\, \partial_x K^\lambda(t)
		&\ds+\,\cos\left(\frac{t}{\lambda}\right)\left(E^\lambda - \partial_x K^\lambda\right)(0)
		\,-\,\frac{1}{\lambda}
		\sin\left(\frac{t}{\lambda}\right)\,j^\lambda (0)
		\,,
		\\[1.5em]
		\ds j^\lambda(t)\,&\ds\underset{\lambda \ll 1}{\simeq}\, \quad\;\;0
		&\ds+\,\lambda\sin\left(\frac{t}{\lambda}\right)\left(E^\lambda - \partial_x K^\lambda\right)(0)
		\,+\,
		\cos\left(\frac{t}{\lambda}\right)\,j^\lambda (0)
		\,.
	\end{array}\right.
\end{equation}
These formal computations show that the quasineutral limit \eqref{neutral:field:flux} is satisfied up to fast oscillations of frequency $1/\lambda$ with amplitude only depending on the initial data $f^\lambda(t=0)$. 

To prepare the design and analysis of our numerical method, let us now make this analysis rigorous in the Hermite framework.

\section{Hermite formulation of the quasineutral regime}
\label{sec2}
\setcounter{equation}{0}
\setcounter{figure}{0}
\setcounter{table}{0} 
The Hermite decomposition arises naturally when studying the Vlasov-Poisson system. It corresponds to a moment method tailored with an additional spectral structure and well suited for accurate numerical approximations. We start by recasting \eqref{vpfp:0} and \eqref{vpfp:neutral} in the Hermite framework, successively we highlight how quasineutral oscillations naturally uncover in this setting. 
\subsection{The Hermite framework}\label{sec21}
We decompose the distribution $f^\lambda$ and its quasineutral counterpart $f$ into their Hermite coefficients $
C^\lambda\,=\,
\left(
C^\lambda_{k}
\right)_{k\in\N}
$ and $
C\,=\,
\left(
C_{k}
\right)_{k\in\N}
$, that is,
\begin{equation}\label{f:decomp}
	f^\lambda\left(t,x,v\right)
	\,=\,
	\sum_{k\in\N}\,
	C^\lambda_{k}
	\left(t,x
	\right)\,\Psi_{k}(v)
	\quad\quad\textrm{and}\quad\quad
	f\left(t,x,v\right)
	\,=\,
	\sum_{k\in\N}\,
	C_{k}
	\left(t,x
	\right)\,\Psi_{k}(v)
	\,,
\end{equation}
where the basis of Hermite functions $\left(\Psi_k\right)_{k\in\N}$ is defined recursively as follows
$\Psi_{-1}=0$, $\Psi_{0}=\cM$ and
\begin{equation}\label{recursive}
	v\,\Psi_{k}(v)\,=\,
	\sqrt{T_0 k}\,\Psi_{k-1}(v)
	\,+\,
	\sqrt{T_0 (k+1)}\,\Psi_{k+1}(v)
	\,,\quad\forall\, k\,\geq\,0\,,
\end{equation}
and where $\cM$ denotes the stationary Maxwellian with fixed temperature $T_0>0$
\begin{equation}
	\label{M:T0}
	\cM(v)\,=\, \frac{1}{\sqrt{2\pi \,T_0}} \, \exp\left( -\frac{|v|^2}{2\,T_0}\right)\,.
\end{equation}
Hermite functions $\left(\Psi_k\right)_{k\in\N}$ constitute an orthonormal system for the inverse Gaussian weight:
\begin{equation}\label{orth}
	\int_{\R}\,
	\Psi_{k}(v)\,\Psi_{l}(v)\,\cM^{-1}(v)\dD v
	\,=\,
	\delta_{k,l}\,,
\end{equation}
where $\delta_{k,l}$ denotes the Kronecker symbol ($\delta_{k,l}=1$ when $k=l$ and
$\delta_{k,l}=0$ otherwise). The decomposition \eqref{f:decomp} may also be interpreted as a moment method, where $C_k$, $k\geq 0$, stands for the moment of the distribution of order $k$, appropriately modified in order to satisfy the orthogonality constraint \eqref{orth}. In particular, from \eqref{f:decomp} one can easily retrieve the macroscopic quantities associated to $f^\lambda$. For example, mass \eqref{vpfp:0}, current density and kinetic energy \eqref{flux:kin} are given by
\[
\rho^\lambda(t,x)
\,=\, C^\lambda_0(t,x)\,;\quad
j^\lambda(t,x)
\,=\, \sqrt{T_0}\,C^\lambda_1(t,x)\,;\quad
K^\lambda(t,x)
\,=\, T_0\left(\sqrt{2}\,\,C^\lambda_2+\,C^\lambda_0\right)(t,x)\,.
\]
Thanks to the above decomposition, one can rewrite \eqref{vpfp:0} and \eqref{vpfp:neutral} as an infinite hyperbolic system with unknowns $C^\lambda\,=\,
\left(
C^\lambda_{k}
\right)_{k\in\N}
$ and $
C\,=\,
\left(
C_{k}
\right)_{k\in\N}
$, where the {discrete} spectral parameter $k \in \N$ now replaces the velocity variable $v\in\R$. More precisely, for each $k\geq 0$, we compute the equation for $C_k^\lambda$ by multiplying \eqref{vpfp:0} with $\Psi_k \,\cM^{-1}$ and integrating over $v\in \R$. We then use \eqref{recursive} to compute the contribution of the free transport operator $v\,\partial_x $ and the following relation to compute the field contribution $\partial_x \phi^\lambda \partial_v $ 
\begin{equation*}
	\partial_v \left(\Psi_k\,\cM^{-1}\right)\,=\,\sqrt{\frac{k}{T_0}}\,\Psi_{k-1}\,\cM^{-1}\,,\quad\forall\, k\,\geq\,0\,.
\end{equation*}
Finally, the orthogonality property \eqref{orth} is used to close the system of equations. This yields to the following coupled system for the coefficients $C^\lambda\,=\,(C^\lambda_{k})_{k \in \N}$,
\begin{equation}
	\label{Hermite:D}
	\left\{
	\begin{array}{l}
		\ds\partial_t C^\lambda_{k} \,+\,
		\sqrt{T_0\,k}\,
		\partial_x C^\lambda_{k-1}\,+ \, \sqrt{T_0\,(k +1)}\,
		\partial_x C^\lambda_{k+1}\,-\,\sqrt{\frac{k}{T_0}}\, E^\lambda\,C^\lambda_{k-1}
		\,  =\,  0\,, \qquad \forall\,k\in\N\,,
		\\[1.2em]
		\ds E^\lambda\,=\, -\partial_{x}\phi^\lambda\,;\;-\lambda^2\partial_{x}^2\phi^\lambda \,=\,C_0^\lambda-1\,,
	\end{array}\right.
\end{equation}
where we have set $C^\lambda_{-1}=0$, whereas the Hermite coefficients $C\,=\,(C_{k})_{k \in \N}$ related to the asymptotic distribution $f(t,x,v)$ satisfy
\begin{equation}
	\label{Hermite:D:neutral}
	\left\{
	\begin{array}{l}
		\ds C_0\,=\, 1\,,
		\\[0.5em]
		\ds C_1\,=\, 0\,,
		\\[0.1em]
		\ds\partial_t C_{k} \,+\,
		\sqrt{T_0\,k}\,
		\partial_x C_{k-1}\,+ \, \sqrt{T_0\,(k +1)}\,
		\partial_x C_{k+1}\,-\,\sqrt{\frac{k}{T_0}}\, E\,C_{k-1}
		\,  =\,  0\,, \qquad \forall\,k\geq 2\,,
		\\[0.7em]
		\ds E\,=\, \sqrt{2}\,T_0\,\partial_x C_2\,.
	\end{array}\right.
\end{equation}
In the next section, using the Hermite framework introduced above, we
justify the formal computations presented in Section \ref{sec12}.
\subsection{Continuous analysis of the quasineutral regime}
\label{sec22}
To compute the oscillatory components of $(f^\lambda,\phi^\lambda)$ within the Hermite framework, we first rewrite the Ampère and flux equations \eqref{amp:eq} as
\begin{equation}\label{Amp:flux:H:0}
\lambda^2 \partial_t E^\lambda \,+\,\sqrt{T_0}\,C_1^\lambda \,=\,0\,,\quad \textrm{and}\quad 
\partial_t C^\lambda_{1} \,+\,
\sqrt{T_0}\,
\partial_x C^\lambda_{0}\,+ \, \sqrt{2\,T_0}\,
\partial_x C^\lambda_{2}\,-\,\frac{1}{\sqrt{T_0}}\, E^\lambda\,C_{0}^\lambda
\,  =\,  0\,.
\end{equation}
Therefore, the system \eqref{vect:Amp:flux} takes the form
\begin{equation}\label{eq:sysh}
	\partial_t
	\begin{pmatrix}
		\ds E^\lambda - E_{\rm slow}^\lambda \\[0.5em]
		\ds \lambda^{-1}\sqrt{T_0}\,C_1^\lambda 
	\end{pmatrix}\,+\,
	\frac{1}{\lambda}\,J\cdot
	\begin{pmatrix}
		\ds E^\lambda - E_{\rm slow}^\lambda \\[0.5em]
		\ds \lambda^{-1}\sqrt{T_0}\,C_1^\lambda 
	\end{pmatrix}
	\,+\,
	\bS^\lambda\,=\,0\,,
\end{equation}
where $J$ is defined in \eqref{vect:Amp:flux} and $\bS^\lambda$ is given by
\begin{equation*}
	\bS^\lambda\,=\,
	\begin{pmatrix}
		\ds\sqrt{2}\,T_0\,\partial_t 
		\partial_x C^\lambda_{2} \\[0.5em]
		\ds \lambda^{-1}(T_0\,
		\partial_x C^\lambda_{0}\,+\,E^\lambda\,(1-C^\lambda_{0}))
	\end{pmatrix}.
\end{equation*}
The quantity $E_{\rm slow}^\lambda$ plays the role of the quasineutral electric field and it is defined by
$$
E_{\rm slow}^\lambda := \sqrt2 \,T_0 \partial_x C_2^\lambda\,.
$$
In the introduction, we interpreted the system \eqref{eq:sysh} as a coupled harmonic oscillator up to a source term $\bS^\lambda$. To justify this statement, we now analyze the behavior of $E^\lambda - E_{\rm slow}^\lambda$ and $C_1^\lambda$ in time. To this aim, we  define $\bU^\lambda$ as
\[
\bU^\lambda(t)\,:=\,
\begin{pmatrix}
	\ds E^\lambda - E_{\rm slow}^\lambda \\[0.5em]
	\ds \lambda^{-1}\sqrt{T_0}\,C_1^\lambda 
\end{pmatrix}(t)
\,-\,
\exp{\left(-\frac{t}{\lambda}J\right)}\cdot\,
\begin{pmatrix}
	\ds E^\lambda - E_{\rm slow}^\lambda  \\[0.5em]
	\ds \lambda^{-1}\sqrt{T_0}\,C_1^\lambda 
\end{pmatrix}(t=0)\,,
\]
where $\exp{(tJ)}$ is the rotation matrix given by \eqref{mat:rot}. The vector $\bU^\lambda$ vanishes at $t=0$ while it satisfies the following equation
\begin{equation}\label{eq:U}
	\partial_t \bU^\lambda
	\,+\,
	\frac{1}{\lambda}\,J\cdot
	\bU^\lambda
	\,+\,
	\bS^\lambda\,=\,0\,.
\end{equation}
Multiplying \eqref{eq:U} by $\exp{(tJ/\lambda)}$ and integrating in time, we obtain the following Duhamel formula
\begin{equation}\label{duhamel}
\bU^\lambda(t)\,=\,\,-\,\,
\int_{0}^{t}
\exp{\left(\frac{s-t}{\lambda}J\right)}\cdot\,\bS^\lambda(s)\,\dD s\,, \qquad \forall t\geq 0\,.
\end{equation}
In Proposition \ref{Prop:dev:asymptotique} below, we prove that if \eqref{hyp:amplitude} is satisfied, the remainder in the former Duhamel formula may be neglected as $\lambda \rightarrow 0$, which leads to the asymptotic expansion \eqref{expansion:field:flux}. This result ensures that the field $ E^\lambda$ and the flux $C_1^\lambda$ converge to  the quasineutral limit up to fast oscillations in time.

To make the formal analysis of Section \ref{sec12} rigorous, we now state precisely, in the Hermite framework, the assumptions on the solution to the Vlasov-Poisson system \eqref{vpfp:0} which are needed in our analysis. The main one is the rigorous counterpart of the formal condition \eqref{hyp:amplitude}: we assume that there exists $0\,\leq\,\alpha\,<\,1$ and a final time $T>0$ such that
\begin{equation}\label{compatibility}
	\sup_{0<\lambda<1}
	\left(
		\lambda^{\alpha-1}
	\left\|
	C^{\lambda}_1(0)
	\right\|_{W^{r_0+1,4}\left(\T\right)}
	+
	\sup_{t\in[0,T]}
	\lambda^{\alpha}\left\| E^\lambda(t) \right\|_{W^{r_0+2,4}\left(\T\right)}\right)
	\;<\;+\infty\,,
\end{equation}
where $r_0 = \lceil 1/(1-\alpha)\rceil $, $\lceil \cdot\rceil$ denotes the upper integer part and where $W^{r_0,4}$ is the Sobolev space of order $r_0$ based on $L^4$.
Let us comment on assumption \eqref{compatibility}. The bound on $E^\lambda$ is exactly the formal condition \eqref{hyp:amplitude}, now measured in a Sobolev norm, and it is required uniformly on a time interval $[0,T]$ which does not depend on $\lambda$. The bound on the initial flux $C_1^\lambda(0)=j^\lambda(0)/\sqrt{T_0}$ is a consequence of the latter one. Indeed, according to the formal expansion \eqref{expansion:field:flux}, the oscillating part of the electric field contains the term $\lambda^{-1}\sin(t/\lambda)\,j^\lambda(0)$, hence the condition $E^\lambda=O(\lambda^{-\alpha})$ imposes $j^\lambda(0)=O(\lambda^{1-\alpha})$. Let us also emphasize that such bounds, on a time interval independent of $\lambda$, have been justified in a very close setting in \cite[Theorem 1.1.2]{Grenier96}.

Furthermore, we assume that 
\begin{equation}
  \label{uniform:reg}
\sup_{0<\lambda<1}\sup_{t\in[0,T]}\max_{k\leq 4}
\left\| C_k^\lambda(t) \right\|_{W^{r_0+3,4}\left(\T\right)}
\;<\;+\infty\,,
\end{equation}
meaning that the solution to the Vlasov-Poisson system remains uniformly regular with respect to $\lambda$.
Under the latter two hypotheses, we are able to prove that the formal expansion \eqref{expansion:field:flux} holds true with explicit convergence rates in $\lambda$, whose sharpness is illustrated by the numerical simulations of Section \ref{sec:4.1}. 
To that aim, we define the following quantities
 \begin{equation}
   \label{cE}
 \left\{
   \begin{array}{l}
     \ds\cE_0^\lambda(t) \,:=\, \left\|E^\lambda(t)-  E^\lambda_{\rm slow}(t)
		-\,E^\lambda_{\textrm{osc}}(t)\right\|_{L^2(\T)}\,, 
     \\[1.1em]
    \ds\cE_1^\lambda(t) \,:=\,  \left\|C_1^\lambda(t)\,-\,C_{1,\textrm{osc}}^\lambda(t)\right\|_{L^2(\T)}\,, 
   \end{array}\right.
 \end{equation}
 where we decompose the electric field into two parts, namely the quasineutral state and its oscillating counterpart, and the same for the current density $C_1^\lambda$. These oscillating quantities read
 \begin{equation}\label{E:C_1:osc}
	\left\{
	\begin{array}{llll}
		\ds E^\lambda_{\textrm{osc}}(t,x) &\ds:=\,\cos\left(\frac{t}{\lambda}\right)\left(E^\lambda - E_{\rm slow}^\lambda\right)(0,x)
		\;-\;\frac{\sqrt{T_0}}{\lambda}
		\sin\left(\frac{t}{\lambda}\right)\,C_1^\lambda (0,x) \,,
		\\[1.3em]
		\ds C_{1,\textrm{osc}}^\lambda(t,x) &\ds:=\,\cos\left(\frac{t}{\lambda}\right)\,C_1^\lambda (0,x)
		\;+\;
		\frac{\lambda}{\sqrt{T_0}}\sin\left(\frac{t}{\lambda}\right)\left(E^\lambda - E_{\rm slow}^\lambda\right)(0,x)
		 \,.
	\end{array}\right.
    \end{equation}
We are now able to prove the following result, which is a quantitative version of the formal expansion \eqref{expansion:field:flux} written in the Hermite framework.
\begin{proposition}\label{Prop:dev:asymptotique}
  Consider a family of solutions $(C^\lambda,\phi^\lambda)_{\lambda>0}$ to \eqref{Hermite:D} with zero total flux
  $$
  \int_{\T}\,C_1^\lambda(t,x)\,\dD x\,=0\,,
  $$
  and global mass $|\T|$. Suppose that they satisfy the compatibility assumption \eqref{compatibility} and the uniform regularity assumption \eqref{uniform:reg}, for some given final time $T>0$. Then, for all $t\in[0,T]$, and all $0<\lambda<1$, the following bounds hold for the electric field and the flux $(E^\lambda, C_1^\lambda)$:
\begin{equation*}
		 \cE_0^\lambda(t)\,\leq\,
		C\,\lambda^{1-\alpha} 
		\qquad {\rm and }\qquad 
		\cE_1^\lambda(t) \,\leq\,
		C\,
		\lambda^{2-\alpha}\,,
\end{equation*}
where the constant $C$ only depends on $\alpha$, $T$, $T_0$, $|\T|$ and the implicit constants in \eqref{compatibility}-\eqref{uniform:reg}.
\end{proposition}
\begin{proof}
The proof is postponed in Appendix \ref{App:prop:dev:asympt}.
\end{proof}

  \begin{remark}
  \label{rem:original:variables}
Each coefficient $C_k^\lambda$ is a linear combination of the moments in $v$ of $f^\lambda$ of order at most $k$. The proof of Proposition \ref{Prop:dev:asymptotique} uses the equations \eqref{Hermite:D} satisfied by $C_k^\lambda$ with $k\leq 3$, that is, the moment equations of the Vlasov-Poisson system up to order three. Hence Proposition \ref{Prop:dev:asymptotique} is also a result on the original system \eqref{vpfp:0}. More precisely, we have
\[
E^\lambda_{\rm slow}\,=\,\sqrt{2}\,T_0\,\partial_x C_2^\lambda\,=\,\partial_x K^\lambda\,-\,T_0\,\partial_x\rho^\lambda\,=\,\partial_x K^\lambda\,-\,\lambda^2\,T_0\,\partial^2_x E^\lambda\,,
\]
where we used the Poisson equation in the last equality. Therefore, under assumption \eqref{compatibility}, $E^\lambda_{\rm slow}$ and $\partial_x K^\lambda$ only differ by a term of order $O(\lambda^{2-\alpha})$, uniformly on $[0,T]$. Replacing $E^\lambda_{\rm slow}$ by $\partial_x K^\lambda$ and $C_1^\lambda$ by $j^\lambda/\sqrt{T_0}$ in \eqref{cE}-\eqref{E:C_1:osc}, Proposition \ref{Prop:dev:asymptotique} gives
\[
\left\|E^\lambda(t)-\partial_x K^\lambda(t)-E^\lambda_{\rm osc}(t)\right\|_{L^2(\T)}\,\leq\, C\,\lambda^{1-\alpha}\,,\qquad
\left\|j^\lambda(t)-j^\lambda_{\rm osc}(t)\right\|_{L^2(\T)}\,\leq\, C\,\lambda^{2-\alpha}\,,
\]
for all $t\in[0,T]$, where $E^\lambda_{\rm osc}$ and $j^\lambda_{\rm osc}$ are the oscillating terms appearing in \eqref{expansion:field:flux}. This corresponds to the rigorous version of the formal expansion \eqref{expansion:field:flux}.
\end{remark}

\begin{remark}
It is worth mentioning that our analysis does not cover the critical case $\alpha\,=\,1$. Actually, in this critical case, the limit of $f^\lambda$ is still valid up to a change of variable, however it is only possible  to characterize the convergence of $E^\lambda$ as 
\[\lambda\left(E^\lambda-E_{\rm osc}^\lambda\right)\longrightarrow\, 0\,,\quad \textrm{as} \quad \lambda \rightarrow 0\,,\]
which means that the slow part of $E^\lambda$ remains not clearly identified, see \cite{Grenier96,Han-Kwan_Hauray15, HK16} for details.
\end{remark}

\section{Fully discrete scheme}
\label{sec3}
\setcounter{equation}{0}
\setcounter{figure}{0}
\setcounter{table}{0} 
In this section, we introduce the numerical scheme used to approximate the Vlasov–Poisson system recast in the Hermite framework \eqref{Hermite:D} and we analyze its properties. We begin by discussing the discretization of the phase space $\T\times \R$ in Section \ref{PS-disc}. Then, in Section \ref{sec:3.2}, we present a first-order time integration scheme based on operator splitting between linear and nonlinear terms. The theoretical properties of the proposed scheme are then analyzed in Section \ref{disc:QN}. Finally, Section \ref{sec:3.4} is devoted to the extension of the first-order method to higher-order accuracy in time. In particular, we propose a second-order time integration scheme based on Strang splitting, combined with an implicit Runge–Kutta method applied to each substep of the splitting.

\subsection{Phase-space discretization}\label{PS-disc}
To discretize the phase space domain, we fix a number of Hermite modes $N_H\in\N^*$. Then, the torus $\T$ is identified with an interval $(a,b)$ of $\mathbb{R}$ with periodic boundary conditions and, for $N_{x}\in\N^\star$, we introduce the discrete set $\J=\{1,\ldots, N_x\}$, the uniform space step $h=(b-a)/N_x$ and the points
$$
x_j  \,=\, a + \left(j-\frac{1}{2}\right) h\,,\qquad j\in\J\,.
$$
We approximate a function $C:(a,b)\to \R$  at point $x_j$,  by $\cC_j$, for $j\in\J$, and we write $C_h=\left(\cC_j\right)_{j\in \cJ}\in \R^\cJ$ the corresponding vector, which is extended by periodicity, that is, $\cC_0=\cC_{N_x}$ and $\cC_{N_x+1}=\cC_1$. We then introduce the discrete operator $\partial_h = (\partial_j)_{j \in\J}$ representing the centered approximation of $\partial_x$ given by
\be
\label{def:Ah}
	\ds\partial_{j} C_h \,=\,  \ds\,\frac{\cC_{j+1} - \cC_{j-1}}{2\, h}  \,, \quad j\in\J\,,
\ee
for all $C_h\,=\,\left(\cC_j\right)_{j\in \cJ}\in \R^\cJ$. Let us notice that $\partial_h$ is skew-symmetric, that is,
  $$
  \sum_{j\in\J} h\,\cC_j\,\partial_j D_h\,=\,-\sum_{j\in\J} h\,\cD_j\,\partial_j C_h
  $$
  for all $C_h =\left(\cC_j\right)_{j\in \cJ}$ and $D_h=\left(\cD_j\right)_{j\in \cJ}\in\R^\cJ$, which is the discrete counterpart of the integration by parts formula on $\T$. Moreover, we assume that $N_x$ is odd, which ensures that discrete Poincaré-Wirtinger and Sobolev inequalities hold uniformly with respect to $h$ (see Lemma \ref{lem:disc:ineq} in Appendix \ref{proof:prop:discrete:QN:limit}). The above detailed choices lead to the following semi-discrete system (where we drop the index $h$ for $C_h$ and only use $C$ for sake of clarity),
\be
\label{semi:discrete}
\left\{
\begin{array}{l}
\ds
\frac{\dD C_{k}}{\dD t}
\,+\,\sqrt{k\, T_0}\,\partial_h\,C_{k-1}
\,+\,\sqrt{(k+1)\,T_0}\,\partial_h\,C_{k+1} \,-\, 
\sqrt{\frac{k}{T_0}}\,E\, C_{k-1}  \,=\,
0\,,
\\[1.em]
\ds
\ds E\, =\,-\partial_h\phi \;;\quad-\lambda^2\partial_h^2\phi\, =\, C_0-1\,,
\end{array}\right.
\ee
for $k\in\{0,\ldots,N_H\}$, where $C_k \in \R^\cJ$ for each $k\in\{0,\ldots,N_H\}$ and with $C_k=0$ when $k>N_H$ and $k=-1$. System \eqref{semi:discrete} is then completed with the condition 
\[
\sum_{j\in\J} h \, \phi_j \,=\,0\,,
\]
ensuring uniqueness of the discrete potential $\left(\phi_j\right)_{j \in \cJ}$. Hence, \eqref{semi:discrete} is a system of ODEs for the $(N_H+1)\,N_x$ unknowns $\left(\cC_{k,j}\right)_{0\leq k\leq N_H,\,j\in\J}$, coupled with the $N_x$ values of the discrete potential $\left(\phi_j\right)_{j\in\J}$. The size of the system is therefore of order $N_H\times N_x$.

\subsection{First order  time splitting scheme}
\label{sec:3.2}
The time discretization is based on a time splitting technique, where the first step consists in solving the Vlasov-Poisson system, linearized around a quasineutral steady state $C^{\rm stat}=\left(\delta_{0,k}\right)_{k\in \N}$. The second step, consists in solving the remaining nonlinear part of the system. More precisely, we first rewrite the semi-discrete system \eqref{semi:discrete} as
\begin{equation}\label{eq:semi}
  	\ds
  	\frac{\dD C_{k}}{\dD t}
  	\,+\,\cL_k\,C\, +\, \cB_k(C) \,=\, 0\,,  
\end{equation}
for each $k\in\{0,\ldots,N_H\}$. We separate the contribution of the linear operator $\cL\,=\,\left(\cL_k\right)_{0\leq k \leq N_H}$ from that of the nonlinear operator $\cB\,=\,\left(\cB_k\right)_{0\leq k \leq N_H}$. The first discrete operator is defined for the modes $C\,=\,\left(C_k\right)_{0\leq k\leq N_H}$ as
  $$
	\ds
	\cL_k\,C\,=\, \sqrt{k\,T_0}\,\partial_h\,C_{k-1}
	\,+\,\sqrt{(k+1) T_0}\,\partial_h\,C_{k+1}\,,\quad \,{\rm for }\,
	k\neq 1\,,
$$  
with $C_k=0$ when $k>N_H$ and $k=-1$ and
$$
	\ds \cL_1\,C\,=\,\sqrt{T_0}\,\partial_h\,C_{0}
	\,+\,\sqrt{2\,T_0}\,\partial_h\,C_{2} \,-\, 
	\frac{1}{\sqrt{T_0}}\,E\,,\quad \textrm{with}\quad E\, =\,-\partial_h\phi\,, \quad \textrm{and}\quad-\lambda^2\partial_h^2\phi\, =\, C_0-1\,.
$$  
The second nonlinear operator $\cB\,=\,\left(\cB_k\right)_{0\leq k \leq N_H}$ is defined by
$$
	\left\{
	\begin{array}{l}
		\ds \cB_{0}(C) \,=\,0\,, \\[1.em]
		\ds \cB_k(C) \,=\,  -\sqrt{\frac{k}{T_0}}\, E \left(C_{k-1}\,-\,\delta_{1,k}\right)\,,\quad{\rm if}\,\,k\geq1\,,\\[1.em]
		\ds E\, =\,-\partial_h\phi\,, \quad \textrm{and}\quad-\lambda^2\partial_h^2\phi\, =\, C_0-1\,,
	\end{array}\right.
$$
where $\delta_{1,k}$ denotes the Kronecker symbol. We then fix a time step $ \Delta t>0 $, set $t^{n}=n \Delta t$ with $n\in\N$, and solve first the linear part on $[t^n, t^{n+1}]$ with a first order fully implicit Euler scheme.
This yields the following linear system to invert 
\be
\label{discrete:step1}
\left\{
\begin{array}{l}
\ds
\frac{C^{(1)}_{1} -C_{1}^{n} }{\Delta t}
\,+\,\sqrt{T_0}\,\partial_h\,C^{(1)}_{0}
\,+\,\sqrt{2\,T_0}\,\partial_h\,C^{(1)}_{2} \,-\, 
\frac{1}{\sqrt{T_0}}\,E^{(1)}\,=\,
0\,,
\\[1.em]
\ds
\frac{C^{(1)}_{k} -C_{k}^{n} }{\Delta t}  \,+\,
\sqrt{k\,T_0}\,\partial_h\,C^{(1)}_{k-1}
\,+\,\sqrt{(k+1) T_0}\,\partial_h\,C^{(1)}_{k+1}
\,=\,0\,,\quad \,{\rm for }\, k\neq 1\,,\\[1.2em]
\ds E^{(1)}\, =\,-\partial_h\phi^{(1)} \;;\quad-\lambda^2\partial_h^2\phi^{(1)}\, =\, C^{(1)}_0-1\,,
\end{array}\right.
\ee
where $C^{(1)}$ indicates the solution obtained after this first splitting-step for $k\in\{0,\ldots,N_H\}$ and, as before, with $C_k^{(1)}=0$ when $k>N_H$ and $k=-1$. To the above system \eqref{discrete:step1}, we add the following uniqueness condition 
\[
\sum_{j\in\J} h \, \phi^{(1)}_j \,=\,0\,.
\]
Let us observe that equations \eqref{discrete:step1} can be recast in matrix-vector form, with a square matrix whose dimension is  $(N_H+2)\, N_x$, where the coefficients of the unknowns are independent of the time index $n$. Therefore, one can perform an $LU$ factorization of the system at $n=0$
and store the resulting $LU$ factors to be reused for all $n\geq 0$, improving computational efficiency. We now focus on the second part of the splitting which involves the nonlinear operator $\cB$. This corresponds to 
$$
\left\{
\begin{array}{l}
  \ds\frac{\dD C_k }{\dD t}   \,+\, \cB_k(C) \,=\, 0\,, \\[1.1em]
  \ds C(t^n) \,=\, C^{(1)}\,.
  \end{array}\right.
  $$
We approximate the above system by using again a fully implicit Euler scheme. This reads
\be
	\label{discrete:step2}
	\left\{
	\begin{array}{lll}
		\ds
		C^{n+1}_{0}\,=\,C^{(1)}_{0}\,,
		\\[1.em]
		\ds
		E^{n+1}\,=\,E^{(1)}\,,
		\\[1.em]
		\ds \frac{C^{n+1}_{k} -C_{k}^{(1)} }{\Delta t}  \,-\,
		\sqrt{\frac{k}{T_0}}\,
		E^{n+1} \left(C^{n+1}_{k-1}\,-\,\delta_{k,1}\right)  
		\,&\ds=\,
	0\,,\quad{\rm if}\,\,k\geq1\,,
	\end{array}\right.
	\ee
for $k\in\{0,\ldots,N_H\}$ and $C_k^{n+1}=0$ when $k>N_H$. Observe that since the electric field $\ds E^{(1)}$ and the density $C_0^{(1)}$ do not change during this second step, the latter system is trivially invertible and hence does not require any particular linear solver. \\

The novelty of this approach lies in discretizing the Poisson equation simultaneously with the equations governing the Hermite coefficients $(C_k)_{0 \leq k \leq N_H}$, as in \cite{abff24, blaustein24}. Let us observe that the strategy above discussed differs from other recent approaches that reformulate the Poisson equation as a harmonic oscillator-type equation for the potential $\phi$, and then build a numerical scheme based on that reformulation (see, for instance, \cite{CDV:07, BCDS:09, ACS, DDNSV:10, DLV, CDV:16, DDD}). In contrast, our method retains the original structure of the Vlasov–Poisson system. However, as we explain later in Proposition~\ref{prop:reformulated_Poisson}, it is still possible to recover a harmonic oscillator formulation for the electric potential from the system \eqref{discrete:step1}–\eqref{discrete:step2} as done in previous works. Ultimately, we emphasize the importance of solving the coefficients $(C_0, C_1, \phi)$ in a completely implicit manner, as they exhibit fast oscillations in the quasineutral regime. 

\subsection{Discrete quasineutral limit}
\label{disc:QN}
In this section, we investigate the theoretical properties of the numerical method \eqref{discrete:step1}–\eqref{discrete:step2}. In Proposition \ref{prop:discrete:QN:limit}, we show that the proposed scheme provides a consistent discretization of the quasineutral limit for a fixed time step $\Delta t$. In addition to this result, in Proposition \ref{prop:reformulated_Poisson}, we demonstrate that the proposed numerical method naturally encodes the harmonic oscillator reformulation of the Poisson equation.
 
We start by proving that a discrete counterpart of Proposition \ref{Prop:dev:asymptotique} holds true with however a fundamental difference: the numerical scheme \eqref{discrete:step1}–\eqref{discrete:step2} filters out the oscillations when the parameter $\lambda$ tends to $0$ and $\Delta t$ is fixed. This is a consequence of the dissipative nature of the implicit scheme. More precisely, we prove that, for a fixed time step and under a discrete version of the assumption \eqref{compatibility}, the numerical method \eqref{discrete:step1}–\eqref{discrete:step2} consistently approximates the quasineutral limit in the sense that:  
\[
E^n \underset{\lambda \rightarrow 0}{=} \sqrt{2}\,T_0\, \partial_h\,C^{n}_{2}+ O(\lambda)\,,\quad \textrm{and}\quad 
C_1^n \underset{\lambda \rightarrow 0}{=}  O(\lambda)\,.
\]
To state our result, we use the discrete $l^2$ and $h^r$ norms of $C\,=\,\left(\cC_j\right)_{j\in \cJ}$, defined as
\[
\left\|C\right\|_{l^2\left(\cJ\right)}^2
\,=\,
\sum_{j\in \cJ}h \,\left|\cC_j\right|^2
\,,\quad \textrm{and}\quad 
\left\|C\right\|_{h^r\left(\cJ\right)}^2\,=\,
\sum_{0\leq s \leq r}
\left\|\partial_h^s C\right\|_{l^2\left(\cJ\right)}^2\,.
\]
We also introduce the discrete analog of $\cE_0$ and $\cE_1$ of Proposition \ref{Prop:dev:asymptotique}, in which we remove the oscillatory part of the solution. They read 
 \begin{equation}
	\label{cE:h}
	\left\{
	\begin{array}{l}
		\ds\cE_0^n \,:=\, \left\|E^n -  \sqrt{2}\,T_0\, \partial_h\,C^{n}_{2} \right\|_{l^2(\cJ)}\,,
		\\[1.1em]
		\ds\cE_1^n \,:=\,  \left\|C_1^n \right\|_{h^1(\cJ)}\,.
	\end{array}\right.
\end{equation}
We are now ready to prove the following result.
\begin{proposition}\label{prop:discrete:QN:limit}
Let us consider a fixed $\Delta t > 0$, an odd number of grid points $N_x$ and a solution $(C^n,\phi^n)_{n\geq 0}$ to \eqref{discrete:step1}-\eqref{discrete:step2} with zero total flux and with total mass $|b-a|$, that is,
\begin{equation}\label{total:flux:mass:h}
\sum_{j\in\J}\cC_{1,j}^{0}\,h\,=0\,,\quad \textrm{and}\quad \sum_{j\in\J}\cC_{0,j}^{0}\,h\, = |b-a| \,.
\end{equation}
Suppose furthermore that $(C^n,\phi^n)_{n\geq 0}$ satisfies the discrete analog of \eqref{compatibility} with $\alpha\,=\, 1$, that is, 
\begin{equation}\label{hyp:E:disc}
	\sup_{\lambda >0}\;
	\sup_{0\,\leq \,n\,\leq\, T/\Delta t}\left(
	\left\|C^n_1\right\|_{h^1\left(\cJ\right)}
	\,+\,	\lambda
	\left\| E^n\right\|_{h^3\left(\cJ\right)}
	\right)
	 < +\infty\,,
\end{equation}
for some given final time $T>0$, and that the discrete analog of \eqref{uniform:reg} is satisfied, that is, 
\begin{equation}\label{hyp:Ck:disc}
	\sup_{\lambda>0}\;
	\sup_{0\leq n\leq T/ \Delta t} \;
	\sum_{2\leq k \leq N_H}
	 \left\|C^n_k\right\|^2_{h^2\left(\cJ\right)} \;< +\infty\,.
\end{equation}
Then, we have for all $\lambda >0$
\begin{equation}\label{result:h}
\ds \sup_{2\leq n\leq T/\Delta t}\;\cE_0^n \,\leq\, C\lambda \,,
\quad {\rm and }\quad
\sup_{1 \leq n \leq T/\Delta t}
\cE_1^n
\,\leq\,
C\lambda\,,
\end{equation}
where $C>0$ depends on $\Delta t$, $T$, $T_0$, $|b-a|$ and on the implicit constants in \eqref{hyp:E:disc}-\eqref{hyp:Ck:disc}. Furthermore, the constant $C$ is uniform with respect to $\lambda$ and the phase-space discretization parameters $h>0$ and $N_H \geq 2$.
\end{proposition}
The proof is detailed in Appendix \ref{proof:prop:discrete:QN:limit}. The key point consists in proving that for $\Delta t > 0 $, the scheme \eqref{discrete:step1}-\eqref{discrete:step2} filters the fast oscillations of $E^n$ and $C_1^n$ around the quasineutral state. This permits to show that the amplitude of $E^n$ and $C_1^n$, prescribed by \eqref{hyp:E:disc}, are reduced by a factor $\lambda$ for all $n\geq 1$. The main mathematical difficulty of the proof arises in the analysis of the nonlinear step \eqref{discrete:step2}. This requires the use of the discrete Sobolev injection $h^1\left(\cJ\right)\hookrightarrow l^\infty\left(\cJ\right)$. Let us point out that Proposition \ref{prop:discrete:QN:limit} is a consistency result: it assumes that the numerical solution satisfies uniform bounds with respect to $\lambda$, that is, a uniform stability property, and proves that, under these bounds, the scheme captures the quasineutral limit. The proof of the uniform stability of the scheme is an open question, which is investigated numerically in Section \ref{sec4}.
\begin{remark}
	We emphasize that our result remains valid in the critical case $\alpha = 1$, which was not treated in the continuous setting. In this critical case, the work \cite{Grenier96} shows that $f^\lambda$ converges to the quasineutral limit up to an oscillating change of variable in velocity. This means that the numerical method filters this change of variable and the numerical solution converges to the quasineutral limit directly. Let us however emphasize that assumptions \eqref{hyp:E:disc}-\eqref{hyp:Ck:disc} are made on the numerical solution itself: they require that the numerical solution satisfies uniform bounds with respect to $\lambda$, which is not guaranteed when $\Delta t$ is of the same order as $\lambda$ (see Section \ref{sec:4.1}).
\end{remark}

  \begin{remark}
    \label{rem:hyp:disc}
    Let us comment on the differences between assumptions \eqref{hyp:E:disc}-\eqref{hyp:Ck:disc} and their continuous counterparts \eqref{compatibility}-\eqref{uniform:reg}.  Proposition \ref{prop:discrete:QN:limit} is established for a fixed time step $\Delta t$, whereas Proposition \ref{Prop:dev:asymptotique} relies on the resolution of the fast time scale $\lambda$.  On the one hand, in the continuous setting, the oscillating parts $E^\lambda_{\rm osc}$ and $C^\lambda_{1,{\rm osc}}$ are not damped and the remainder is controlled only when $\alpha<1$. In the discrete setting, the time step $\Delta t$ is fixed and the fast time scale is not resolved: observing that
    $$
    C_0^n-1=\lambda^2\,\partial_h E^n=O(\lambda)
    $$
    for all $n$, the discrete time derivative $(C_0^{n+1}-C_0^n)/\Delta t$ is of order $\lambda/\Delta t$, which, through the implicit treatment of the continuity and flux equations in \eqref{discrete:step1}, provides a control of $C_1^{(1)}$ and $E^{n+1}$ without computing the oscillating part of the solution. Hence the critical case $\alpha=1$ can be treated theoretically. This is consistent with the linear behavior of the implicit Euler step \eqref{discrete:step1}, which damps the oscillations with frequency $\lambda^{-1}$ by a factor $(1+\Delta t^2/\lambda^2)^{-1/2}$ at each time step. On the other hand,  in the continuous setting, the bound on $C_1^\lambda$ is only assumed at the initial time, but it must be small, that is, of order $\lambda^{1-\alpha}$, whereas for $t>0$ the behavior of $C_1^\lambda$ is a consequence of the analysis. In the discrete setting, we only assume that $C_1^n$ is bounded uniformly with respect to $\lambda$, without any smallness assumption, but this bound is required for all $n$ since the analysis is performed step by step, for a fixed $\Delta t$, without any time integration over the fast scale. As the uniform bound on $(C_k^n)_{k\geq 2}$ in \eqref{hyp:Ck:disc}, it should be understood as a uniform regularity assumption on the numerical solution.
\\
Finally, in the continuous setting, the iteration of the Duhamel formula requires $r_0+3$ derivatives, with $r_0=\lceil 1/(1-\alpha)\rceil$, and the control of quadratic terms such as $|E^\lambda_{\rm osc}|^2$ in $L^2$, which leads to the use of $L^4$-based Sobolev norms. In the discrete setting, no iteration is needed: the nonlinear terms are estimated with the discrete Sobolev injection $h^1(\cJ)\hookrightarrow l^\infty(\cJ)$, which holds uniformly with respect to $h$ in dimension one when $N_x$ is odd, hence the $l^2$-based norms $h^r(\cJ)$ with $r\leq 3$ are sufficient.
\end{remark}

As already mentioned, most of the asymptotic preserving schemes proposed in the literature for the quasineutral limit \cite{CDV:07, BCDS:09, DDNSV:10, CDV:16, DD:17} rely on a reformulation of the Poisson equation as a harmonic oscillator equation. The key point is that, contrary to the original Poisson equation which degenerates when $\lambda=0$, this reformulated equation remains non degenerate in the limit $\lambda\rightarrow 0$ and provides the quasineutral electric field. In the continuous Hermite framework, this reformulated Poisson equation can be derived by differentiating the Poisson equation (second line of \eqref{Hermite:D}) twice in time and substituting $C^\lambda_{0}$ and $C^\lambda_{1}$ using their corresponding evolution equations (first line of \eqref{Hermite:D} for $k=0$ and $1$). These computations lead to the so-called reformulated Poisson equation:
\be\label{continuous:reformul}
\lambda^2\, \partial^2_t(\partial_{x}E^\lambda)+\partial_{x}(C^\lambda_{0}\,E^\lambda)=\partial^2_x (\sqrt{2}\,T_0\,C^\lambda_{2}+T_0\,C^\lambda_{0})\,.
\ee
Indeed, when $\lambda = 0$ and $C_0^\lambda = 1$, equation \eqref{continuous:reformul} reduces to $\partial_x E = \partial_x^2 (\sqrt{2}\,T_0\,C_2+T_0\,C_0)$, which is the quasineutral relation $E=\sqrt{2}\,T_0\,\partial_x C_2$ given in \eqref{Hermite:D:neutral}. It is then natural to ask whether an analogous discrete reformulated Poisson equation can be derived from the proposed scheme, even if it is not explicitly used in its construction. This is the object of the following result.
\begin{proposition}\label{prop:reformulated_Poisson} Let $(C^n)_{n\geq 0}$ be a solution to \eqref{discrete:step1}-\eqref{discrete:step2} with $E=-\partial_h \phi$. Then, the following reformulated discrete Poisson equation for all $n\geq 1$ holds true:
	\be\label{discrete:reformul}
	\begin{split}
		\lambda^2\,\frac{\partial_hE^{n+1}-2\,\partial_hE^{n}+\partial_hE^{n-1}}{\Delta t^2}\,+\,\partial_h\,\Bigl(E^{n+1}\,C^{n+1}_{0}\Bigl) \hspace{6cm}\\
		\hspace{3cm}\,=\,\partial^2_h
		\left(\sqrt{2}\,T_0\,C^{(1,n+1)}_{2}+T_0\,C^{n+1}_{0}\right)\,+\,\lambda^2\,\Delta
		t\,\partial_h\,\left(\frac{E^{n+1}\, \partial_hE^{n+1}-E^{n}\,
			\partial_hE^{n}}{\Delta t}\right),
	\end{split}
	\ee
	where $C^{(1,n+1)}_{2}$ corresponds to the Hermite coefficient $C^{(1)}_{2}$ computed during the step \eqref{discrete:step1} of the time splitting scheme to get $(C_k^{n+1})_{k\geq 0}$ from~$(C_k^{n})_{k\geq 0}$.
\end{proposition}
The proof is postponed in Appendix \ref{proof:prop:reformulated_Poisson}.

Proposition \ref{prop:reformulated_Poisson} is used neither in the implementation of the scheme nor in the proof of Proposition \ref{prop:discrete:QN:limit}. Its interest is to establish a link between our approach and the previous ones \cite{CDV:07, BCDS:09, DDNSV:10, CDV:16, DD:17}: although the scheme \eqref{discrete:step1}-\eqref{discrete:step2} only discretizes the original Poisson equation, the discrete reformulated Poisson equation \eqref{discrete:reformul} is automatically satisfied. Equation \eqref{discrete:reformul} is a consistent discretization of \eqref{continuous:reformul}, up to the additional term of order $O(\lambda^2\,\Delta t)$ in the right hand side, where the field $E^{n+1}$ is treated implicitly. Moreover, it provides another explanation of the asymptotic preserving property of our scheme. Indeed, setting formally $\lambda=0$ in \eqref{discrete:reformul}, we get $C_0^{n+1}=1$ and
\[
\partial_h E^{n+1}\,=\,\partial_h^2\left(\sqrt{2}\,T_0\,C_2^{(1,n+1)}\,+\,T_0\,C_0^{n+1}\right)\,=\,\sqrt{2}\,T_0\,\partial_h^2\,C_2^{(1,n+1)}\,,
\]
which is a consistent discretization of the quasineutral relation $E=\sqrt{2}\,T_0\,\partial_x C_2$ in \eqref{Hermite:D:neutral}, in agreement with Proposition \ref{prop:discrete:QN:limit}.


\subsection{Second order time splitting scheme}
\label{sec:3.4}
In this section, we extend the first-order time method described in \eqref{discrete:step1}–\eqref{discrete:step2} to a second-order discretization. To this end, observe that each step of the time-splitting scheme \eqref{discrete:step1}–\eqref{discrete:step2} consists of a first-order implicit method. Therefore, this approach can be generalized to higher-order schemes by employing existing Runge–Kutta techniques \cite{BFR, pr:imex} and high-order splitting schemes. In particular, we apply a second-order time-splitting strategy (Strang splitting), combined with a second-order stiffly accurate implicit Runge–Kutta method. Thus, given the equation
\begin{equation}\label{eq:semi1}
	\frac{\dD C_k}{\dD t} + \mathcal{L}_k C + \mathcal{B}_k(C) = 0\,,  
\end{equation}
for all coefficients $k \in \{0, \ldots, N_H\}$ of the Hermite expansion, we proceed as follows. For each time step starting from $C^n$, we first solve the linear part on half time step $\Delta t / 2$ by taking as initial data $C^{(0)}=C^n$:
\begin{equation*}
	\frac{\dD C_k}{\dD t} + \mathcal{L}_k C = 0, 
\end{equation*}
for all coefficients $k \in \{0, \ldots, N_H\}$. We set $C^{(1)}$ the solution obtained after this half time step and then we solve the nonlinear part on a full time step $\Delta t$ using $C^{(1)}$ as initial data. This reads
\begin{equation*}
	\frac{\dD C_k}{\dD t} + \mathcal{B}_k (C) = 0\,.
\end{equation*}
Finally, by setting $C^{(2)}$ the solution obtained after this second step as a new initial data, we conclude the splitting procedure by solving the linear part on another half time step $\Delta t/2$:
\begin{equation*}
	\frac{\dD C_k}{\dD t} + \mathcal{L}_k C = 0.
\end{equation*}
This gives $C^{(3)}=C^{n+1}$, the solution at time $t^{n+1}$. Moreover, each of the substeps above is solved using a second-order, diagonally implicit, stiffly accurate Runge–Kutta scheme, defined as follows:  
\[
\left\{
\begin{array}{l}
	\displaystyle K_1 \,=\, \mathcal{F}\left( C^{(i)} - \gamma\,\overline{\Delta t} \, K_1 \right), \\[1em]
	\displaystyle K_2 \,=\, \mathcal{F}\left( C^{(i)} - (1 - \gamma)\,\overline{\Delta t} \, K_1 - \gamma\,\overline{\Delta t} \, K_2 \right), \\[1em]
	\displaystyle C^{(i\,+\,1)} \,=\, C^{(i)} \,-\, (1 - \gamma)\,\overline{\Delta t} \,  K_1 \,-\, \gamma\,\overline{\Delta t} \, K_2\,,
\end{array}
\right.
\]
where $\mathcal{F} \in \{ \mathcal{L}, \mathcal{B} \}$, $\gamma = 1 - \frac{1}{\sqrt{2}}$, and $\overline{\Delta t}$ denotes $\Delta t$ in the case of the nonlinear component and $\Delta t/2$ in the case of the two half-steps for the linear component. Finally the superscripts $(i)$ and $(i+1)$ indicate respectively the value of the Hermite coefficients at the beginning and at the end of the Strang splitting substep.

\begin{remark}\label{rem:second:order}
The second order scheme shares the structure of the first order one, which is responsible for its asymptotic behavior as $\lambda\rightarrow 0$. Indeed, since the Runge-Kutta method is diagonally implicit, each stage of each substep is an implicit Euler step. We therefore expect the second order scheme to capture the quasineutral limit as well. However, a rigorous proof is beyond the scope of this paper, but this property is supported by the numerical simulations of Section \ref{sec4}.
\end{remark}

\section{Numerical simulations}
\label{sec4}
\setcounter{equation}{0}
\setcounter{figure}{0}
\setcounter{table}{0}

The numerical experiments are conducted using the second order method presented in the previous Section \ref{sec:3.4}, which shares the structure of the first order scheme analyzed in Proposition \ref{prop:discrete:QN:limit} (see Remark \ref{rem:second:order}). The temperature $T_0$ is fixed to $T_0=1$ for all tests, while the other numerical parameters are chosen based on the specific test case in order to accurately capture the physical phenomena under investigation.

\subsection{Near equilibrium}
\label{sec:4.1}
We first consider an initial distribution for the Vlasov–Poisson system \eqref{vpfp:0}, consisting of a perturbation from a global Maxwellian state of order $\delta \lambda^{2 - \alpha}$, with $\alpha\in [0,1)$, given by
\begin{equation}\label{ini:cond:equi}
	\ds f^{\lambda}_{\rm in}(x,v)\,=\,
	\ds
	\frac{1}{\sqrt{2\pi}}\,\left(1+\delta\lambda^{2-\alpha}\cos{\left(k_x x\right)}\right)\exp\left(
			-\frac{|v|^2}{2}\right)\,.
                  \end{equation}
The spatial domain is $x\in[-10,10]$, the frequency $k_x = \pi/10$ while we fix the size of the perturbation to $\delta=0.1$\,. For this initial data and when $\alpha\in [0,1)$, it has been shown in \cite{ Han-Kwan_Hauray15} that in the limit $\lambda\rightarrow 0$, the solution $f^\lambda$ to \eqref{vpfp:0} converges toward $f^0$ defined as
\begin{equation}\label{stat:state}
	f^0(v)\;=\, \frac{1}{\sqrt{2\pi}}\,\exp\left(-\frac{|v|^2}{2}\right),
\end{equation}
whereas $E^\lambda$  weakly converges to $E^0=0$, which corresponds to the stationary solution of the
quasineutral model \eqref{vpfp:neutral}.  \\
First, a reference solution is computed to assess the numerical order of convergence of the quantities $\cE_0^\lambda$ and $\cE_1^\lambda$, defined in \eqref{cE}-\eqref{E:C_1:osc}, describing the oscillatory components of the
electric field $E^\lambda$ and $C_1^\lambda$, respectively. This solution uses $N_H = 128$ Hermite modes, $N_x =2049$ mesh points in the physical space and $\Delta t = 10^{-4}$ is chosen to resolve the smallest value of the Debye length considered in this test, that is, $\lambda = 10^{-2}$ .
 \begin{figure}[ht!]
 \centering
	\begin{tabular}{cc}
          \includegraphics[width=7.5cm]{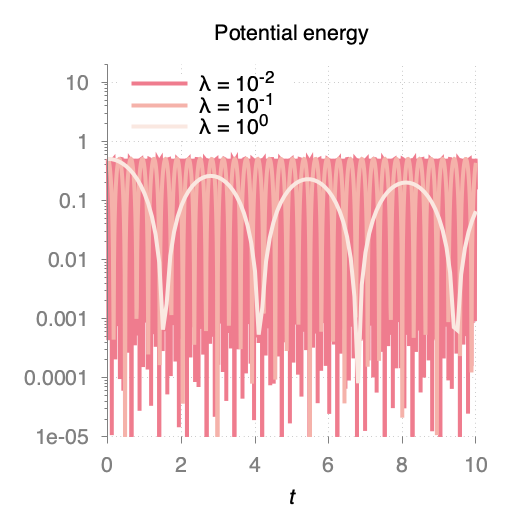}&
          \includegraphics[width=7.5cm]{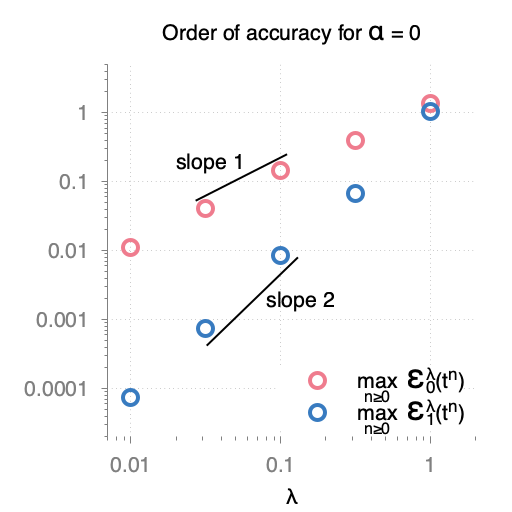}
                                                           \\
          (a) &    (b)
          \\
          \includegraphics[width=7.5cm]{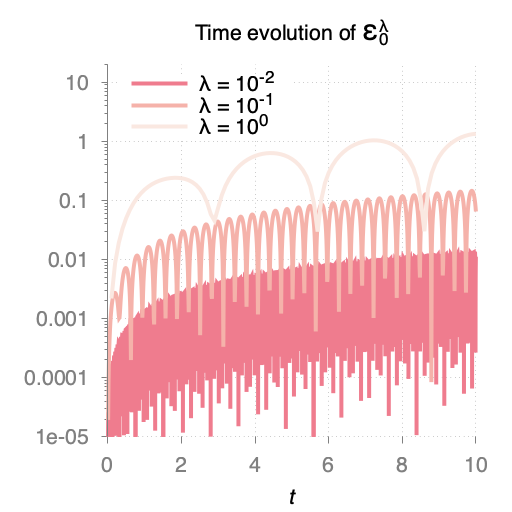}&
          \includegraphics[width=7.5cm]{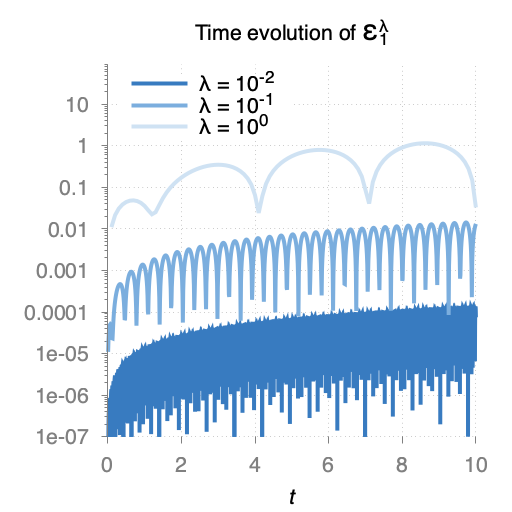}
          \\
          (c) &    (d)
         \end{tabular}
	\caption{{\bf Near equilibrium test case with $\alpha=0$ :}
         (a) time evolution of the scaled potential energy
         $\dfrac12\sum_{j\in\J} h \,|E_j^n|^2$;  (b) order of
         convergence for $\max_{n\geq 0}{\mathcal E}^\lambda_0(t^n)$ and
         $\max_{n\geq 0}{\mathcal E}^\lambda_1(t^n)$ ; 
          (c) time evolution of
	  ${\mathcal E}^\lambda_0(t)$ and  
        $(d)$ time evolution of ${\mathcal E}^\lambda_1(t)$, as defined in
\eqref{cE}-\eqref{E:C_1:osc}, for different values of $\lambda$.}
      \label{fig:10}
    \end{figure}

In Figure \ref{fig:10}-(a) the time evolution of the potential
energy
$$\frac{1}{2}\,\sum_{j\in\J}h\,
|E^n_j|^2\,,
$$
is shown in the case $\alpha=0$. The oscillatory nature of the electric field with a frequency inversely proportional to $\lambda$ appears clearly. Figure \ref{fig:10}-(b) represents the quantities
$$
\max_{n\geq 0} \cE_0^\lambda(t^n)  \qquad {\rm and } \qquad \max_{n\geq 0} \cE_1^\lambda(t^n)\,,
$$
defined in \eqref{cE}-\eqref{E:C_1:osc} for different values of $\lambda$. This illustrates Proposition \ref{Prop:dev:asymptotique}, which established first and second order of convergence for $\cE_0^\lambda$ and $\cE_1^\lambda$ respectively. Moreover, these figures exhibit that $E^\lambda$ converges to its quasineutral limit $E^0$ up to the fast oscillations in time of $E^\lambda_{\rm osc}$ which are explicitly given in terms of the initial condition by \eqref{E:C_1:osc}.

\begin{figure}[ht!]
 \centering
 \begin{tabular}{cc}
          \includegraphics[width=7.5cm]{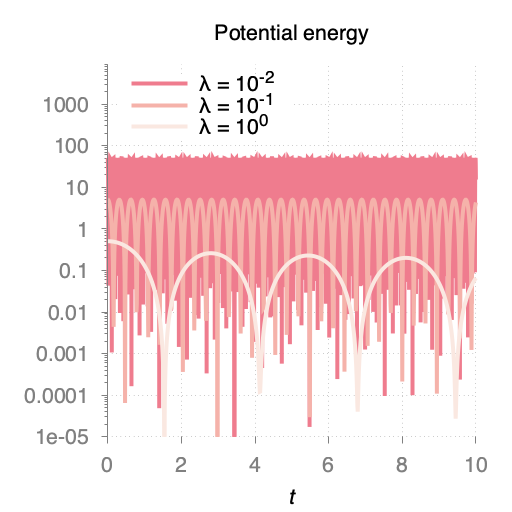}&
          \includegraphics[width=7.5cm]{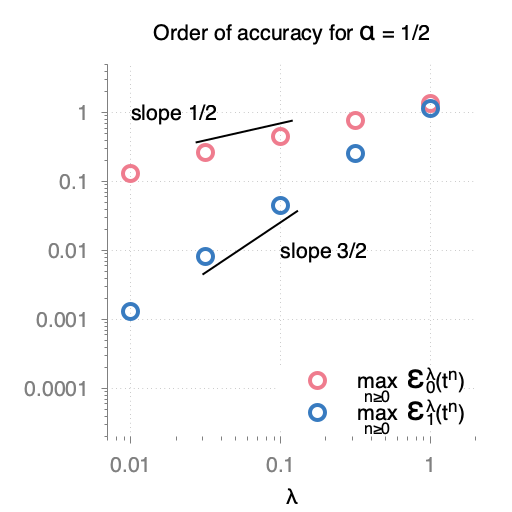}
            \\
          (a) &    (b)\\
          \includegraphics[width=7.5cm]{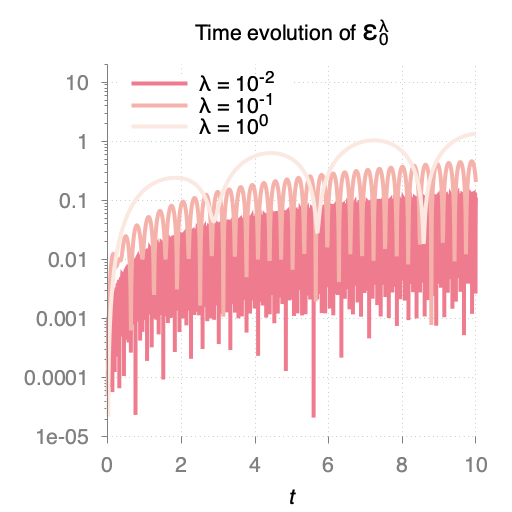}&
          \includegraphics[width=7.5cm]{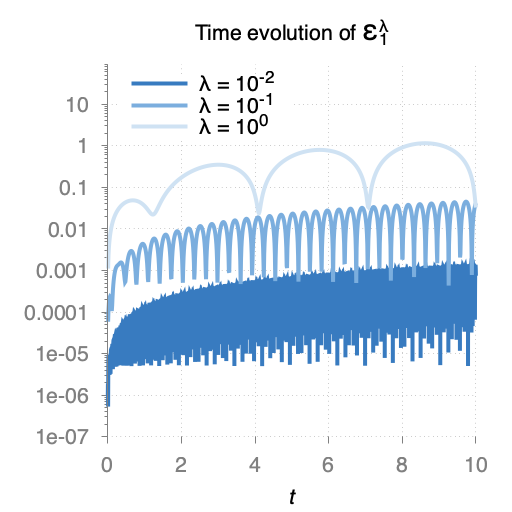}
                                                           \\
          (c) &    (d)
 \end{tabular}
	\caption{{\bf Near equilibrium test case with $\alpha=1/2$ :}
	(a) time evolution of the scaled potential energy
	$\dfrac12\sum_{j\in\J} h \,|E_j^n|^2$;  (b) order of
	convergence for $\max_{n\geq 0}{\mathcal E}^\lambda_0(t^n)$ and
	$\max_{n\geq 0}{\mathcal E}^\lambda_1(t^n)$ ; 
	(c) time evolution of
	${\mathcal E}^\lambda_0(t)$ and  
	$(d)$ time evolution of ${\mathcal E}^\lambda_1(t)$, as defined in
	\eqref{cE}-\eqref{E:C_1:osc}, for different values of $\lambda$.}
      \label{fig:11}
\end{figure}
Figure \ref{fig:11} displays the same quantities as Figure \ref{fig:10}, in the case of $\alpha=1/2$. We observe that, as predicted by the theory, the amplitude of the electric field is now of order $O(\lambda^{-1/2})$ as $\lambda$ approaches zero (Figure \ref{fig:11}-(a)). Our numerical outcomes correspond to the theoretical results of Proposition \ref{Prop:dev:asymptotique} regarding the oscillation frequency in time and the order of convergence of the quantity $\cE^\lambda_0$ and $\cE^\lambda_1$ (Figure \ref{fig:11}-(b)). Let us observe that, while the norm of $E^\lambda$ grows as $\lambda\rightarrow 0$, the electric field
$E^\lambda$ strongly converges to its quasineutral limit $E^0$ up to explicitly known fast oscillations in time.\\

\begin{figure}[ht!]
 \centering
	\begin{tabular}{cc}
          \includegraphics[width=7.5cm]{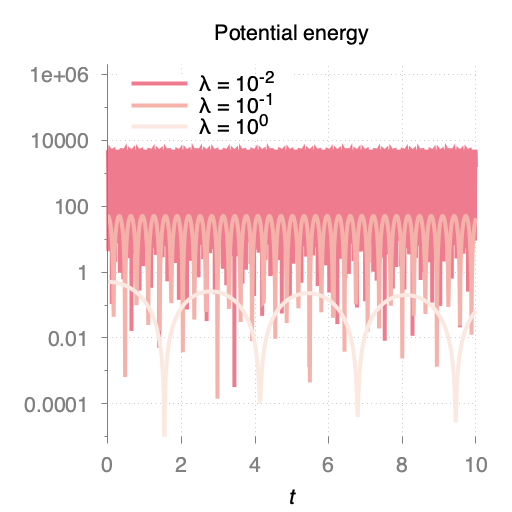}&
         \includegraphics[width=7.5cm]{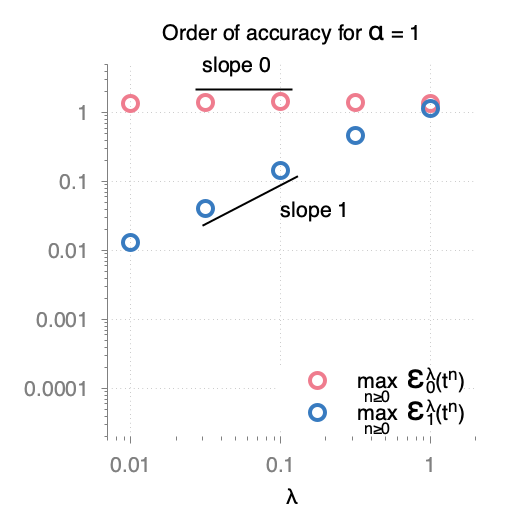}
                                                           \\
          (a) &    (b)  \\
          \includegraphics[width=7.5cm]{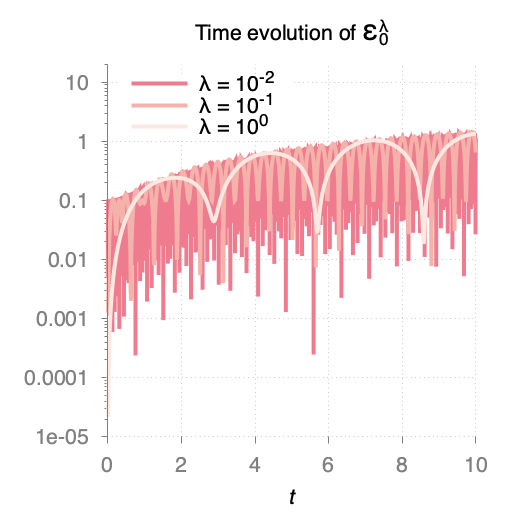}&
         \includegraphics[width=7.5cm]{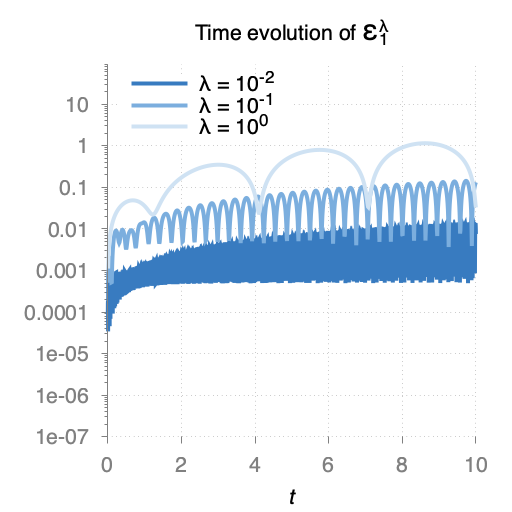}
                                                           \\
          (c) &    (d)
        \end{tabular}
	\caption{{\bf Near equilibrium test case with $\alpha=1$ :}
	(a) time evolution of the scaled potential energy
	$\dfrac12\sum_{j\in\J} h \,|E_j^n|^2$;  (b) order of
	convergence for $\max_{n\geq 0}{\mathcal E}^\lambda_0(t^n)$ and
	$\max_{n\geq 0}{\mathcal E}^\lambda_1(t^n)$ ; 
	(c) time evolution of
	${\mathcal E}^\lambda_0(t)$ and  
	$(d)$ time evolution of ${\mathcal E}^\lambda_1(t)$, as defined in
	\eqref{cE}-\eqref{E:C_1:osc}, for different values of $\lambda$.}
      \label{fig:12}
\end{figure}
We then consider the case $\alpha=1$. This situation is not in the scope of our theoretical analysis but it allows us to test its range of validity. The results are shown in Figure \ref{fig:12}. The amplitude of the electric field is now of order $O(\lambda^{-1})$ while $\max_{n\geq 0}\cE_0^\lambda(t^n)$ no longer
converges to zero as $\lambda\rightarrow 0$ (Figure \ref{fig:12}-(b) and (c)). This is in line with theoretical investigations \cite{Grenier96,Han-Kwan_Hauray15, HK16}, which show that $\lambda\, \cE^\lambda_0 \rightarrow 0$. In contrast, the quantity $\cE_1^\lambda$ remains of the order of $\lambda$ (Figure \ref{fig:12}-(b) and (d)). These results suggest that the estimate on $\cE_1^\lambda$ in Proposition \ref{Prop:dev:asymptotique} may remain true with $\alpha=1$.

\begin{figure}[ht!]
 \centering
	\begin{tabular}{ccc}
          \includegraphics[width=5.cm]{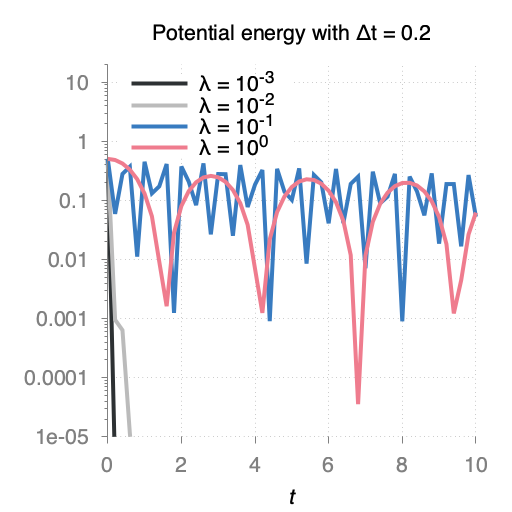}&
           \includegraphics[width=5.cm]{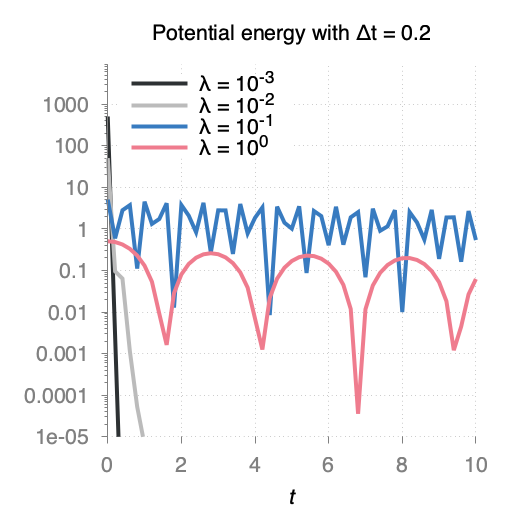}&
        \includegraphics[width=5.cm]{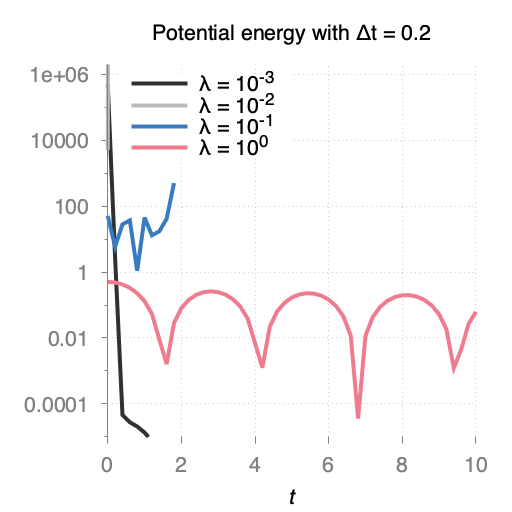}
                                                           \\
          (a) $\alpha=0$ &    (b) $\alpha=1/2$   & (c) $\alpha=1$
        
         \end{tabular}
	\caption{{\bf Near equilibrium test case.} Time evolution of the scaled potential energy  $\dfrac12\sum_{j\in\J} h \,|E_j^n|^2$ with $\Delta t=0.2$, $N_x=65$ and $N_H=128$ for different $\lambda$: (a)
            $\alpha=0$ ;  (b) $\alpha=1/2$ ;  (c) $\alpha=1$. } 
      \label{fig:13}
\end{figure}

Finally, we investigate the behavior of our numerical scheme when the time step $\Delta t$ and the space mesh size are fixed while the quasineutral parameter $\lambda$ tends to zero. The aim is to study the ability of the scheme to capture the asymptotic behavior of the solution and to illustrate the asymptotic preserving property shown in Proposition \ref{prop:discrete:QN:limit}. The numerical parameters are $\Delta t= 0.2$, $N_x=65$ and $N_H=128$ and Figure \ref{fig:13} displays the time evolution of the scaled potential energy $\|E^n\|_{l^2}^2/2$ for different values of $\alpha$ and $\lambda$. For $\alpha=0$ and $\alpha=1/2$, corresponding to Figure \ref{fig:13}-(a) and (b), the time evolution of the scaled potential energy $\|E^n\|_{l^2}^2/2$ gives results in line with the theoretical findings even for a large time step $\Delta t= 0.2$, both in terms of amplitude and frequency of oscillations. However, when $\lambda$ becomes smaller ($\lambda =10^{-2}$ and  $10^{-3}$), these discretization parameters are no longer sufficient to describe the small time scales. Hence, the electric field converges toward the expected weak limit. This behavior illustrates Proposition \ref{prop:discrete:QN:limit}: in this test case, the slow component $\sqrt{2}\,T_0\,\partial_h C_2^n$ remains very small, so that estimate \eqref{result:h} indicates that the electric field $E^n$ becomes of order $O(\lambda)$ after the first time steps, and the scaled potential energy of order $O(\lambda^2)$, which is indeed observed in Figure \ref{fig:13}-(a), (b) for $\lambda=10^{-2}$ and $10^{-3}$, corresponding to $\alpha= 0$ and $1/2$. In the case $\alpha=1$ and $\lambda=10^{-1}$ or $\lambda=10^{-2}$, corresponding to Figure \ref{fig:13}-(c), we observe that the numerical solution becomes unstable and blows up. In this critical regime, where $\lambda$ is not small enough compared to $\Delta t$, the uniform bounds \eqref{hyp:E:disc}-\eqref{hyp:Ck:disc} are not satisfied by the numerical solution, and the time step and spatial mesh need
	to be refined to obtain stable and accurate results. These cases are not covered by Proposition \ref{prop:discrete:QN:limit}, which holds for a fixed $\Delta t$ : the constant in \eqref{result:h} depends on $\Delta t$ and is not expected to be uniformly bounded in the joint limit $(\lambda,\Delta t)\rightarrow 0$, since the estimate \eqref{result:h} is not consistent with its continuous analog in Proposition \ref{Prop:dev:asymptotique}. Furthermore, it is interesting to notice that for $\lambda$ small enough ($\lambda=10^{-3}$ in this case) and $\Delta t=0.2$, the numerical solution remains bounded and the scaled potential energy is instantaneously damped, which means that the approximate electric field $E^{n}$ converges to zero, in agreement with Proposition \ref{prop:discrete:QN:limit}.

\subsection{Smooth perturbation of equilibrium}
\label{sec:4.2}
We now consider an initial distribution with non-homogeneous temperature driven by a perturbation of order $O(1)$ as $\lambda \rightarrow 0$. This reads
\begin{equation}
  \label{ini:cond:qn:damping}
	\ds f_{\rm in}(x,v)\,=\,
	\ds
	\frac{1}{\sqrt{2\pi\,T_{\rm in}(x)}}\, \exp\left(
          -\frac{|v|^2}{2\,T_{\rm in}(x)} \right)\,,
                  \end{equation}
with $x\in[-10,10]$, frequency $k_x =\pi/10$ and $T_{\rm in}$ given by
$$
T_{\rm in}(x) = 1+\delta\,\cos\left(k_x x\right),
$$
with $\delta = 0.1$ . This gives for the density and the flux the following relations
$$
\rho_{\rm in}(x) = \int_\R f_{\rm in}(x,v) \dD v \,=\, 1 \quad{\rm and}\quad
j_{\rm in}(x) = \int_\R f_{\rm in}(x,v) v\,\dD v \,=\, 0\,, \qquad \forall x
\in (-10,10)\,.
$$
In this situation, the quasineutral model \eqref{vpfp:neutral} has non trivial slow dynamics in time that we wish to describe using the proposed numerical approximation. 

\begin{figure}[ht!]
 \centering
	\begin{tabular}{cc}
          \includegraphics[width=7.5cm]{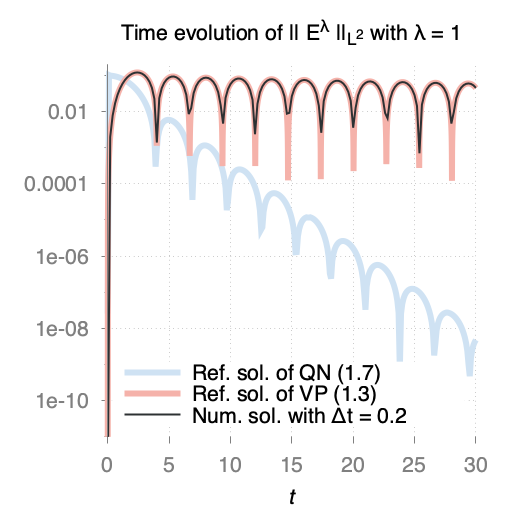}&
          \includegraphics[width=7.5cm]{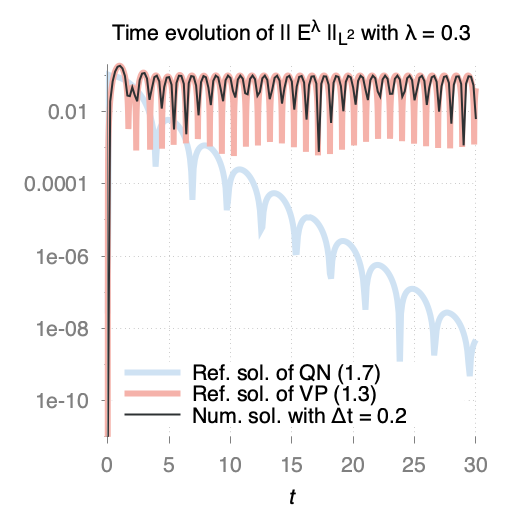}
                                                           \\
          (a) &    (b)
          \\
          \includegraphics[width=7.5cm]{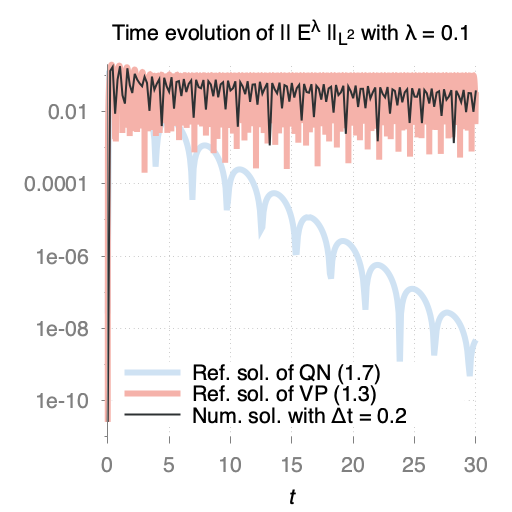}&
          \includegraphics[width=7.5cm]{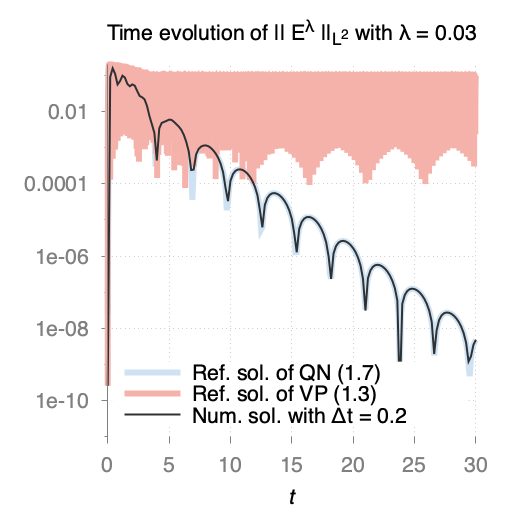}
          \\
          (c) &    (d)
        \end{tabular}
        \caption{{\bf Smooth perturbation of equilibrium test case:} time evolution of  $\|E^\lambda\|_{L^2}$ in logarithmic
          scale with $(a)$ $\lambda= 1$;  $(b)$ $\lambda =
        3. \;10^{-1}$;  $(c)$ $\lambda = 10^{-1}$ and $(d)$ $\lambda = 3. 10^{-2}$. Red (thick) curve: reference solution of the Vlasov-Poisson system \eqref{vpfp:0} computed on a fine mesh; blue (thick) curve: reference solution of the quasineutral model \eqref{vpfp:neutral}; black (thin) curve: numerical solution with $\Delta t=0.2$.}
      \label{fig:20}
    \end{figure}

\begin{figure}[ht!]
 \centering
	\begin{tabular}{cc}
          \includegraphics[width=7.5cm]{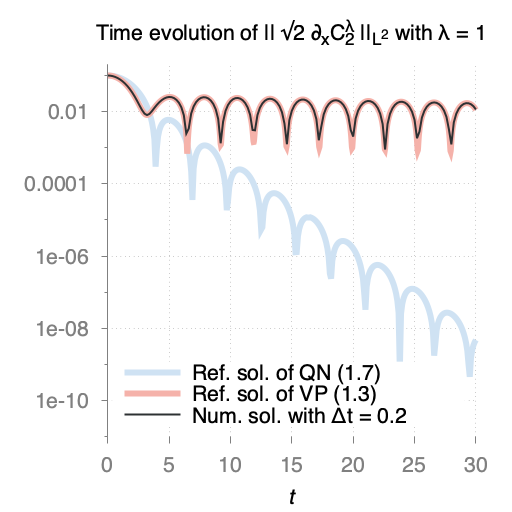}&
          \includegraphics[width=7.5cm]{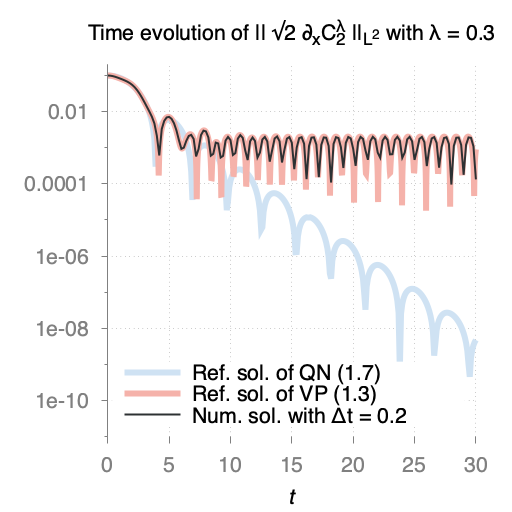}
                                                           \\
          (a) &    (b)
          \\
          \includegraphics[width=7.5cm]{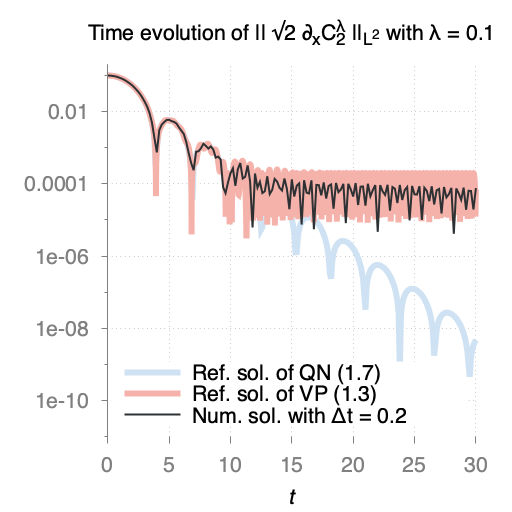}&
          \includegraphics[width=7.5cm]{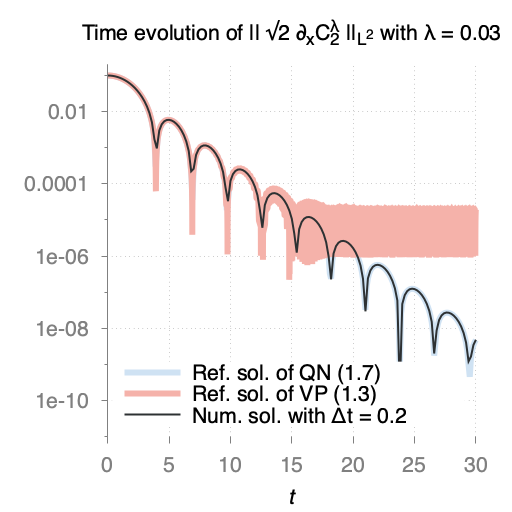}
          \\
          (c) &    (d)
        \end{tabular}
        \caption{{\bf Smooth perturbation of equilibrium test case:} time evolution of  $\|\sqrt{2} \partial_x C_2^\lambda\|_{L^2}$ in logarithmic
          scale  with $(a)$ $\lambda= 1$;  $(b)$ $\lambda =
        3. \;10^{-1}$;  $(c)$ $\lambda = 10^{-1}$ and $(d)$ $\lambda = 3. 10^{-2}$. Red (thick) curve: reference solution of the Vlasov-Poisson system \eqref{vpfp:0} computed on a fine mesh; blue (thick) curve: reference solution of the quasineutral model \eqref{vpfp:neutral}; black (thin) curve: numerical solution with $\Delta t=0.2$.}
      \label{fig:21}
    \end{figure}

The reference solution is obtained using the second order scheme provided in Section \ref{sec:3.4} with a fine mesh, $N_H=2000$, $N_x=1001$ and $\Delta t = 10^{-3}$, which permits to describe both the fast and the slow scales for different values of the Debye length, namely $\lambda\in \{1,\, 0.3,\, 0.1,\, 0.03\}$. The time
evolution of the $L^2$ norm of the electric field $E^\lambda$ and the slow component $E^\lambda_{\rm slow} = \sqrt 2 \partial_x C_2^\lambda$ are shown in 
Figure \ref{fig:20} and in Figures \ref{fig:21} respectively (red curves). When $\lambda=1$,
a damping of the electric field is observed as opposite to the case $\lambda \ll 1$, where the electric field oscillates with a frequency inversely proportional to the Debye length $\lambda$. At the same time, the slow component
$E_{\rm slow}^\lambda$ is first strongly damped with slow oscillations independent of $\lambda$, and then, when its amplitude reaches the order of
$\delta\lambda^2$, the damping stops and it starts oscillating rapidly with a frequency inversely
proportional to $\lambda$. The observed results stay in good agreement with the analytical investigations of \cite{Grenier96, HK16, Bobylev_Potapenko19}.

\begin{figure}[ht!]
 \centering
	\begin{tabular}{cc}
          \includegraphics[width=7.5cm]{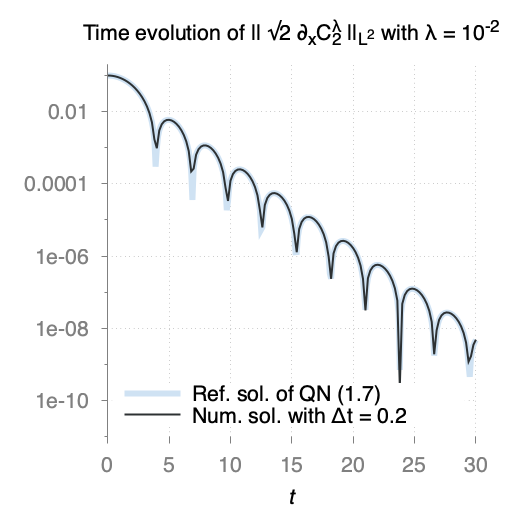}&
          \includegraphics[width=7.5cm]{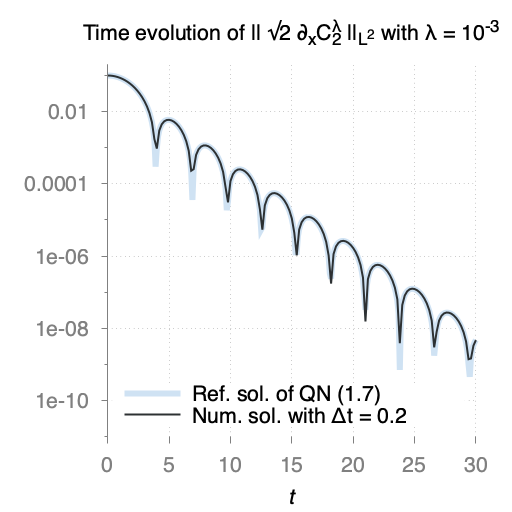}
                                                           \\
          (a) &    (b)
         \end{tabular}
        \caption{{\bf Smooth perturbation of equilibrium test case:} time evolution of  $\|\sqrt{2} \partial_x C_2^\lambda\|_{L^2}$ in logarithmic
          scale  with $(a)$ $\lambda= 10^{-2}$ and $(b)$ $\lambda =
        10^{-3}$. Blue (thick) curve: reference solution of the quasineutral model \eqref{vpfp:neutral}; black (thin) curve: numerical solution with $\Delta t=0.2$.}
      \label{fig:22}
    \end{figure}

To conclude, we also measure the accuracy and the robustness of the second order scheme using a coarser mesh,  $\Delta t=0.2$,  $N_x=101$, $N_H=400$ and for a wide range of values for the Debye length, $\lambda\in\{1,\,3\cdot 10^{-1},\,10^{-1},\,3\cdot 10^{-2}, 10^{-2},10^{-3}\}$. Figures \ref{fig:20}, \ref{fig:21} and \ref{fig:22} display the time evolution of the $L^2$ norm of the approximate electric field $E^{n}$ and its slow counterpart $E_{\rm slow}^{n}\,=\, \sqrt{2}\, \partial_h C^{n}_{2}$ (thin black curves).\\
When $\lambda$ and $\Delta t$ are of the same order, the discrete system is not stiff and a coarse mesh allows to describe precisely the time
oscillations as seen in Figures \ref{fig:20} and~\ref{fig:21} (a)-(b). When~$\lambda$ becomes smaller, the time step $\Delta t=0.2$ is too large to provide a good approximation of
the fast scales as one can notice in Figures~\ref{fig:20} and~\ref{fig:21} (c)-(d). However, we observe on Figure~\ref{fig:21}-(c) that when the norm of $E^\lambda_{\rm slow}$ is larger than the threshold $\delta\,\lambda^2$, the numerical scheme is still able to describe correctly the slow scale dynamics of $E_{\rm slow}^\lambda$. When
$\lambda = 3\cdot 10^{-2}$ with a fixed with $\Delta t=0.2$, Figures \ref{fig:20} and \ref{fig:21} -(d) show that the numerical method eliminates the fast physical oscillations and projects the solution onto the quasineutral slow dynamics limit. In this setting, and for even smaller values of $\lambda$ (see Figure \ref{fig:22} for $\lambda = 10^{-2}$ and $10^{-3}$), the numerical scheme captures
the slow dynamics precisely with a large time step, which illustrates its
asymptotic preserving property. In particular, the black curve in Figure \ref{fig:20}-(d) illustrates Proposition \ref{prop:discrete:QN:limit}, which ensures that for a fixed $\Delta t$, the numerical scheme damps the oscillatory component of the field $E^\lambda$ in the limit $\lambda \rightarrow 0$.

In summary, the numerical results definitively show the
stability and consistency of the numerical scheme
\eqref{discrete:step1}-\eqref{discrete:step2} in the quasineutral
regime,  without a prohibitive computational cost. A trade-off inherent to this approach is the filtering of high-frequency oscillations for coarser resolution in time.

\subsection{Oscillatory perturbation from equilibrium}\label{sec:4.3}
We now consider an oscillatory initial distribution by modifying the
previous one. More
precisely, we choose $f_{\rm in}$ as
\begin{equation}
  \label{ini:cond:osc}
	\ds f_{\rm in}(x,v)\,=\,
	\ds
	\frac{1}{\sqrt{2\pi\,T_{\rm in}(x)}}\, \left( 1 + \delta
          \cos(k_x x)\sin(3\pi v)\right) \,\exp\left(-\frac{|v|^2}{2\,T_{\rm in}(x)}\right)\,,
                  \end{equation}
                 with  $\delta=0.05$, $x\in[-10,10]$, $k_x =
                 \pi/10$ and $T_{\rm in}$ given by
                  $$
T_{\rm in}(x) = 1+\delta\,\cos\left(k_x x\right).
$$ 
As for the previous test case, we have $\rho_{\rm in}(x) = 1$ and, since
  $$
  \int_\R \frac{v\sin(3\pi v)}{\sqrt{2\pi\,T}}\,e^{-|v|^2/(2T)}\,\dD v=3\pi\,T\,e^{-9\pi^2T/2},
  $$
the initial flux satisfies $|j_{\rm in}(x)|\leq 10^{-15}$, hence it is negligible, while high order moments are not spatially homogeneous. This configuration is such that the initial perturbation induces oscillations in the velocity space and, in this situation, one does not expect a fast convergence to zero of the slow component $E_{\rm slow}^\lambda$. In fact, even in the quasineutral regime, velocity oscillations produce a sort of "echo" in the solution, so that the electric field exhibits slowly decaying oscillations with respect to time.

\begin{figure}[ht!]
 \centering
	\begin{tabular}{cc}
          \includegraphics[width=7.5cm]{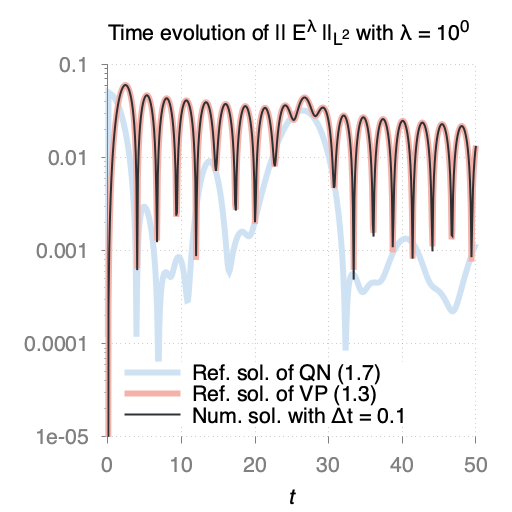}&
          \includegraphics[width=7.5cm]{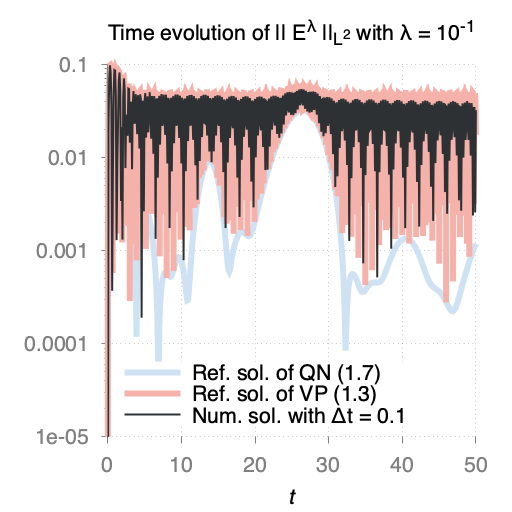}
                                                           \\
          (a) &    (b)
          \\
          \includegraphics[width=7.5cm]{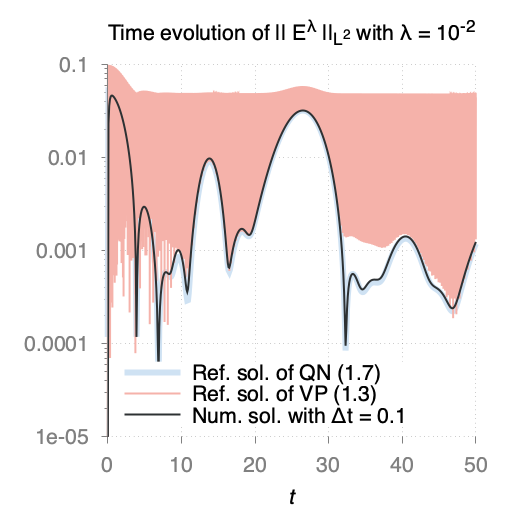}&
          \includegraphics[width=7.5cm]{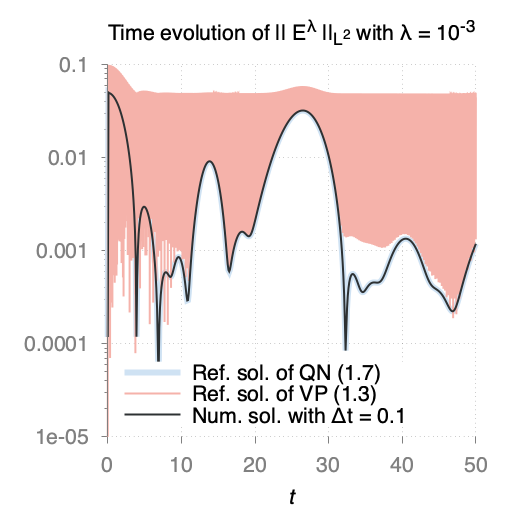}
          \\
          (c) &    (d)
        \end{tabular}
        \caption{{\bf Oscillatory perturbation from equilibrium test case:} time evolution of  $\|E^\lambda\|_{L^2}$ in logarithmic
          scale  with $(a)$ $\lambda= 1$;  $(b)$ $\lambda =
        10^{-1}$;  $(c)$ $\lambda = 10^{-2}$ and $(d)$ $\lambda = 10^{-3}$. Red (thick) curve: reference solution of the Vlasov-Poisson system \eqref{vpfp:0} computed on a fine mesh; blue (thick) curve: reference solution of the quasineutral model \eqref{vpfp:neutral}; black (thin) curve: numerical solution with $\Delta t=0.1$.}
      \label{fig:30}
    \end{figure}

\begin{figure}[ht!]
 \centering
	\begin{tabular}{cc}
          \includegraphics[width=7.5cm]{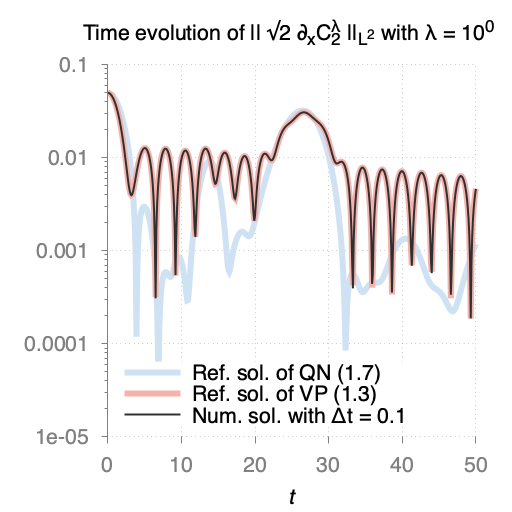}&
          \includegraphics[width=7.5cm]{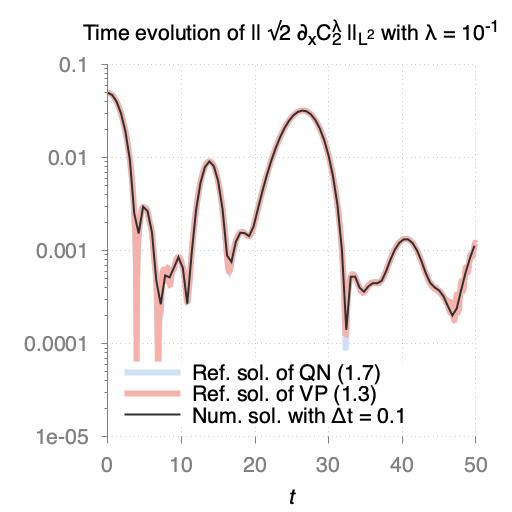}
                                                           \\
          (a) &    (b)
          \\
          \includegraphics[width=7.5cm]{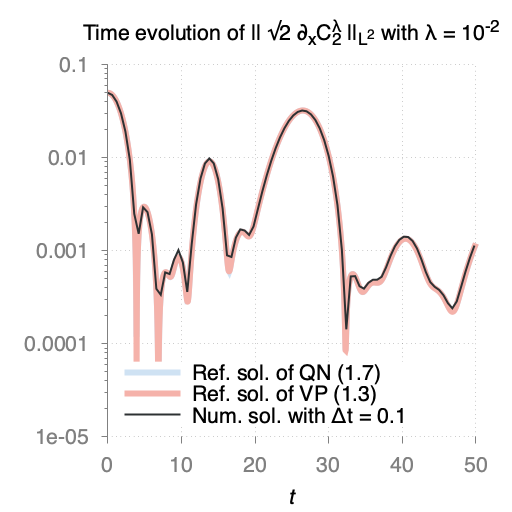}&
          \includegraphics[width=7.5cm]{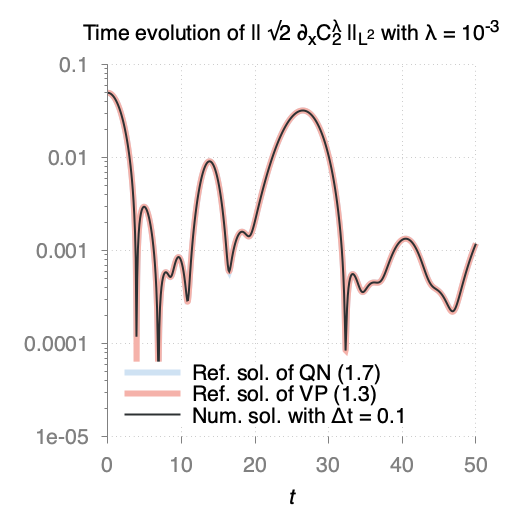}
          \\
          (c) &    (d)
        \end{tabular}
        \caption{{\bf Oscillatory perturbation from equilibrium test case: } time evolution of  $\|\sqrt{2} \partial_x C_2^\lambda\|_{L^2}$ in logarithmic
          value with $(a)$ $\lambda= 1$;  $(b)$ $\lambda =
        10^{-1}$;  $(c)$ $\lambda = 10^{-2}$ and $(d)$ $\lambda = 10^{-3}$. Red (thick) curve: reference solution of the Vlasov-Poisson system \eqref{vpfp:0} computed on a fine mesh; blue (thick) curve: reference solution of the quasineutral model \eqref{vpfp:neutral}; black (thin) curve: numerical solution with $\Delta t=0.1$.}
      \label{fig:31}
    \end{figure}

The reference solution is obtained using a fine mesh with $N_H=1200$, $N_x=501$, $\Delta t =
10^{-3}$ and the second order in time scheme presented in Section \ref{sec:3.4}. This permits to capture both fast and
slow scales phenomena for different values of the Debye length: $\lambda\in\{ 1,\, 10^{-1},
\,10^{-2},\,10^{-3}\}$. In Figure \ref{fig:30}, the time
evolution of the $L^2$ norm of the electric field $E^\lambda$ is shown while in Figure \ref{fig:31}, we  report the evolution of the slow component $E^\lambda_{\rm slow} =
\sqrt 2 \partial_x C_2^\lambda$ (red curves). When $\lambda=1$, the $L^2$ norm of the electric field is first damped and oscillates slowly, then, when $t\simeq 25$, the ``echo'' is observed. When $\lambda \ll 1$, the electric field $E^\lambda$ starts to oscillate with a frequency inversely proportional to $\lambda$, while the slow component~$E_{\rm slow}^\lambda=\sqrt 2\,\partial_xC_2^\lambda$  is first strongly damped, then its amplitude increases again to reach a maximum around time $t\simeq 12.5$ and $t\simeq 25$ (see Figure \ref{fig:31}). The fast oscillations of $E_{\rm slow}^\lambda$ have an amplitude of the order of $\delta \lambda^2$ (\textit{cf.} Section \ref{sec:4.2}) and they are not visible on Figure \ref{fig:31} since the amplitude of $\|\sqrt 2\,\partial_xC_2^\lambda\|_{L^2}$ is above this threshold.

\begin{figure}[ht!]
 \centering
	\begin{tabular}{cc}
          \includegraphics[width=7.5cm]{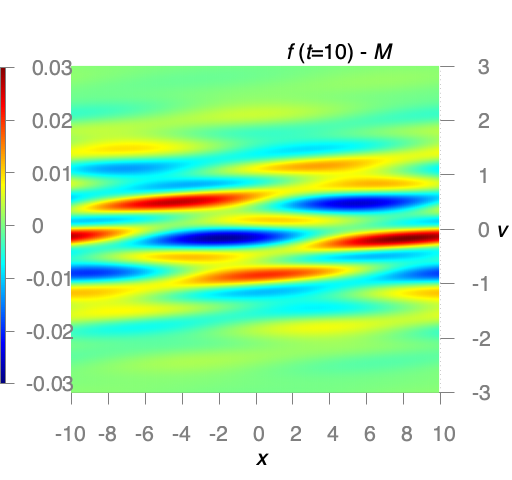}&
          \includegraphics[width=7.5cm]{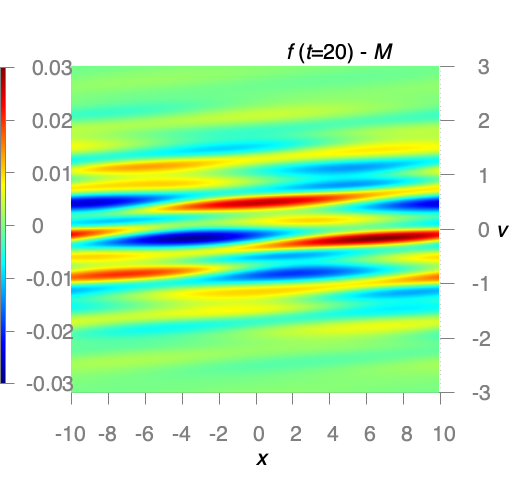}
                                                           \\
          \includegraphics[width=7.5cm]{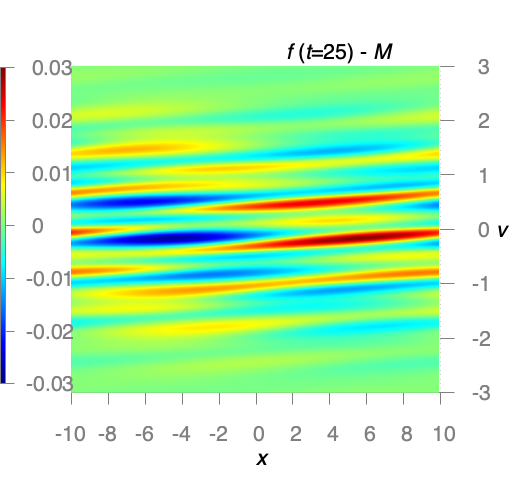}&
          \includegraphics[width=7.5cm]{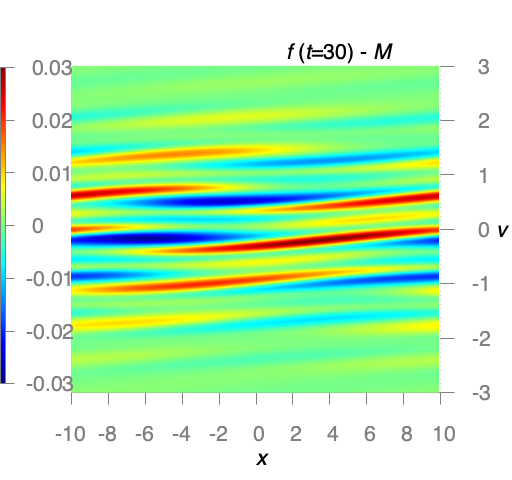}
                                                           \\
          \includegraphics[width=7.5cm]{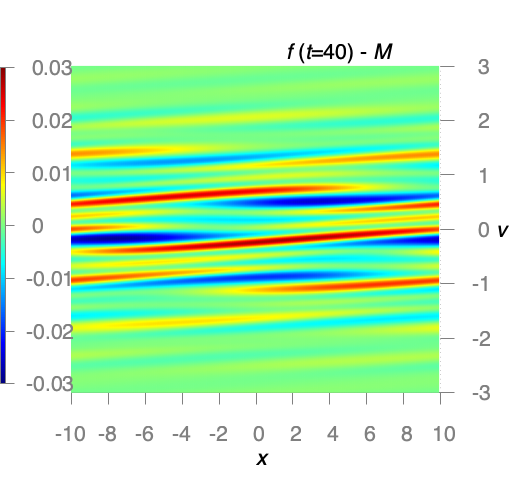}&
          \includegraphics[width=7.5cm]{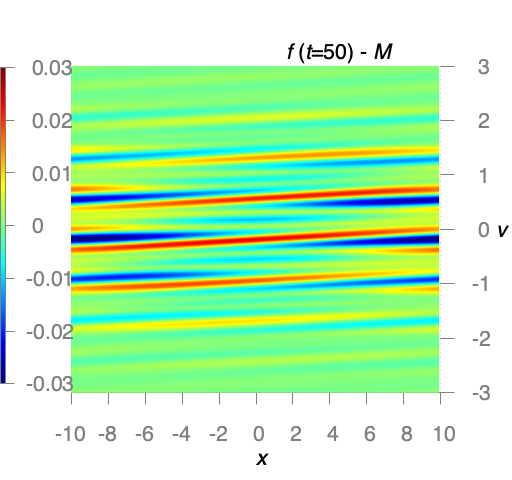}
        \end{tabular}
        \caption{{\bf Oscillatory perturbation from equilibrium test case:}  snapshots of the
          distribution function $f^\lambda- \cM$ at time $t=10$ ;  $20$ ;
          $25$ ;
          $30$ ; $40$ and $50$ for $\lambda=10^{-3}$.}
      \label{fig:32}
    \end{figure}

To investigate the asymptotic preserving property of the scheme we also performed numerical simulations using a coarse mesh, namely $\Delta t=0.1$, $N_x=101$ and $N_H=1200$ for a wide variety of
$\lambda \ll 1$. Notice that, in this case, a large number of Hermite modes is
required to adequately describe oscillations in velocity. The numerical results are reported in Figures \ref{fig:30} and \ref{fig:31} using black curves and compared with the reference
solution. For $\lambda=1$ and  $10^{-1}$, the coarse mesh allows to describe precisely the time
oscillations (Figures~\ref{fig:30} and~\ref{fig:31} (a)-(b)). When $\lambda$ becomes smaller, the
time step $\Delta t=0.1$ is too large to provide an approximation of the fast scales,  but  the slow
scale, corresponding to the quasineutral asymptotic model, is still approximated precisely
(see Figures \ref{fig:30}-(c) and (d) and \ref{fig:31}-(c) and (d)). In particular, Figures \ref{fig:30}-(c) and (d) illustrate the asymptotic preserving property in Proposition \ref{prop:discrete:QN:limit}, since they show that the numerical method projects the approximate electric field $E^n$ onto its quasineutral limit as $\lambda\rightarrow 0$, for a fixed $\Delta t$.
Finally, in Figure \ref{fig:32}, we present several time snapshots of the perturbation $f^\lambda- \cM$, where $\cM$ is the homogeneous Maxwellian
equilibrium with $\lambda=10^{-3}$. When $t=25$,
we observe the oscillation of $f^\lambda$ around the equilibrium corresponding to the ``echo'' of the electric field previously discussed. 

\subsection{Two-stream instability}
\label{sec:4.4}
In this last example, we consider the two stream instability problem with initial distribution
\beq
\label{2stream}
	\ds f_{\rm in}(x,v)\,=\,
	\frac{1}{6\sqrt{2\pi\,T_{\rm in}(x)}}\, \left(1+\frac{5\,|v|^2}{T_{\rm in}(x)}\right)\, \exp\left(
          -\frac{|v|^2}{2T_{\rm in}(x)}\right)\,,
                  \end{equation}
                  and with $T_{\rm in}$ given by
                  $$
T_{\rm in}(x) = 1+\delta\,\cos\left(k_x x\right).
$$
The other parameters are $\delta=0.01$, $x\in[-6,6]$ and $k_x =
\pi/6$. For this initial data, the Vlasov-Poisson system \eqref{vpfp:0} first
develops an instability (called two-stream instability) which then
stabilizes due to nonlinear effects. However, for the quasineutral
limiting model, the situation is more intricate (see for instance
\cite{Bobylev_Potapenko19} for the ill-posedness of  the quasineutral model
with the Landau collisional operator for a two-stream initial distribution or \cite{Grenier96,Han-Kwan_Hauray15}). Thus, the
quasineutral limit cannot be valid for all time intervals.  It is important to notice that obtaining consistent numerical results  for the
Vlasov-Poisson system  over long times is
very difficult when $\lambda$ is small,  this is why we only took $\lambda
=0.04$ as the smallest value. However, even for these intermediate values, we get
meaningful results, which help us to better understand the dynamics. \\

\begin{figure}[ht!]
 \centering
	\begin{tabular}{cc}
          \includegraphics[width=7.5cm]{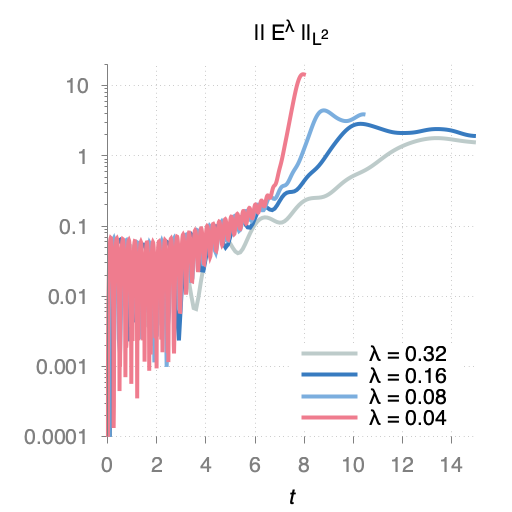}&
          \includegraphics[width=7.5cm]{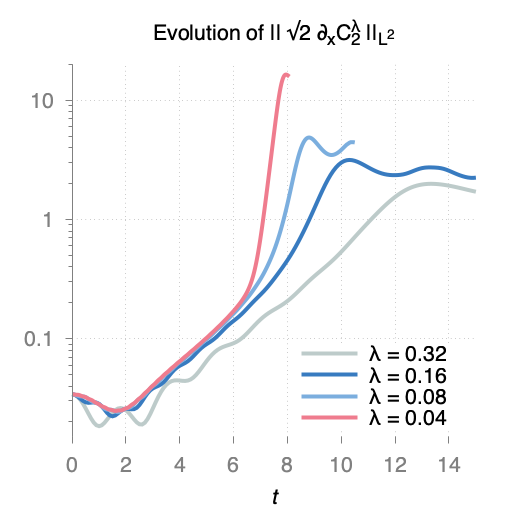}
                                                           \\
          (a) &    (b)
        \end{tabular}
        \caption{{\bf Two-stream instability test case:} time evolution of  (a)
          $\|E^\lambda\|_{L^2}$ and  (b) $\|\sqrt{2} \partial_x
          C_2^\lambda\|_{L^2}$  in logarithmic
          scale.}
      \label{fig:40}
    \end{figure}

\begin{figure}[ht!]
 \centering
	\begin{tabular}{cc}
          \includegraphics[width=7.5cm]{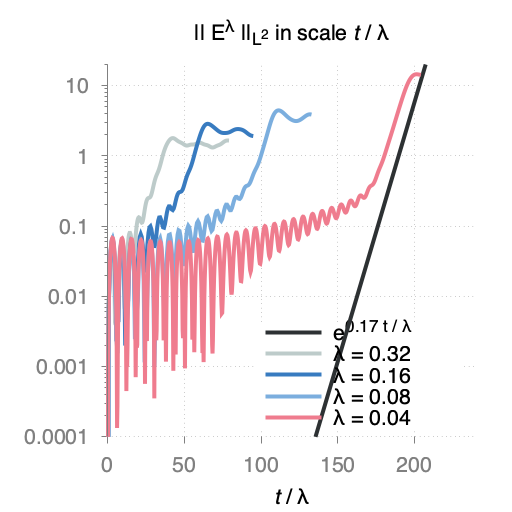}&
          \includegraphics[width=7.5cm]{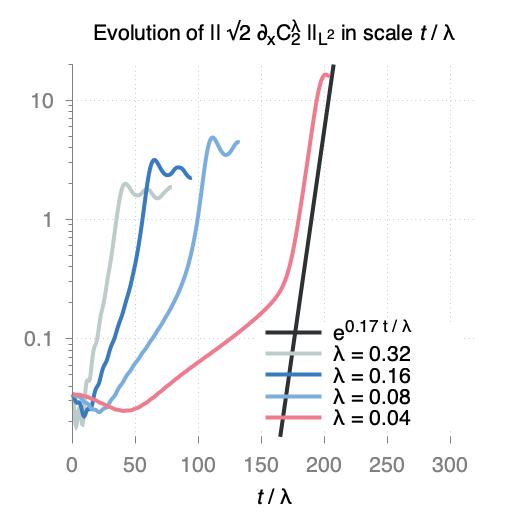}
                                                           \\
          (a) &    (b)
        \end{tabular}
        \caption{{\bf Two-stream instability test case:} evolution of  (a)
          $\|E^\lambda\|_{L^2}$ and  (b) $\|\sqrt{2} \partial_x
          C_2^\lambda\|_{L^2}$  in logarithmic
          scale with respect to $t/\lambda$.}
      \label{fig:41}
    \end{figure}

A reference solution is obtained with the second order scheme of Section \ref{sec:3.4} on a fine mesh with $N_H=2000$, $N_x=2001$ and $\Delta t =10^{-4}$. This permits, as before, to observe both fast and slow scales of the solution for different values of $\lambda$. The time evolution of the $L^2$ norm of the electric field $E^\lambda$ and of the slow component $E^\lambda_{\rm slow} =\sqrt 2 \partial_x C_2^\lambda$ are again computed and shown in Figure \ref{fig:40} for different values of $\lambda$. We observe that the norm of $E^\lambda$ on the time interval $[0,6]$ strongly oscillates with a frequency inversely proportional to $\lambda$ and grows exponentially with a rate independent of $\lambda$. For larger times, the instability rate increases as $\lambda$ decreases: a second, much faster, exponential growth occurs, with a rate proportional to $\lambda^{-1}$. The evolution of $E_{\rm slow}^\lambda$ follows a similar growth dynamics, without exhibiting an oscillatory behavior. In the situation depicted, the rise of the instability is due to the slow part of the electric field, illustrating that the limiting system is not globally well-posed. Figure \ref{fig:41} presents the evolution of the same quantities as of Figure \ref{fig:40} with a scale of $t/\lambda$. This permits to understand that during this second phase, for $t\geq 6$, the instability is of the order of $O(e^{\kappa t/\lambda})$, with $\kappa\simeq 0.17$ for $\lambda=0.04$. These results show that when addressing the full scales, the numerical scheme allows for a better understanding of the instability phenomenon and the fact that the quasineutral limit is not valid on large time intervals for this initial data, which is consistent with the theoretical results of E. Grenier \cite{Grenier96} and D. Han-Kwan and M. Hauray \cite{Han-Kwan_Hauray15}.

\begin{figure}[ht!]
 \centering
	\begin{tabular}{cc}
          \includegraphics[width=7.5cm]{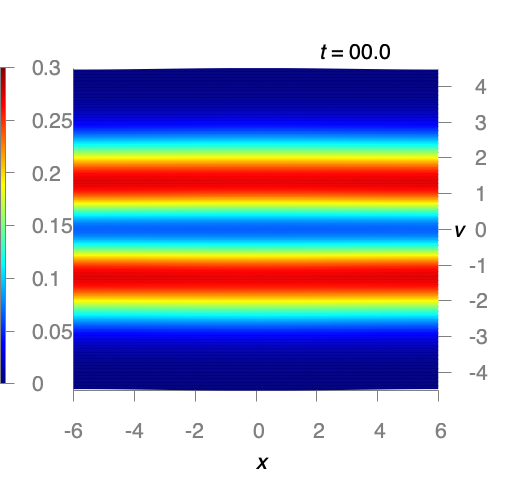}&
          \includegraphics[width=7.5cm]{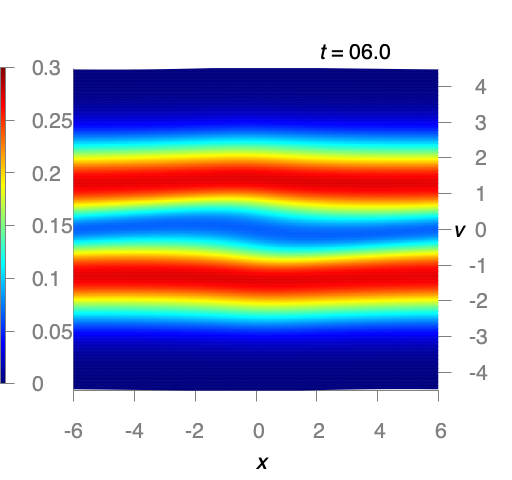}
                                                           \\
          \includegraphics[width=7.5cm]{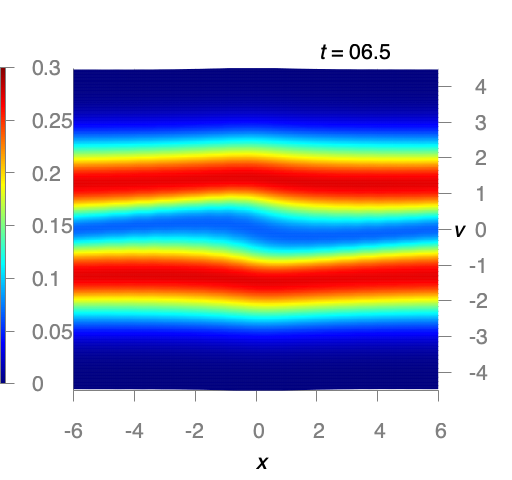}&
          \includegraphics[width=7.5cm]{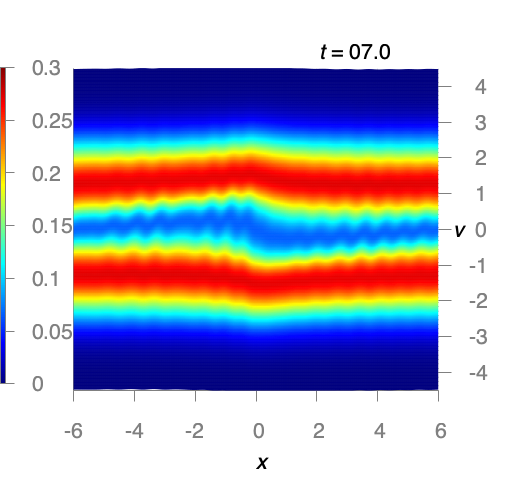}
                                                           \\
          \includegraphics[width=7.5cm]{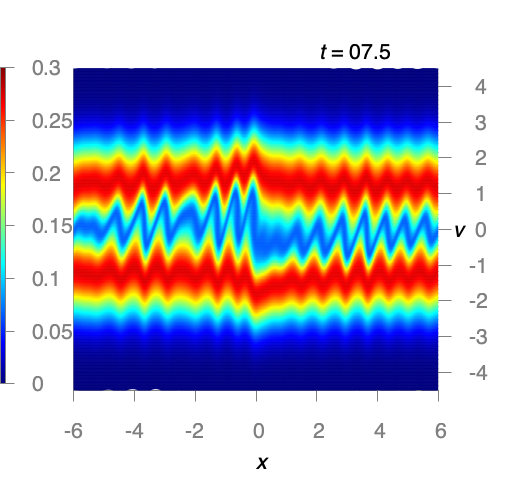}&
          \includegraphics[width=7.5cm]{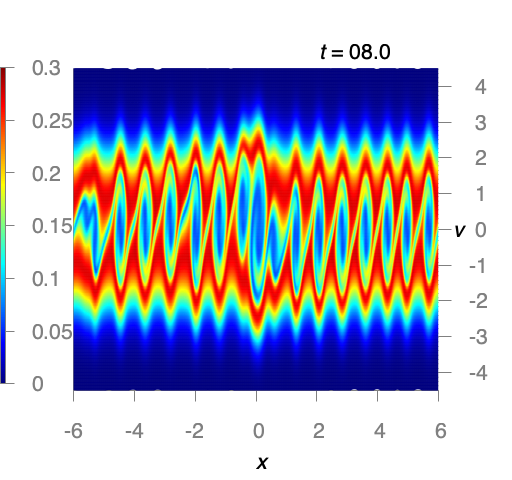}
        \end{tabular}
        \caption{{\bf Two-stream instability test case:}  snapshots of the
          distribution function $f^\lambda$ at time $t=0$ ;  $6$ ;
          $6.5$ ;
          $7$ ; $7.5$ and $8$ for $\lambda=0.04$.}
      \label{fig:42}
    \end{figure}

Finally, in Figure \ref{fig:42}, we present the snapshots of the distribution
function $f^\lambda$ with $\lambda=0.04$ at different time, when the
instability occurs. At time $t=8$, we observe several vortices in
phase space $(x,v)$, which illustrates the  complexity of the dynamics
when $\lambda$ is small.    

\section{Conclusion and perspectives}
\label{sec5}
\setcounter{equation}{0}
\setcounter{figure}{0}
\setcounter{table}{0}

In this work, we proposed a new numerical scheme for the
Vlasov-Poisson system able to overcome the strong
restrictions due to the fast scale dynamics related to
quasineutrality in plasmas. This discretization is quite different
from those previously proposed in the literature, as it does not seek
to reformulate the Poisson equation into an equivalent harmonic-type oscillator, but rather to discretize the Poisson equation simultaneously
with the equations of the first moments. A significant
difficulty in this study is that the quasineutral limit is not always
valid for all initial conditions and all time intervals. On the one hand,  our approach expands the distribution function in velocity space using Hermite functions, this enables to study the quasineutral limit by explicitly separating the oscillatory part from the rest of the solution and permits to obtain error
estimates with respect to $\lambda$ both at the continuous and the discrete levels.  On the other hand, we propose a time-splitting method,
where the first step involves solving the linearized system around the
Maxwellian stationary state and the second step solves the nonlinear remaining part of the system.

In a second part, we have performed several numerical simulations for stable and unstable initial conditions. We have recovered the theoretical convergence order estimates for well-prepared initial data. In these situations and when the time step is very large compared to $\lambda$, the numerical scheme filters out the rapid oscillations and captures the asymptotic limit. Moreover, our numerical simulations indicate two different sources of instability. On the one hand, when the quasineutral limit is not valid on large time intervals, as for the two-stream instability, the instability is physical and it is correctly captured by the scheme when the fast scales are resolved. On the other hand, in the critical case $\alpha=1$, when the time step is of the same order as $\lambda$, the numerical solution may become unstable. Understanding the quasineutral limit remains an open problem within the framework of
kinetic theory, our numerical results are consistent with those available in the mathematical analysis literature \cite{Grenier96, Han-Kwan_Hauray15, Bobylev_Potapenko19} and we expect to provide novel insights that could guide the demonstration of the required mathematical results.

Concerning the analysis of the numerical scheme, it appears crucial to gain a better understanding of the stability of the scheme with respect to the discretization parameters and the Debye length.  Furthermore, extending the proposed method to the multi-dimensional case \cite{blaustein24}  represents an important step for realistic applications.

%
\section*{Acknowledgement}

Francis Filbet and Marie-Hélène Vignal are partially funded by the ANR Project Muffin
(ANR-19-CE46-0004). Alain Blaustein  acknowledges the support of the CDP C2EMPI, as well as the French State under the France-2030 programme, the University of Lille, the Initiative of Excellence of the University of Lille, the European Metropolis of Lille for their funding and support of the R-CDP-24-004-C2EMPI project.

\appendix
\section{Proof of Proposition \ref{Prop:dev:asymptotique}}\label{App:prop:dev:asympt}

We fix some time $t\in[0,T]$ and prove that the right hand side in the Duhamel formula \eqref{duhamel} is of order $\lambda^{1-\alpha}$. To this aim, we distinguish linear and nonlinear terms in $\bS^\lambda$. Then, to identify the nonlinear terms and their scaling in $\lambda$, we replace $\partial_t C_2^\lambda$ according to \eqref{Hermite:D} with $k=2$ and $C_0^\lambda$ according to the second line of \eqref{Hermite:D} in $\bS^\lambda$, we obtain
\begin{equation}\label{Duhamel:norm}
	\bU^\lambda(t)
	\,=\,-
	\int_{0}^{t}
	\exp{\left(\frac{s-t}{\lambda}J\right)}\,(\bN^\lambda+\bL^\lambda)(s)\,\dD s
	\,,
\end{equation}
where $\bN^\lambda$ and $\bL^\lambda$ are given as follows
\[
\bN^\lambda\,=\,
2\,T_0^{\,1/2}\partial_x
\begin{pmatrix}
	\ds\,C_1^\lambda E^\lambda  \\[0.5em]
	\ds 0
\end{pmatrix}
\,-\,\frac{\lambda}{2}\,\partial_x
\begin{pmatrix}
	\ds0 \\[0.5em]
	\ds |E^\lambda|^2
\end{pmatrix};\quad\quad
\bL^\lambda\,=\,
-\,T_0^{\,3/2}\partial^2_x 
\begin{pmatrix}
	\ds\,
	\sqrt{6}\,C_3^\lambda + 2\,C_1^\lambda  \\[0.5em]
	\ds 0
\end{pmatrix}\,+\,
\lambda \,T_0\,
\partial^2_x
\begin{pmatrix}
	\ds0 \\[0.5em]
	\ds 
	E^\lambda
\end{pmatrix}.
\]
To control the contribution of $\bL^\lambda$, we integrate by parts in time, which is a standard method for oscillating integrals, and use assumptions \eqref{compatibility}-\eqref{uniform:reg} to control the time derivatives of $C_1^\lambda$ and $C_3^\lambda$. The nonlinear contribution $\bN^\lambda$ is more intricate since $|E^\lambda|^2$ may behave like $O(\lambda^{-2\alpha})$. However, we show that the singular terms are solely due to fast oscillations, whose time integral are small as $\lambda\rightarrow 0$.\\
First, we replace $C_1^\lambda$ in $\bN^\lambda$ according to the Amp\`ere equation \eqref{Amp:flux:H:0}, which yields
\[
\bN^\lambda(s)\,=\,-\,
\lambda^2\,\partial_{t}\partial_{x}
\begin{pmatrix}
	\ds |E^\lambda(s)|^2\\[0.5em]
	\ds 0
\end{pmatrix}
\,-\,\frac{\lambda}{2}\,\partial_x
\begin{pmatrix}
	\ds0 \\[0.5em]
	\ds |E^\lambda(s)|^2
\end{pmatrix}\,,
\]
where we also applied the relation $2E^\lambda\partial_{t}E^\lambda\,=\, \partial_{t} |E^\lambda|^2$. Multiplying the latter formula by $\exp{\left(sJ/\lambda\right)}$, applying Leibniz rule for products and using the definition of $J$ below \eqref{vect:Amp:flux}, yields
\begin{equation}\label{N:lambda:technical}
	\exp{\left(\frac{s}{\lambda}J\right)}\cdot\bN^\lambda(s)\,=\,
	-\lambda^2\,\partial_s
	\left(\exp{\left(\frac{s}{\lambda}J\right)}\cdot\partial_{x}
	\begin{pmatrix}
		\ds |E^\lambda(s)|^2\\[0.5em]
		\ds 0
	\end{pmatrix}\right)
	-\,\frac{3\lambda}{2}
	\exp{\left(\frac{s}{\lambda}J\right)}\cdot\partial_{x}
	\begin{pmatrix}
		\ds 0\\[0.5em]
		\ds |E^\lambda(s)|^2
	\end{pmatrix}
	.
\end{equation}
The main difficulty comes from the second term in the latter right
hand side.  Hence, we estimate it by  computing explicitly the contribution of
the singular terms in the expansion of $E^\lambda$ as $\lambda
\rightarrow 0$. More precisely, we substitute $E^\lambda$ thanks to
the relation
$$
E^\lambda\,=\,\bU_1^{\lambda}+E^\lambda_{\textrm{slow}}+E^\lambda_{\textrm{osc}}\,,
$$
where $\bU^\lambda_1$ denotes the first component of $\bU^\lambda$ and where  $E^\lambda_{\textrm{osc}}$, given by \eqref{E:C_1:osc}, contains the singular terms in the expansion of $E^\lambda$. Hence, we have
\[
\frac{1}{2} \,\partial_x |E^\lambda|^2\,=\,
\bU_1^{\lambda}\, \partial_x
E^\lambda
+
E^\lambda_{\textrm{osc}}\,\partial_x\bU_1^{\lambda}
+
\frac{1}{2} \partial_x |E^\lambda_{\textrm{osc}}|^2
+
E^\lambda_{\textrm{osc}}\, \partial_x E^\lambda_{\textrm{slow}}
+
E^\lambda_{\textrm{slow}}\, \partial_x E^\lambda\,,
\]
We replace $|E^\lambda|^2$ according to the latter relation in the second term on the right hand side of \eqref{N:lambda:technical} and obtain
  \begin{equation}\label{reform:N}
	\exp{\left(\frac{s}{\lambda}J\right)}\cdot\bN^\lambda\,=\,
	-\lambda^2 \,\partial_s
	\left(\exp{\left(\frac{s}{\lambda}J\right)}\cdot
	\bN^\lambda_1\right)
	\,-\,
	\exp{\left(\frac{s}{\lambda}J\right)}\cdot\left(A_1\bU^\lambda_1+A_2\,\partial_x\bU^\lambda_1+\	\lambda \,\bN^\lambda_2 \right),
\end{equation}
where $(\bN^\lambda_1,\bN^\lambda_2)$ and $(A_1,A_2)$ are given as follows
\[
\bN_1^\lambda(s)\,=\,
\partial_{x}
\begin{pmatrix}
	\ds |E^\lambda(s)|^2\\[0.5em]
	\ds 0
\end{pmatrix}
;\quad\quad
A_1(s)\,=\,
3\,\lambda\,
\begin{pmatrix}
	\ds0 \\[0.5em]
	\ds \partial_xE^\lambda(s)
\end{pmatrix}\;;\quad\quad
A_2(s)\,=\,
3\,\lambda\,
\begin{pmatrix}
	\ds0 \\[0.5em]
	\ds E^\lambda_{\textrm{osc}}(s)
\end{pmatrix}\;
;
\]
whereas 
\[
\bN^\lambda_2(s)\,=\,3\,
\begin{pmatrix}
	\ds0 \\[0.5em]
	\ds \frac{1}{2} \partial_x |E^\lambda_{\textrm{osc}}(s)|^2
	+
	E^\lambda_{\textrm{osc}}(s)\,	
	\partial_x E^\lambda_{\textrm{slow}}(s)
	+
	E^\lambda_{\textrm{slow}}(s)  \partial_x E^\lambda(s)
\end{pmatrix}(s)\,
.
\]
Then, we substitute $\bN^\lambda$ in \eqref{Duhamel:norm} according to \eqref{reform:N} and obtain
\begin{equation}\label{duham:reform}
	\bU^\lambda(t)
	=\,
	\int_{0}^{t}
	\exp{\left(\frac{s-t}{\lambda}J\right)}\cdot
	\left(
	A_1(s)\bU_1^\lambda(s)
	+
	A_2(s)\partial_x\bU_1^\lambda(s)\right)
	\dD s
	+B(t)\,,
\end{equation}
where $A_1$ and $A_2$ are given below \eqref{reform:N} and $B$ gathers all the other terms, that is, 
\[
B(t)
=\,
\int_{0}^{t}
\exp{\left(\frac{s-t}{\lambda}J\right)}\cdot(\lambda\,\bN^\lambda_2(s)-\bL^\lambda(s))\,\dD s
+
\lambda^2\,\bN^\lambda_1(t)
- \lambda^2\exp{\left(-\frac{t}{\lambda}J\right)}\cdot\bN^\lambda_1(0)\,,
\]
where $\bL^\lambda$ is given in \eqref{Duhamel:norm} whereas
$\bN_1^\lambda$, $\bN_2^\lambda$ are given in \eqref{reform:N}.

Our strategy is to prove that $A_1$, $A_2$ and $B$ are of order
$O(\lambda^{1-\alpha})$ and iterate the Duhamel formula
\eqref{duham:reform} to gain powers of $\lambda$.
\\

\paragraph{\bf \underline{Step 1.} Upper bound of $A_1$, $A_2$ and $B$.}
Let us  estimate $A_1$, $A_2$ and $B$, given in \eqref{duham:reform}, starting with $B$, which contains the main difficulty, due to the contribution of $\bN^\lambda_2$. To estimate $B$, we reformulate the first term of $\bL^\lambda$, defined in \eqref{Duhamel:norm}, which is given for $s\geq 0$ by the function
\[
\bL^\lambda_1(s)\,=\,
-\,T_0^{\,3/2}\partial^2_x 
\begin{pmatrix}
	\ds\,
	\sqrt{6}\,C_3^\lambda + 2\,C_1^\lambda  \\[0.5em]
	\ds 0
\end{pmatrix}(s,\cdot).
\]
We multiply $\bL^\lambda_1(s)$ by $\exp{\left(sJ/\lambda\right)}$ and apply Leibniz rule for products, which yields
\begin{equation*}
	\exp{\left(\frac{s}{\lambda}J\right)}\cdot\bL^\lambda_1(s)\,=\,
	\lambda \,J^{-1}\,\partial_s
	\left(\exp{\left(\frac{s}{\lambda}J\right)}\cdot
	\bL^\lambda_1(s)\right)
	-
	\lambda \,J^{-1}
	\exp{\left(\frac{s}{\lambda}J\right)}\cdot
	\partial_s\,
	\bL^\lambda_1(s)\,.
\end{equation*}
Next, since $J^{-1} = -J$ and since $J$ and $\exp{\left(\frac{s}{\lambda}J\right)}$ commute, we obtain
  \begin{eqnarray*}
	\exp{\left(\frac{s}{\lambda}J\right)}\cdot\bL^\lambda_1(s) &=&
	-\,\lambda \,\partial_s
	\left(\exp{\left(\frac{s}{\lambda}J\right)}\cdot T_0^{\,3/2}
	\partial^2_x 
	\begin{pmatrix}
		\ds 0 \\[0.5em]
		\ds\,
		\sqrt{6}\,C_3^\lambda + 2\,C_1^\lambda 
	\end{pmatrix}\right)(s)
    \\
    && +\,
	\lambda 
	\exp{\left(\frac{s}{\lambda}J\right)}\cdot\,T_0^{\,3/2}
	\partial_{s}\partial^2_x 
	\begin{pmatrix}
		\ds 0 \\[0.5em]
		\ds\,
		\sqrt{6}\,C_3^\lambda + 2\,C_1^\lambda 
	\end{pmatrix}(s).
\end{eqnarray*}
According to the latter relation, the product between $\exp{\left(sJ/\lambda\right)}$ and $\bL^\lambda$ satisfies
  \begin{equation}\label{reform:L}
	\exp{\left(\frac{s}{\lambda}J\right)}\cdot\bL^\lambda(s)\,=\,
	-\lambda \,\partial_s
	\left(\exp{\left(\frac{s}{\lambda}J\right)}\cdot
	\bL^\lambda_2(s)\right)
	\,+\,
	\lambda
	\exp{\left(\frac{s}{\lambda}J\right)}\cdot\,
	\bL^\lambda_3(s)\,,
\end{equation}
where $\bL^\lambda_2$ and $\bL^\lambda_3$ are given as follows
  \[
\bL_2^\lambda\,=\,T_0^{\,3/2}
\partial^2_x 
\begin{pmatrix}
	\ds 0 \\[0.5em]
	\ds\,
	\sqrt{6}\,C_3^\lambda + 2\,C_1^\lambda 
\end{pmatrix};\quad\quad
\bL^\lambda_3\,=\,
T_0^{\,3/2}
\partial_{t}
\partial^2_x 
\begin{pmatrix}
	\ds 0 \\[0.5em]
	\ds\,
	\sqrt{6}\,C_3^\lambda + 2\,C_1^\lambda 
\end{pmatrix}\,+\,
T_0
\begin{pmatrix}
	\ds0 \\[0.5em]
	\ds \partial^2_x
	E^\lambda
\end{pmatrix}.
\]
We replace $\bL^\lambda$ according to \eqref{reform:L} in the definition of $B$ in \eqref{duham:reform}, we find 
\[
B(t)
=\,\lambda
\int_{0}^{t}
\exp{\left(\frac{s-t}{\lambda}J\right)}\cdot(\bN^\lambda_2-\bL^\lambda_3)(s)\,\dD s
\,+\,
\lambda\,(\lambda\,\bN^\lambda_1+\bL^\lambda_2)(t)
\,-\, \lambda\,\exp{\left(-\frac{t}{\lambda}J\right)} \cdot(\lambda\,\bN^\lambda_1+ \bL^\lambda_2)(0)\,.
\]
We take the $L^2$ norm of the $l$-th derivative and apply the triangle inequality in the latter relation
\begin{equation}\label{estim:B:0}
\begin{split}
\left\|
\partial_x^l
B(t)\right\|_{L^2}
\leq\;&\ds\lambda
\left\|
\partial_x^l
\int_{0}^{t}
\exp{\left(\frac{s}{\lambda}J\right)}\cdot\bN^\lambda_2(s)\,\dD s\right\|_{L^2}
+
\lambda
\left\|
\partial_x^l
\int_{0}^{t}
\exp{\left(\frac{s}{\lambda}J\right)}\cdot\bL^\lambda_3(s)\,\dD s\right\|_{L^2}\\[0.8em]
+\,&\ds
\lambda^2
\left\|
\partial_x^l\bN^\lambda_1(t)\right\|_{L^2}
+\lambda
\left\|
\partial_x^l\bL^\lambda_2(t)\right\|_{L^2}
+ \lambda^2
\left\|
\partial_x^l
\bN^\lambda_1(0)\right\|_{L^2}
+
\lambda
\left\|
\partial_x^l\bL^\lambda_2(0)\right\|_{L^2}.
\end{split}
\end{equation}
To estimate the first term in the right hand side of \eqref{estim:B:0}, we substitute $\bN_2^\lambda$ according to its definition below \eqref{reform:N} and apply the triangle inequality, which yields
\[
\lambda
\left\|
\partial_x^l
\int_{0}^{t}
\exp{\left(\frac{s}{\lambda}J\right)}\cdot\bN^\lambda_2\dD s\right\|_{L^2}
\,\leq\,\cN_{21}(t)\,+\,\cN_{22}(t)\,,
\]
where 
\begin{equation*}
	\left\{
	\begin{array}{l}
		\ds\cN_{21}(t) 
		\,=\,
		\frac{3\lambda}{2}
		\left\|
		\partial_x^l
		\int_{0}^{t}
		\exp{\left(\frac{s}{\lambda}J\right)}\cdot\partial_x
		\begin{pmatrix}
			\ds0 \\[0.5em]
			\ds |E^\lambda_{\textrm{osc}}|^2
		\end{pmatrix}
		\dD s\,\right\|_{L^2},
		\\[1.5em]
		\ds\cN_{22}(t) 
		\,=\,
		3\lambda
		\left\|
		\partial_x^l
		\int_{0}^{t}
		\exp{\left(\frac{s}{\lambda}J\right)}\cdot
		\begin{pmatrix}
			\ds0 \\[0.5em]
			\ds E^\lambda_{\textrm{osc}}\,
			\sqrt{2}\,T_0\,
			\partial_x^2 C^\lambda_{2}
			+
			\sqrt{2}\,T_0\,
			\partial_x C^\lambda_{2} \partial_x
			E^\lambda
		\end{pmatrix}
		\dD s\,\right\|_{L^2}.
	\end{array}\right.
\end{equation*}
The estimate for $\cN_{21}(t)$ is the most intricate since $E^\lambda_{\textrm{osc}}$, given by \eqref{E:C_1:osc}, is of order $O(\lambda^{-\alpha})$. We compute the time integral in $\cN_{21}(t)$ explicitly and show that the singularity cancels due to the fast oscillations of $E^\lambda_{\textrm{osc}}$. To make these computations more tractable, we identify 
$\exp{\left(sJ/\lambda\right)}$ and $E^\lambda_{\textrm{osc}}$ in $\cN_{21}(t)$ with their complex representation
\[
\exp{\left(\frac{s}{\lambda}J\right)}\,=\,e^{-i\frac{s}{\lambda}}\quad;\quad
\begin{pmatrix}
	\ds0 \\[0.5em]
	\ds |E^\lambda_{\textrm{osc}}|^2
\end{pmatrix}\,=\,i
\left|E^\lambda_{\textrm{osc}}\right|^2\quad;\quad 
E^\lambda_{\textrm{osc}}(s)\,=\,
\frac{ze^{i\frac{s}{\lambda}}+ \bar{z}e^{-i\frac{s}{\lambda}}}{2}\;,
\]
where the time independent complex number $z$ is given by
\beq
\label{def:z}
z \,=\, \left(E^\lambda - \sqrt{2}\,T_0\, \partial_x C^\lambda_{2}\right)(0)
\,+\,i\frac{\sqrt{T_0}}{\lambda}
\,C_1^\lambda (0).
\eeq
We reformulate $\cN_{21}(t)$ according to the latter relations
\[
\cN_{21}(t) 
\,=\,
\frac{3\lambda}{2}
\left\|
\partial_x^{l+1}
\int_{0}^{t}
e^{-i\frac{s}{\lambda}}\,i
\left|\frac{ze^{i\frac{s}{\lambda}}+ \bar{z}e^{-i\frac{s}{\lambda}}}{2}\right|^2
\dD s\,\right\|_{L^2}
\]
and expand the square in the latter expression and compute the time integral explicitly:
\[
\cN_{21}(t) 
\,=\,
\frac{3\lambda^2}{2}
\left\|
\partial_x^{l+1}
\left(
\frac{|z|^2}{2}
\left(1-e^{-i\frac{t}{\lambda}}\right)
+
\frac{z^2}{4}
\left(e^{i\frac{t}{\lambda}}-1\right)
+
\frac{\bar{z}^2}{12}
\left(1-e^{-3i\frac{t}{\lambda}}\right)\right)\right\|_{L^2}.
\]
Then, we apply Leibniz rule to develop the derivatives in $x$ of order $l+1$ and estimate each product thanks to Cauchy-Schwarz inequality, which yields
\[
\cN_{21} (t)
\,\leq\,
C\,\lambda^2
\,\left\|
z\right\|^2_{W^{l+1,4}},
\]
for some constant $C$ depending on $l\geq 0$. We replace $z$ according to its definition \eqref{def:z}  and apply the triangle inequality in the latter relation and obtain
\[
\cN_{21} (t)
\,\leq\,
C\,
\left(
\lambda^2
\left\|
C^\lambda_{2}(0)
\right\|^2_{W^{l+2,4}}
\,+\,
\lambda^2
\left\|
E^{\lambda}(0)
\right\|^2_{W^{l+1,4}}
\,+\,
\left\|
C^{\lambda}_1(0)
\right\|^2_{W^{l+1,4}}
\right).
\]
To estimate $\cN_{22}(t)$, we bound the time integral thanks to Jensen inequality and then take the supremum in time, this yields
\[
\cN_{22}(t) \,\leq\,
3\,\lambda\,
t
\sup_{0\leq s\leq T}\left(
\left\|
\partial_x^{l}
\left(E^\lambda_{\textrm{osc}}(s)\,
\sqrt{2}\,T_0\,
\partial_x^2 C^\lambda_{2}(s)\right)
\right\|_{L^2}
+
\left\|
\partial_x^{l}
\left(\sqrt{2}\,T_0\,
\partial_x C^\lambda_{2}(s) \partial_x
E^\lambda(s)\right)\right\|_{L^2}\right)
\,.
\]
Then, we apply Leibniz rule to develop the derivatives of order $l$ and estimate each product thanks to Cauchy-Schwarz inequality, which yields
\[
\cN_{22}(t) \,\leq\,
C\,\lambda
\sup_{0\leq s\leq T}
\left(
\left\|
C^\lambda_{2}(s)
\right\|_{W^{l+2,4}}
\left\|
E^\lambda_{\textrm{osc}}(s)
\right\|_{W^{l,4}}
+
\left\|
C^\lambda_{2}(s)
\right\|_{W^{l+1,4}}
\left\|
E^{\lambda}(s)
\right\|_{W^{l+1,4}}
\right),
\]
for some constant $C$ depending on $l\geq 0$, on the final time $T$ and the temperature $T_0$. We estimate $E^\lambda_{\textrm{osc}}$ with the triangular inequality, which yields
\[
\cN_{22}(t)\,\leq\,
C\,\lambda
\sup_{0\leq s\leq T}
\left(
		\left\|
C^\lambda_{2}(s)
\right\|_{W^{l+2,4}}
(
\|
E^\lambda(0)
\|_{W^{l,4}}
+
\|
C^\lambda_{2}(0)
\|_{W^{l+1,4}}
+\lambda^{-1}
\|
C^\lambda_{1}(0)
\|_{W^{l,4}}
+
\|
E^{\lambda}(s)
\|_{W^{l+1,4}})
\right),
\]
and then apply Young inequality to estimate each product, with a weight adapted to the scaling of each term. More precisely, setting
  $$
  a:=\sup_{0\leq s\leq T}\|C^\lambda_{2}(s)\|_{W^{l+2,4}},
  $$
  we use that for $0<\lambda<1$,
  $$
  \lambda\, a\, b\,\leq\,\frac{1}{2}\left(\lambda^{1-\alpha}\,a^2\,+\,\lambda^{1+\alpha}\,b^2\right)\,,
  $$
  with $b=\|E^\lambda(0)\|_{W^{l,4}}$ or $b=\|E^\lambda(s)\|_{W^{l+1,4}}$. Then, we also use that for $0<\lambda<1$,
\[
\left\{
\begin{array}{l}
  \ds\lambda\,a\,\|C^\lambda_{2}(0)\|_{W^{l+1,4}}\,\leq\,\lambda\, a^2\,\leq\,\lambda^{1-\alpha}\,a^2\,,
\\[1.1em]
\ds a\,\|C^\lambda_{1}(0)\|_{W^{l,4}}\,\leq\,\frac{1}{2}\left(\lambda^{1-\alpha}\,a^2\,+\,\lambda^{\alpha-1}\,\|C^\lambda_{1}(0)\|_{W^{l,4}}^2\right)\,.
\end{array}\right.
\]
In particular, the contribution of $C^\lambda_{2}(0)$ is absorbed in the term involving the supremum in time of $C^\lambda_2$. We obtain
\[
\cN_{22}(t)\,\leq\,
C
\sup_{0\leq s\leq T}
\left(
\lambda^{1-\alpha}
\left\|
C^\lambda_{2}(s)
\right\|_{W^{l+2,4}}^2
+
\lambda^{\alpha-1}
\left\|
C^\lambda_{1}(0)
\right\|_{W^{l,4}}^2
+
\lambda^{1+\alpha}
\left\|
E^{\lambda}(s)
\right\|_{W^{l+1,4}}^2
\right),
\]
for all $0<\lambda<1$. Gathering our estimates on $\cN_{21}$ and $\cN_{22}$, and using that $\lambda^2\leq \min(\lambda^{1-\alpha},\lambda^{1+\alpha})$ and $1\leq\lambda^{\alpha-1}$ for $0<\lambda<1$ to bound the terms involving $C_2^\lambda(0)$, $E^\lambda(0)$ and $C_1^\lambda(0)$ in the estimate of $\cN_{21}$, we obtain
\begin{equation}\label{estim:N4}
	\begin{split}
		\lambda
		\left\|
		\partial_x^l
		\int_{0}^{t}
		\exp{\left(\frac{s}{\lambda}J\right)} \cdot\bN^\lambda_2\dD s\right\|_{L^2}
		\,\leq\;
		&C\left(\lambda^{\alpha-1}
		\left\|
		C^\lambda_{1}(0)
		\right\|_{W^{l+1,4}}^2
		+
		\lambda^{1+\alpha}
		\sup_{0\leq s\leq T}
		\left\|
		E^{\lambda}(s)
		\right\|_{W^{l+1,4}}^2
		\right)\\
		\,& +\,C\,\lambda^{1-\alpha}
		\sup_{0\leq s\leq T}
		\left\|
		C^\lambda_{2}(s)
		\right\|_{W^{l+2,4}}^2,
	\end{split}
\end{equation}
for all $0<\lambda<1$ and for some constant $C$ depending on $l\geq 0$, on the final time $T$ and the temperature $T_0$.\\
We now estimate the second term in the right hand side \eqref{estim:B:0} thanks to Jensen inequality
\begin{equation*}
	\lambda
	\left\|
	\partial_x^l
	\int_{0}^{t}
	\exp{\left(\frac{s}{\lambda}J\right)}\cdot\bL^\lambda_3\dD s\right\|_{L^2}
	\leq\,\lambda \,t\sup_{0\leq s\leq T}\left\|\partial_x^l\bL^\lambda_3(s)\right\|_{L^2}
	\,,
\end{equation*}
where $\bL_3^\lambda$ is defined below \eqref{reform:L}.
Then, we substitute the time derivatives of $C_1^\lambda$ and $C_3^\lambda$ in the expression of $\bL_3^\lambda$ thanks to the first line of \eqref{Hermite:D} with $k=1$ and $k=3$ respectively, which yields
\begin{equation*}
\begin{split}
	&\lambda
	\left\|
	\partial_x^l
	\int_{0}^{t}
	\exp{\left(\frac{s}{\lambda}J\right)}\cdot\bL^\lambda_3\dD s\right\|_{L^2}
	\leq\\
	&\lambda T_0
	t\sup_{0\leq s\leq T}\left\|\partial^{l+2}_x
	\left(
	\sqrt{18}  \, E^\lambda C_2^\lambda
	-\sqrt{18}\,T_0\,\partial_x C_2^\lambda
	-\sqrt{24}\,T_0\,\partial_x C_4^\lambda
	+
	2E^\lambda C_{0}^\lambda
	-
	2T_0\,
	\partial_x C^\lambda_{0}\,- \, 2\sqrt{2}\,T_0\,
	\partial_x C^\lambda_{2}+E^\lambda
	\right)\right\|_{L^2}.
	\end{split}
\end{equation*}
We apply the triangle inequality and obtain
\begin{equation*}
	\begin{split}
		\lambda
		\left\|
		\partial_x^l
		\int_{0}^{t}
		\exp{\left(\frac{s}{\lambda}J\right)}\cdot\bL^\lambda_3\dD s\right\|_{L^2}
		\,\leq\,
		&C
\lambda \sup_{0\leq s\leq T}
\left(\left\|\partial^{l+2}_x E^\lambda(s)\right\|_{L^2}
+
\left\|\partial^{l+2}_x
(E^\lambda C_2^\lambda
)(s)\right\|_{L^2}+
\left\|\partial^{l+2}_x
(E^\lambda C_0^\lambda
)(s)\right\|_{L^2}
		\right)\\
		+
		&C
		\lambda \sup_{0\leq s\leq T}\left(
		\left\|\partial^{l+3}_x  C_2^\lambda(s)\right\|_{L^2}+
		\left\|\partial^{l+3}_x C_4^\lambda(s)\right\|_{L^2}+
		\left\|\partial^{l+3}_x C_0^\lambda(s)\right\|_{L^2}
		\right)
		\,,
	\end{split}
\end{equation*}
for some constant $C$ depending on the final time $T$ and the temperature $T_0$.
Then, we apply Leibniz rule to develop the nonlinear terms and estimate each product thanks to Young inequality
\begin{equation}\label{estim:L3}
	\begin{split}
		\lambda
		\left\|
		\partial_x^l
		\int_{0}^{t}
		\exp{\left(\frac{s}{\lambda}J\right)}\cdot\bL^\lambda_3\dD s\right\|_{L^2}
		\,\leq\,
		&C
		\sup_{0\leq s\leq T}
		\left(\lambda \left\|\partial^{l+2}_x E^\lambda(s)\right\|_{L^2}
		,\;\lambda^{1+\alpha}\left\| E^\lambda(s)\right\|_{W^{l+2,4}}^2
		\right)\\
		+
		&\,C \sup_{\substack{0\leq s\leq T\\k=0,2,4}}\left(
		\lambda\left\|\partial^{l+3}_x C_k^\lambda(s)\right\|_{L^2}
		,\;
		\lambda^{1-\alpha}
		\left\|  C_k^\lambda(s)\right\|_{W^{l+2,4}}^2
		\right)
		\,,
	\end{split}
\end{equation}
for some constant $C$ depending on $l\geq 0$, on the final time $T$ and the temperature $T_0$, where, here and in the sequel, for two nonnegative functions $a$ and $b$, we use the notation $\sup_{s}\left(a(s),\,b(s)\right)$ for $\sup_{s}\max\left(a(s),\,b(s)\right)$. We use the same method to estimate the last four terms on the right hand side of \eqref{estim:B:0} and obtain
\begin{equation}\label{estim:last4}
	\begin{split}
\lambda^2
\left\|
\partial_x^l\bN^\lambda_1(t)\right\|_{L^2}
+\lambda
\left\|
\partial_x^l\bL^\lambda_2(t)\right\|_{L^2}
&+ \lambda^2
\left\|
\partial_x^l
\bN^\lambda_1(0)\right\|_{L^2}
+\,
\lambda\,
\left\|
\partial_x^l\bL^\lambda_2(0)\right\|_{L^2}\\[0.8em]
&\leq\,C\,\lambda^2
\sup_{0\leq s\leq T}\,
\left(
\| E^\lambda(s)\|_{W^{l+1,4}}^2
\right)+\,C\,\lambda
 \sup_{\substack{0\leq s\leq T\\k=1,3}}\,
\left(
\|\partial^{l+2}_x  C_k^\lambda(s)\|_{L^2}
\right),
\end{split}
\end{equation}
for some constant $C$ depending on $l\geq 0$ and on the temperature $T_0$. Gathering estimates \eqref{estim:N4}, \eqref{estim:L3} and \eqref{estim:last4} in \eqref{estim:B:0}, we find
\begin{equation*}
	\begin{split}
		\left\|
		\partial_x^l
		B(t)\right\|_{L^2}
		\leq\;
		&C\left(\lambda^{\alpha-1}
		\left\|
		C^\lambda_{1}(0)
		\right\|_{W^{l+1,4}}^2
		+
		\sup_{0\leq s\leq T}
		\left(\lambda \left\|\partial^{l+2}_x E^\lambda(s)\right\|_{L^2}
		,\;\lambda^{1+\alpha}\left\| E^\lambda(s)\right\|_{W^{l+2,4}}^2
		\right)\right)
		\\
		+
		&\ds\,C \sup_{\substack{0\leq s\leq T\\0\leq k\leq 4}}\left(
		\lambda\left\| C_k^\lambda(s)\right\|_{H^{l+3}}
		,\;
		\lambda^{1-\alpha}
		\left\|  C_k^\lambda(s)\right\|_{W^{l+2,4}}^2
		\right)
		\,,
	\end{split}
\end{equation*}
for all $0<\lambda<1$ and for some constant $C$ depending on $l\geq 0$, on the final time $T$ and the temperature $T_0$. We sum the latter estimate over all $l\leq n$ and take the supremum in $t$ on the left hand side, which yields
\begin{equation*}
	\begin{split}
		\sup_{0\leq t\leq T}
		\left\|
		B(t)\right\|_{H^{n}}
		\leq\;
		&C\left(\lambda^{\alpha-1}
		\left\|
		C^\lambda_{1}(0)
		\right\|_{W^{n+1,4}}^2
		+
		\sup_{0\leq t\leq T}
		\left(\lambda \left\|E^\lambda(t)\right\|_{H^{n+2}}
		,\;\lambda^{1+\alpha}\left\| E^\lambda(t)\right\|_{W^{n+2,4}}^2
		\right)\right)
		\\
		+
		&\ds\,C \sup_{\substack{0\leq t\leq T\\0\leq k\leq 4}}\left(
		\lambda\left\|C_k^\lambda(t)\right\|_{H^{n+3}}
		,\;
		\lambda^{1-\alpha}
		\left\|  C_k^\lambda(t)\right\|_{W^{n+2,4}}^2
		\right)
		\,,
	\end{split}
\end{equation*}
for all $0<\lambda<1$ and $n\geq 0$, where the constant $C$ depends on $n\geq 0$, on the final time $T$ and the temperature $T_0$. We bound $L^2\left(\T\right)$ norms with their $L^4\left(\T\right)$ counterparts, which yields
\begin{equation*}
	\begin{split}
		\sup_{0\leq t\leq T}
		\left\|
		B(t)\right\|_{H^{n}}
		\leq\;
		&C\left(\lambda^{\alpha-1}
		\left\|
		C^\lambda_{1}(0)
		\right\|_{W^{n+1,4}}^2
		+
		\sup_{0\leq t\leq T}
		\left(\lambda \left\|E^\lambda(t)\right\|_{W^{n+2,4}}
		,\;\lambda^{1+\alpha}\left\| E^\lambda(t)\right\|_{W^{n+2,4}}^2
		\right)\right)
		\\
		+
		&\ds\,C \sup_{\substack{0\leq t\leq T\\0\leq k\leq 4}}\left(
		\lambda\left\|C_k^\lambda(t)\right\|_{W^{n+3,4}}
		,\;
		\lambda^{1-\alpha}
		\left\|  C_k^\lambda(t)\right\|_{W^{n+3,4}}^2
		\right)
		\,.
	\end{split}
\end{equation*}
To conclude, we fix some $n\in\N^\star$ such that $n\leq r_0$ and estimate the norms of $E^\lambda$ and $C_1^\lambda(0)$ on the first line thanks to assumption \eqref{compatibility} and the norms of coefficient $C_k^\lambda$, $0\leq k\leq 4$, on the second line thanks to assumption \eqref{uniform:reg}. This gives
\begin{equation}\label{estim:B:final}
	\begin{split}
		\sup_{0\leq t\leq T}
		\left\|
		B(t)\right\|_{H^{n}}
		\leq\;
		C
		\lambda^{1-\alpha}\,,
	\end{split}
\end{equation}
for all $0<\lambda<1$ and $0\leq n\leq r_0$, where $C$ depends on $\alpha$, on the final time $T$, the temperature $T_0$, the size of the domain $\T$ and the implicit constants in \eqref{compatibility}-\eqref{uniform:reg}, and where $r_0$ is given in \eqref{compatibility}.\\
Using the same approach, we find the following estimates for $A_1$ and $A_2$ (defined by \eqref{reform:N})
\begin{equation*}
	\sup_{0\leq t\leq T}\left(
	\left\|A_1(t)\right\|_{W^{n,\infty}},\,
	\left\|A_2(t)\right\|_{W^{n,\infty}}\right)\leq\;
	C
	\sup_{0\leq t\leq T}\left(
	\lambda\left\|E^\lambda(t)\right\|_{W^{{n+1},\infty}},\,\left\|C_1^\lambda(0)\right\|_{W^{{n},\infty}},\,\lambda\left\|C^\lambda_2(0)\right\|_{W^{{n+1},\infty}}\right).
\end{equation*}
We use Sobolev injections to bound the $L^\infty(\T)$ norms with their $W^{1,4}(\T)$ counterpart, which yields
\begin{equation*}
	\sup_{0\leq t\leq T}\left(
	\left\|A_1(t)\right\|_{W^{n,\infty}},\,
	\left\|A_2(t)\right\|_{W^{n,\infty}}\right)\leq\;
	C
	\sup_{0\leq t\leq T}\left(
	\lambda\left\|E^\lambda(t)\right\|_{W^{{n+2},4}},\,\left\|C_1^\lambda(0)\right\|_{W^{n+1,4}},\,\lambda\left\|C^\lambda_2(0)\right\|_{W^{{n+2},4}}\right).
\end{equation*}
Using assumption \eqref{compatibility} to estimate $E^\lambda$ and $C_1^\lambda$ and assumption \eqref{uniform:reg} for $C_2^\lambda$, we find 
\begin{equation}\label{estim:A:final}
\sup_{0\leq t\leq T}\left(
\left\|A_1(t)\right\|_{W^{n,\infty}},\,
\left\|A_2(t)\right\|_{W^{n,\infty}}\right)\leq\;
C
\lambda^{1-\alpha}\,,
\end{equation}
for all $0<\lambda<1$ and $0\leq n\leq r_0$, where $C$ depends on $\alpha$, on the final time $T$, the temperature $T_0$, the size of the domain $\T$ and the implicit constants in \eqref{compatibility}-\eqref{uniform:reg}, and where $r_0$ is given in \eqref{compatibility}.
\\

\paragraph{\bf \underline{Step 2.} iteration of the Duhamel formula.}  We now focus on the iteration process. We consider some integer $l\geq 0$, take the $l$-th derivative in $x$ and then the $L^2(\T)$-norm in \eqref{duham:reform}. After applying Jensen inequality to bound the time integral, using that $\exp((s-t)J/\lambda)$ is an isometry of $\R^2$ and using Leibniz rule to estimate the products between $A_i$, $i\in \{1,2\}$, and $\bU_1^\lambda$, we obtain
\begin{equation*}
	\sup_{0\leq t\leq T}
	\left\|\partial^l_x\bU^\lambda(t)\right\|_{L^2}
	\leq\,
	C
	\sup_{0\leq t\leq T}\left(
	\left\|A_1(t)\right\|_{W^{l,\infty}},\,
	\left\|A_2(t)\right\|_{W^{l,\infty}}\right)
	\sup_{0\leq t\leq T}\left\|\bU_1^\lambda(t)\right\|_{H^{l+1}}
	+\sup_{0\leq t\leq T}\left\|\partial^l_x B(t)\right\|_{L^2},
\end{equation*}
for some constant $C$ depending on $l\geq0$ and $T$. We sum the latter estimate over all integers $l\geq0$ less than $n$, for any $n\in\N$ such that $n\leq r_0$, where $r_0$ is given in \eqref{compatibility}
\begin{equation*}
	\sup_{0\leq t\leq T}\left\|\bU^\lambda(t)\right\|_{H^{n}}
	\leq\,
	C\left(
	\sup_{0\leq t\leq T}\left(
	\left\|A_1(t)\right\|_{W^{n,\infty}},\,
	\left\|A_2(t)\right\|_{W^{n,\infty}}\right)
	\sup_{0\leq t\leq T}\left\|\bU_1^\lambda(t)\right\|_{H^{n+1}}
	+\sup_{0\leq t\leq T}\left\| B(t)\right\|_{H^{n}}\right),
\end{equation*} 
for some constant $C$ depending on $n\geq0$ and $T$. We bound the norms of $(A_1,A_2)$ according to \eqref{estim:A:final} and the norm of $B$ according to \eqref{estim:B:final}, which yields
\begin{equation*}
	\sup_{0\leq t\leq T}\left\|\bU^\lambda(t)\right\|_{H^{n}}
	\leq\,
	C\lambda^{1-\alpha}\left(
	\sup_{0\leq t\leq T}\left\|\bU_1^\lambda(t)\right\|_{H^{n+1}}
	+1\right),
\end{equation*} 
for all $0<\lambda<1$ and $0\leq n\leq r_0$, where the constant $C$ depends on $\alpha$, on the final time $T$, the temperature $T_0$, the size of the domain $\T$ and the implicit constants in \eqref{compatibility}-\eqref{uniform:reg}. Next, we point out that $\left\|\bU_1^\lambda\right\|_{H^{n+1}}\leq \left\|\bU^\lambda\right\|_{H^{n+1}}$, which allows to iterate the latter formula with respect to $n$ and obtain
\begin{equation*}
	\sup_{0\leq t\leq T}\left\|\bU^\lambda(t)\right\|_{L^{2}}
	\leq\,
	C\left(
	\lambda^{n\left(1-\alpha\right)}
	\sup_{0\leq t\leq T}\left\|\bU_1^\lambda(t)\right\|_{H^{n+1}}
	+\lambda^{1-\alpha}\right),
\end{equation*} 
for all $0<\lambda<1$ and $0\leq n\leq r_0$. Now, we fix $n\,=\,r_0$, and since $r_0\geq 1/(1-\alpha) $, we obtain
\begin{equation*}
	\sup_{0\leq t\leq T}\left\|\bU^\lambda(t)\right\|_{L^{2}}
	\leq\,
	C\left(
	\lambda
	\sup_{0\leq t\leq T}\left\|\bU_1^\lambda(t)\right\|_{H^{r_0+1}}
	+\lambda^{1-\alpha}\right),
\end{equation*} 
for all $0<\lambda<1$, where $C$ depends on $\alpha$, $T$, $T_0$, $\T$ and the implicit constants in \eqref{compatibility}-\eqref{uniform:reg}.
To bound the norm of $\bU_1^\lambda$, we substitute $\bU_1^\lambda$
according to the following relation $\bU^\lambda_1\,=\,E^\lambda - E^\lambda_{\textrm{slow}} 
-\,E^\lambda_{\textrm{osc}}$ and apply the triangle inequality
\begin{equation*}
	\sup_{0\leq t\leq T}\left\|\bU^\lambda(t)\right\|_{L^{2}}
	\leq\,
	C\left(
	\lambda
	\sup_{0\leq t\leq T}\left(\left\|E^\lambda(t)\right\|_{H^{r_0+1}}
	+\left\|C_2^\lambda(t)\right\|_{H^{r_0+2}}+\left\|E^\lambda_{\textrm{osc}}(t)\right\|_{H^{r_0+1}}
	\right)
	+
	\lambda^{1-\alpha}
	\right).
\end{equation*} 
To bound the norm of $E^\lambda_{\textrm{osc}}$, we substitute $E^\lambda_{\textrm{osc}}$ according to \eqref{E:C_1:osc} and apply the triangle inequality
\begin{equation*}
	\sup_{0\leq t\leq T}\left\|\bU^\lambda(t)\right\|_{L^{2}}
	\leq\,
	C\left(
	\sup_{0\leq t\leq T}\left(\lambda\left\|E^\lambda(t)\right\|_{H^{r_0+1}}
	,\,\left\|C^{\lambda}_1(0)\right\|_{H^{r_0+1}},\,\lambda\left\|C_2^\lambda(t)\right\|_{H^{r_0+2}}
	\right)
	+
	\lambda^{1-\alpha}
	\right).
\end{equation*} 
We bound $E^\lambda$ and $C_1^\lambda$ thanks to \eqref{compatibility} and $C_2^\lambda$ thanks to \eqref{uniform:reg}, which yields the result
\begin{equation*}
	\sup_{0\leq t\leq T}\left\|\bU^\lambda(t)\right\|_{L^{2}}
	\leq\,
	C
	\lambda^{1-\alpha}\,,
\end{equation*} 
for all $0<\lambda<1$, where the constant $C$ depends on $\alpha$, on the final time $T$, the temperature $T_0$, the size of the domain $\T$ and the implicit constants in \eqref{compatibility}-\eqref{uniform:reg}.

\section{Proof of Proposition \ref{prop:discrete:QN:limit}}\label{proof:prop:discrete:QN:limit}
In this proof, we analyze the limit $\lambda \rightarrow 0$ for a fixed $\Delta t>0$. Without loss of generality, we assume that $0<\lambda\leq 1$. Indeed, when $\lambda>1$, estimates \eqref{result:h} directly follow from assumptions \eqref{hyp:E:disc}-\eqref{hyp:Ck:disc}, since in this case $\cE_1^n\leq C\leq C\lambda$ and $\cE_0^n\leq \|E^n\|_{l^2(\cJ)}+\sqrt{2}\,T_0\|C^n_2\|_{h^1(\cJ)}\leq C/\lambda+C\leq 2C\lambda$. We denote by $C_k^{(1,n+1)}$ the coefficient $C^{(1)}_{k}$ computed in
\eqref{discrete:step1} to get $C^{n+1}$ from $C^n$. We first recall the discrete functional inequalities used in the proof.
\begin{lemma}\label{lem:disc:ineq}
Suppose that $N_x$ is odd. Then there exists a constant $K>0$, only depending on $|b-a|$, such that for all $C_h=(\cC_j)_{j\in\J}\in\R^\cJ$ with $\sum_{j\in\J}h\,\cC_j=0$,
\begin{equation}\label{lem:poincare:sobolev}
\left\|C_h\right\|_{l^2(\cJ)}\,\leq\,\frac{|b-a|}{2}\,\left\|\partial_h C_h\right\|_{l^2(\cJ)}\,,
\end{equation}
and for all $C_h\in\R^\cJ$,
\begin{equation}\label{lem:sobolev}
\left\|C_h\right\|_{l^\infty(\cJ)}\,\leq\,K\,\left\|C_h\right\|_{h^1(\cJ)}\,.
\end{equation}
\end{lemma}
These inequalities are classical, see for instance \cite[Lemma $3.3$]{BF_09_22}.\\
Our proof is divided into $2$ steps. In the first one, we solve the linearized scheme \eqref{discrete:step1} to prove that $C_0^{n}$ converges to $1$, that $C_1^{(1,n)}$ converges to $0$ and that $C_2^{(1,n)}$ remains bounded uniformly in $\lambda$. In the second step, we analyze the nonlinear scheme \eqref{discrete:step2} and use our estimate from the first step to prove that $C_1^{n}$ converges to $0$ and that $E^{n}$ is close to $\sqrt{2}\,T_0\, \partial_h\,C^{n}_{2}$ as $\lambda\rightarrow 0$. \\

\paragraph{\bf \underline{Step 1.} Analysis of the linear step \eqref{discrete:step1}.}
The convergence of $C_0^n$ is a consequence of the Poisson coupling between $E^{n}$ and $C_0^{n}$. Indeed, taking the $h^2(\cJ)$ norm of the third line of \eqref{discrete:step1} and replacing $C_0^{(1)}$ and $E^{(1)}$ with $C_0^{n+1}$ and $E^{n+1}$ thanks to the first two lines in \eqref{discrete:step2}, we deduce
\[
\left\| C_0^{n+1} - 1\right\|_{h^2\left(\cJ\right)}\,=\, 
\lambda ^2\,\left\|
\partial_h E^{n+1}
\right\|_{h^2\left( \cJ \right) }\,\leq\,\lambda ^2\,\left\|
E^{n+1}
\right\|_{h^3\left( \cJ \right) }\,,
\]
for all $0\leq n+1\leq T/\Delta t$. The latter relation is also valid when $n+1=0$. Hence, we bound the right hand side according to \eqref{hyp:E:disc}, which yields
\begin{equation}\label{estim:C0:h}
	\sup_{0\leq n\leq T/\Delta t}\;\left\|C_0^{n} - 1\right\|_{h^2\left(\cJ\right)} \,\leq\, C\lambda \,.
\end{equation}
To demonstrate that $\left(C_1^{(1,n)}\right)_{1\leq n \leq T/\Delta t}$ converges, we take the $h^2(\cJ)$ norm in the second line of \eqref{discrete:step1} with $k=0$, and use that $C_0^{(1,n+1)}=C_0^{n+1}$ according to the first line of \eqref{discrete:step2}, which yields
\[
\sqrt{T_0}
\left\|
\,\partial_h\,C^{(1,n+1)}_{1}\right\|_{h^2\left(\cJ\right)}
\,=\,
\left\|\frac{C^{n+1}_{0} -C_{0}^{n} }{\Delta t}  \right\|_{h^2\left(
	\cJ \right)}
\,\leq\,\frac{1}{\Delta t}\left(\left\|C^{n+1}_{0} -1\right\|_{h^2\left(\cJ \right)}+\left\|C^{n}_{0} -1\right\|_{h^2\left(\cJ \right)}\right)\,,
\]
for all $0\leq n\leq T/\Delta t-1$. We bound the right hand side according to \eqref{estim:C0:h}, take the supremum over $n$ and shift the index $n+1\rightarrow n$, which gives 
\begin{equation}
  \label{bound:C1n:0}
\sup_{1\leq n\leq T/\Delta t}\left\|
\partial_h C^{(1,n)}_{1}\right\|_{h^2\left(\cJ\right)}
\,\leq\,C\lambda\,,
\end{equation}
where the constant $C>0$ depends on $\Delta t$. To control the full $h^3(\cJ)$ norm of $C_1$, we bound the $l^2(\cJ)$ norm of $C_1$ with the $l^2(\cJ)$ norm of $\partial_h C_1$ thanks to the discrete Poincaré-Wirtinger inequality \eqref{lem:poincare:sobolev} of Lemma \ref{lem:disc:ineq}, which reads
\begin{align*}
	\left\|C_1^{n}
	\right\|_{l^2\left(\cJ\right)}
	\,\leq\,\ds\,
	C
	\left\|
	\partial_h C_{1}^{n}
	\right\|_{l^2\left(\cJ\right)}
	\,+\,C\left|
	\sum_{j\in\J}h\,\cC^{n}_{1,j}
	\right|
	\,.
\end{align*}
We check that the total flux on the right hand side is $0$ thanks to assumption \eqref{total:flux:mass:h}. Indeed, summing the first line of \eqref{discrete:step1} over $j\in\J$, the terms involving $\partial_h$ vanish by periodicity, as well as the sum of $E^{(1)}=-\partial_h\phi^{(1)}$. Then, summing \eqref{discrete:step2} with $k=1$ over $j\in\J$ and using the Poisson equation, the nonlinear term reads
  $$
  \lambda^2\sum_{j\in\J}h\,E_j^{n+1}\,\partial_jE^{n+1},
  $$ which vanishes by skew-symmetry of $\partial_h$. Hence, by induction on $n\geq 0$,
\[
\sum_{j\in\J}h \,\cC^{n}_{1,j}\,
\,=\,
\sum_{j\in\J}h\,\cC^{(1,n+1)}_{1,j}\, \,=\,0\,.
\]
Therefore, we obtain
\begin{equation}\label{disc:poincaré}
	\left\|C_1^{n}
	\right\|_{l^2\left(\cJ\right)}
	\,\leq\,\ds\,
	C
	\left\|
	\partial_h C_{1}^{n}
	\right\|_{l^2\left(\cJ\right)}
	\,,\quad \textrm{and}\quad 
	\left\|C_1^{(1,n+1)}
	\right\|_{l^2\left(\cJ\right)}
	\,\leq\,\ds\,
	C
	\left\|
	\partial_h C_{1}^{(1,n+1)}
	\right\|_{l^2\left(\cJ\right)}.
\end{equation}
Together with \eqref{bound:C1n:0}, the last inequality yields
\begin{equation}
	\label{bound:C1n}
	\sup_{1\leq n\leq T/\Delta t}\left\|
C^{(1,n)}_{1}\right\|_{h^3\left(\cJ\right)}
	\,\leq\,C\lambda\,.
\end{equation}
Next, we show that $C_2^{(1,n)}$ is bounded for all $1\leq n \leq T/\Delta t$. To do so, we fix $0\leq r\leq 2$ and set $D_k:=\partial_h^r C^{(1,n+1)}_{k}$ and $D^n_k:=\partial_h^r C^{n}_{k}$. Since $\partial_h$ commutes with translations, we may apply $\partial_h^r$ to the second line of \eqref{discrete:step1}, which gives for $2\leq k\leq N_H$
\[
\frac{D_k-D_k^n}{\Delta t}\,+\,\sqrt{k\,T_0}\,\partial_h D_{k-1}\,+\,\sqrt{(k+1)\,T_0}\,\partial_h D_{k+1}\,=\,0\,,
\]
with $D_{N_H+1}=0$. We take the discrete $l^2(\cJ)$ scalar product $\langle C,D\rangle=\sum_{j\in\J}h\,\cC_j\,\cD_j$ of this equation with $D_k$ and sum over $k \in\{2,\dots , N_H\}$. Thanks to the skew-symmetry of $\partial_h$, we have $\langle \partial_h D_{k},D_{k-1}\rangle=-\langle D_{k},\partial_h D_{k-1}\rangle$, hence the contributions of the terms $\partial_h D_{k-1}$ and $\partial_h D_{k+1}$ cancel two by two, except the one involving $D_1$ in the equation for $k=2$. Therefore, we get
\[
\sum_{k=2}^{N_H}\left\langle D_k-D_k^n,D_k\right\rangle\,=\,-\Delta t\,\sqrt{2\,T_0}\,\left\langle \partial_h D_{1},D_2\right\rangle\,.
\]
Then, we use that $\langle D_k-D_k^n,D_k\rangle\geq \frac{1}{2}\left(\|D_k\|^2_{l^2(\cJ)}-\|D^n_k\|^2_{l^2(\cJ)}\right)$ and the Cauchy-Schwarz inequality, which yields
\[
\sum_{k=2}^{N_H}
\left\| \partial_h^r C^{(1,n+1)}_{k}\right\|_{l^2\left(\cJ\right)}^2
\leq\,
2\,\Delta t\,\sqrt{2\,T_0}
\left\|\partial_h^{r+1} C^{(1,n+1)}_{1}\right\|_{l^2\left(\cJ\right)}
\left\| \partial_h^r C^{(1,n+1)}_{2}\right\|_{l^2\left(\cJ\right)}\,+\,
\sum_{k=2}^{N_H}
\left\| \partial_h^r C^{n}_{k}\right\|_{l^2\left(\cJ\right)}^2\,.
\]
We bound $\partial_h^{r+1}C^{(1,n+1)}_{1}$ on the right hand side according to \eqref{bound:C1n}, since $r+1\leq 3$ and $\lambda\leq 1$, and we estimate the sum thanks to assumption \eqref{hyp:Ck:disc}. Keeping only the term $k=2$ on the left hand side, this yields 
\[
\left\| \partial_h^r C^{(1,n+1)}_{2}\right\|_{l^2\left(\cJ\right)}^2
\leq\,
C\left(
\left\|\partial_h^r C^{(1,n+1)}_{2}\right\|_{l^2\left(\cJ\right)}\,+\,
1\right)\,,
\]
where the constant $C$ depends on $\Delta t$ and $T_0$, but neither on $h$ nor on $N_H$. Applying Young inequality to the first term of the right hand side and summing over $0\leq r\leq 2$, we deduce that
\begin{equation}\label{bound:C2n}
	\sup_{1\leq n\leq T/\Delta t}\left\|
	\,C^{(1,n)}_{2}\right\|_{h^2\left(\cJ\right)}
	\,\leq\,C
	\,.
\end{equation}
\paragraph{\bf \underline{Step 2.} Analysis of the nonlinear step \eqref{discrete:step2}.}
The last step consists in proving that $C_1^{n}$ is of order $\lambda$  for all $1\leq n\leq T/\Delta t$ and that $E^{n}$ is close to $\sqrt{2}\,T_0\, \partial_h\,C^{n}_{2}$ for all $2\leq n\leq T/\Delta t$ as $\lambda\rightarrow 0$. 
We proceed in three stages:
\begin{itemize}
	\item we first prove that $E^{n}$ is uniformly bounded in $\lambda$
	for all $1\leq n\leq T/\Delta t$ ;
	\item we deduce that $C_1^{n}$ is of order $\lambda$  for all $1\leq
	n\leq T/\Delta t$ ;
	\item we obtain that $E^{n}$ is close to $\sqrt{2}\,T_0\,
	\partial_h\,C^{n}_{2}$ for all $2\leq n\leq T/\Delta t$ as
	$\lambda\rightarrow 0$.
\end{itemize}
Let us first prove that $E^{n}$ is uniformly bounded in $\lambda$ for all $n\geq 1$. We take the $h^1(\cJ)$ norm in the first line of \eqref{discrete:step1}, which yields
\[
\left\|E^{(1,n+1)}
\right\|_{h^1\left(\cJ\right)}
\,\leq\,
\sqrt{2}\,T_0 \left\|
 \partial_h\,C^{(1,n+1)}_{2}
\right\|_{h^1\left(\cJ\right)}
\,+\,
\frac{\sqrt{T_0}}{\Delta t}
\left(
\left\|
C^{(1,n+1)}_{1} 
\right\|_{h^1\left(\cJ\right)}
+
\left\|
C_{1}^{n} 
\right\|_{h^1\left(\cJ\right)}
\right)
\,+\,
T_0
\left\|\partial_h C^{(1,n+1)}_{0}
\right\|_{h^1\left(\cJ\right)}.
\]
On the right hand side, we estimate $C_0$ according to \eqref{estim:C0:h}, $C_1^{(1,n+1)}$ according to \eqref{bound:C1n}, $C^{(1,n+1)}_{2}$ according to \eqref{bound:C2n} and assumption \eqref{hyp:E:disc} to estimate $C_1^n$. We obtain that for all $\lambda>0$, it holds
\begin{equation}\label{bound:En:1}
	\sup_{1\leq n\leq T/\Delta t}\left\|
	E^{n}\right\|_{h^1\left(\cJ\right)}
	\,\leq\,C
	\,.
\end{equation}
Next, we deduce that $C_1^n$ is of order $\lambda$  for all $1\leq n\leq T/\Delta t$. To do so, we  take the $l^2(\cJ)$ norm of the derivative of \eqref{discrete:step2} with $k=1$, which yields 
\[
\left\|\partial_h C^{n+1}_{1}\right\|_{l^2\left(\cJ\right)}
\,\leq\,
\left\| \partial_h C^{(1,n+1)}_{1}\right\|_{l^2\left(\cJ\right)}
\,+\,
\frac{\Delta t}{\sqrt{T_0}}\,
\left\| \partial_h \left(E^{n+1}\left(C^{n+1}_{0}\,-\,1\right) \right) \right\|_{l^2\left(\cJ\right)} \,.
\]
We estimate $C^{(1,n+1)}_{1}$ thanks to \eqref{bound:C1n}, which gives
\begin{equation*}
\left\|\partial_h C^{n+1}_{1}\right\|_{l^2\left(\cJ\right)}
\,\leq\,
C
\left(
\left\| \partial_h\left(E^{n+1}\left(C^{n+1}_{0}\,-\,1\right)\right)  \right\|_{l^2\left(\cJ\right)}
+\lambda 
\right)
\,.
\end{equation*}
We estimate the derivative of the product between $E^{n+1}$ and $C_0^{n+1}-1$, as follows 
\begin{align*}
	\left\|\partial_h C^{n+1}_{1}\right\|_{l^2\left(\cJ\right)}
	\,\leq\,\ds \,C
	\left(
	\left\| E^{n+1}\right\|_{l^\infty\left(\cJ\right)}
	\left\| C_0^{n+1}-1\right\|_{h^1\left(\cJ\right)}
	+
	\left\| E^{n+1}\right\|_{h^1\left(\cJ\right)}
	\left\| C_0^{n+1}-1\right\|_{l^\infty\left(\cJ\right)}
	+\lambda 
	\right).
\end{align*}
To estimate the $l^\infty(\cJ)$ norms, we use the following discrete Sobolev inequality, given by \eqref{lem:sobolev} in Lemma \ref{lem:disc:ineq}, with a constant $C$ only depending on $|b-a|$
\begin{equation}\label{sob:ineq:E}
	\left\| E^{n+1}\right\|_{l^\infty\left(\cJ\right)} \,\leq\,C
	\left\| E^{n+1}\right\|_{h^1\left(\cJ\right)}\,.
\end{equation}
We deduce
\begin{align*}
	\left\|\partial_h C^{n+1}_{1}\right\|_{l^2\left(\cJ\right)}
	\,\leq\,\ds \,C
	\left(\left\| E^{n+1}\right\|_{h^1\left(\cJ\right)}
	\left\| C_0^{n+1}-1\right\|_{h^1\left(\cJ\right)}
	+\lambda 
	\right).
\end{align*}
We estimate the norm of $E^{n+1}$ thanks to \eqref{bound:En:1}, the norm of $C_0^{n+1}-1$ thanks to \eqref{estim:C0:h} and take the supremum over all $1\leq n+1 \leq T/\Delta t$, we find 
\begin{equation*}
	\sup_{1\leq n\leq T/\Delta t}\left\|\partial_h 
	C_1^{n}\right\|_{l^2\left(\cJ\right)}
	\,\leq\,C\lambda
	\,,\quad \forall \lambda >0\,.
\end{equation*}
To control the full $h^1(\cJ)$ norm of $C_1^n$, we apply the discrete Poincaré inequality \eqref{disc:poincaré}, which yields the second estimate in \eqref{result:h}
\begin{equation}\label{estim:dh:C1}
	\sup_{1\leq n\leq T/\Delta t}\left\|
	C_1^{n}\right\|_{h^1\left(\cJ\right)}
	\,\leq\,C\lambda
	\,,\quad \forall \lambda >0\,.
\end{equation}
To conclude, we prove that $E^n$ is close to $\sqrt{2}\,T_0\, \partial_h\,C^{n}_{2}$ for all $2\leq n\leq T/\Delta t$ as $\lambda\rightarrow 0$. To do so, we take the $l^2(\cJ)$ norm in the first line of \eqref{discrete:step1}, which yields
\[
\left\|E^{(1,n+1)}
-
\sqrt{2}\,T_0\, \partial_h C^{(1,n+1)}_{2}
\right\|_{l^2\left(\cJ\right)}
\,\leq\,
\frac{\sqrt{T_0}}{\Delta t}
\left(
\left\|
C^{(1,n+1)}_{1} 
\right\|_{l^2\left(\cJ\right)}
+
\left\|
C_{1}^{n} 
\right\|_{l^2\left(\cJ\right)}
\right)
\,+\,
T_0
\left\|\partial_h C^{(1,n+1)}_{0}
\right\|_{l^2\left(\cJ\right)}.
\]
On the right hand side, we estimate $C_0$ according to \eqref{estim:C0:h}, $C_1^{(1,n+1)}$ according to \eqref{bound:C1n}, and $C_1^n$ according to the latter estimate in \eqref{estim:dh:C1}, which yields 
\begin{equation}\label{estim:E1:C2}
	\left\|E^{n+1}
	-
	\sqrt{2}\,T_0\, \partial_h C^{(1,n+1)}_{2}
	\right\|_{l^2\left(\cJ\right)}
	\,\leq\,
	C\lambda
	\,,
\end{equation}
for all $n\geq 1$. Then, we prove that $\partial_h C^{(1,n+1)}_{2}$ is close to $\partial_h C^{n+1}_{2}$ taking the $l^2(\cJ)$ norm in \eqref{discrete:step2} with $k=2$, which yields
\[
\left\|\partial_h \left(C^{n+1}_{2} -C_{2}^{(1,n+1)}\right)\right\|_{l^2\left(\cJ\right)}\, \leq  \,
C\left\|\partial_h
\left(E^{n+1} C^{n+1}_{1}\right)\right\|_{l^2\left(\cJ\right)}
\]
We estimate the derivative of the product $\partial_h
\left(E^{n+1} C^{n+1}_{1}\right)$ as follows
\[
\left\|\partial_h \left(C^{n+1}_{2} -C_{2}^{(1,n+1)}\right)\right\|_{l^2\left(\cJ\right)}\, \leq  \,
C\left(\left\|E^{n+1} \right\|_{l^\infty\left(\cJ\right)}
\left\|\partial_h C^{n+1}_{1}\right\|_{l^2\left(\cJ\right)}
+
\left\|\partial_h E^{n+1} \right\|_{l^2\left(\cJ\right)}
\left\| C^{n+1}_{1}\right\|_{l^\infty\left(\cJ\right)}
\right)\,.
\]
Then, we apply the Sobolev inequality \eqref{sob:ineq:E} to estimate the $L^\infty$ norms in the last inequality and deduce 
\[
\left\|\partial_h\left(C^{n+1}_{2} -C_{2}^{(1,n+1)}\right)\right\|_{l^2\left(\cJ\right)}\, \leq  \,
C\left\|
E^{n+1}\right\|_{h^1\left(\cJ\right)} \left\|C^{n+1}_{1}\right\|_{h^1\left(\cJ\right)}\,.
\]
We estimate the norm of $E^{n+1}$ thanks to \eqref{bound:En:1} and the norm of $C_1^{n+1}$ thanks to \eqref{estim:dh:C1}, which yields 
\[
\left\|\partial_h\left(C^{n+1}_{2} -C_{2}^{(1,n+1)}\right)\right\|_{l^2\left(\cJ\right)}\, \leq  \,
C\lambda\,,
\]
for all $n\geq 1$. Plugging the latter estimate in \eqref{estim:E1:C2} and taking the supremum over all $1 \leq n \leq T/\Delta t - 1$, we deduce the first estimate in \eqref{result:h}
\begin{equation*}\sup_{2 \leq n \leq T/\Delta t}
	\left\|E^{n}
	-
	\sqrt{2}\,T_0\, \partial_h\,C^{n}_{2}
	\right\|_{l^2\left(\cJ\right)}
	\,\leq\,
	C\lambda
	\,,\quad \forall \lambda >0\,,
\end{equation*}
which concludes the proof.

\section{Proof of Proposition \ref{prop:reformulated_Poisson}}
\label{proof:prop:reformulated_Poisson}

Using~\eqref{discrete:step1} and \eqref{discrete:step2}, the equations for $C_0^{n+1}$, $C_1^{n+1}$ and $E^{n+1}$ can be re-written
\be\label{scheme:nplusun}
\left\{
\begin{array}{l}
	\ds\frac{C^{n+1}_{0} -C_{0}^{n} }{\Delta t}
	\,+\,\sqrt{T_0}\,\partial_h\,\Bigl(C^{n+1}_{1}-\Delta t\, \frac{1}{\sqrt{T_0}}\,E^{n+1}\, (C^{n+1}_{0}-1)\,\Bigl)=0\,,\\[1.2em]
	\ds\frac{C^{n+1}_{1} -C_{1}^{n} }{\Delta t}
	\,+\,\sqrt{T_0}\,\partial_h\,C^{n+1}_{0}
	\,+\,\sqrt{2\,T_0}\,\partial_h\,C^{(1,n+1)}_{2} \,-\, 
	\frac{1}{\sqrt{T_0}}\,E^{n+1}\,C^{n+1}_{0}=\,
	0\,,\\[1.2em]
	\ds\lambda^2\partial_hE^{n+1}\, =\, C^{n+1}_0-1\,,
\end{array}\right.
\ee
where $C^{(1,n+1)}_{2}$ corresponds to $C^{(1)}_{2}$ calculated during the Step~1 of the calculus of $(C_k^{n+1})_{k,\geq 0}$ knowing~$(C_k^{n})_{k\geq 0}$.

Then, inserting the first equation of~\eqref{scheme:nplusun} into the third one, we get
$$\displaylines{
	\lambda^2\,\frac{\partial_hE^{n+1}-2\,\partial_hE^{n}+\partial_hE^{n-1}}{\Delta t^2}\, =\,\frac{C^{n+1}_{0}-2\,C^{n}_{0}+C^{n-1}_{0}}{\Delta t^2}
	\hfill\cr
	\hfill =\,-\sqrt{T_0}\,\frac{\partial_h\,C^{n+1}_{1}-\partial_h\,C^{n}_{1}}{\Delta t}+\partial_h\,\Bigl(E^{n+1}\, (C^{n+1}_{0}-1)-E^{n}\,(C^{n}_{0}-1)\Bigl).}$$
Now, using the second and third equations of~\eqref{scheme:nplusun} yields
$$\displaylines{
	\lambda^2\,\frac{\partial_hE^{n+1}-2\,\partial_hE^{n}+\partial_hE^{n-1}}{\Delta t^2}\, \hfill\cr
	=\,\partial^2_h \Bigl(\sqrt{2}\,T_0\,C^{(1,n+1)}_{2}+T_0\,C^{n+1}_{0}\Bigl)-\partial_h\,\Bigl(E^{n+1}\,C^{n+1}_{0}\Bigl)+\lambda^2\,\partial_h\,\Bigl(E^{n+1}\, \partial_hE^{n+1}-E^{n}\, \partial_hE^{n}\Bigl).}$$
\bibliographystyle{abbrv}
\bibliography{refer}
\end{document}